\documentclass[a4paper,10pt]{amsart}
\usepackage{amssymb,amsmath,epsfig}
\usepackage[english]{babel}

\input xy
\xyoption{all}

\newtheorem{proposition}{Proposition}[section]
\newtheorem{theorem}[proposition]{Theorem}
\newtheorem{lemma}[proposition]{Lemma}
\newtheorem{definition}[proposition]{Definition}
\newtheorem{corollary}[proposition]{Corollary}

\newtheorem{example}[proposition]{Example}

\newcommand{\C}{\mathbf{C}}
\newcommand{\Z}{\mathbf{Z}}
\newcommand{\N}{\mathbf{N}}
\newcommand{\R}{\mathbf{R}}
\newcommand{\Q}{\mathbf{Q}}

\newcommand{\D}{\mathbf{D}}
\newcommand{\I}{\mathbf{I}}
\newcommand{\A}{\mathbf{A}}

\begin{document}

\title[Continued fractions and surface singularities]{The geometry
  of continued fractions and the topology of 
  surface singularities}  
\author{Patrick Popescu-Pampu}
\address{Univ. Paris 7 Denis Diderot, Inst. de
  Maths.-UMR CNRS 7586, {\'e}quipe "G{\'e}om{\'e}trie et dynamique" \\case
  7012, 2, place Jussieu, 75251 Paris cedex 05, France.}
\email{ppopescu@math.jussieu.fr}

\subjclass{Primary 52 C 05
    ; Secondary 14 M 25, 32 S 25, \linebreak 32 S 50, 57 N 10.}


\begin{abstract}{We survey the use of continued
fraction expansions in the algebraical and topological study of
complex analytic singularities. We 
also prove new results, firstly concerning a geometric duality
with respect to a lattice 
between plane supplementary cones and secondly concerning the existence of a
canonical 
plumbing structure on the \textit{abstract boundaries} (also called
\textit{links}) of 
normal surface singularities. The duality between
supplementary cones gives in particular a geometric interpretation of a 
duality discovered by Hirzebruch between the continued fraction
expansions of two numbers $\lambda >1$ and $\frac{\lambda}{\lambda -1}$. }
    
\end{abstract}

\maketitle
\thispagestyle{empty}

\tableofcontents

\section{Introduction}

Continued fraction expansions appear naturally when one resolves germs
of plane curves by sequences of plane blowing-ups, or Hirzebruch-Jung
(that is, cyclic quotient) surface singularities by toric
modifications. 

They also appear when one passes from the natural
plumbing decomposition of the abstract boundary of a normal surface
singularity to its minimal JSJ decomposition. In this case it is very
important to keep track of natural orientations. In general, as
was shown by Neumann \cite{N 81}, if one changes the orientation of
the boundary, the resulting 3-manifold is no more
orientation-preserving diffeomorphic to the boundary of an isolated
surface singularity. The only exceptions are {\em Hirzebruch-Jung
singularities} and {\em cusp-singularities}. This last class of singularities
got its name from its appearance in Hirzebruch's work \cite{H
  73} as germs at the compactified cusps of Hilbert modular
surfaces. For both classes of singularities, one gets an involution on
the set of analytical isomorphism types of the singularities in the
class, by changing the orientation of the boundary. From the viewpoint
of computations, Hirzebruch saw that both types of singularities have
structures which can be encoded in continued fraction expansions of
positive integers, and that the previous involution manifests itself in
a duality between such expansions. 

In the computations with continued fractions alluded to before, there appear
in fact two kinds of continued fraction expansions. Some are constructed
using only additions - we call them in the sequel {\em Euclidean continued
fractions} - and the others using only subtractions - we call them
{\em Hirzebruch-Jung continued fractions}. There is a simple formula, also
attributed to Hirzebruch, which allows to pass from one type of continued
fraction expansion of a number to the other one. Both types of
expansions have geometric interpretations in terms of polygonal lines
$P(\sigma)$. If $(L, \sigma)$ is a pair consisting of a 2-dimensional lattice
$L$ and a strictly convex cone $\sigma$ in the associated real vector
space, $P(\sigma)$ denotes the boundary of the convex hull of the set
of lattice points situated inside  $\sigma$ and different from the
origin. 

For Euclidean continued fractions this
interpretation is attributed to Klein \cite{K 96}, while for the
Hirzebruch-Jung 
ones it is attributed to Cohn \cite{C 73}. 

It is natural to try to understand how both geometric interpretations
fit together. By superimposing the corresponding drawings, we were led to
consider {\em two supplementary cones} in a real plane, in the presence of a
lattice. By \textit{supplementary cones} we mean two closed strictly
convex cones which have a common edge and whose union is a
half-plane. 

Playing with examples made us understand that the
algebraic duality between continued fractions alluded to before has as
geometric counterpart a duality between two supplementary cones in the
plane with respect to a lattice. This duality is easiest to express in
the case where the edges of the cones are irrational: 

\medskip
\textit{Suppose that the edges of the supplementary cones $\sigma$ and
  $\sigma'$ are
  irrational. Then the edges of each polygonal line $P(\sigma)$ and
  $P(\sigma')$ correspond bijectively in a 
  natural way to the vertices of the other one.}
\medskip

When at least one of the edges is rational, the correspondence is
slightly more complicated (there is a defect of bijectivity near the
intersection points of the polygonal lines with the edges of the
cones), as explained in Proposition \ref{parallel}. In this 
duality, points correspond to lines and conversely, as in the
classical polarity relation between points and lines with respect to a
conic. But the duality relation described in this paper is more
elementary, in the sense that it uses only parallel transport in the
plane. For this reason it can be explained very simply by drawing on a
piece of cross-ruled paper.

The duality between supplementary cones gives a simple way to
think about the 
relation between the pair $(L, \sigma)$ and its dual pair $(\check{L},
\check{\sigma})$ and in particular about the
relations between various invariants of toric surfaces (see section
\ref{toric}). Indeed (see Proposition \ref{isom}):

\medskip
  \textit{The supplementary cone of $\sigma$ is canonically isomorphic
  over the integers with the dual cone $\check{\sigma}$, once an
  orientation of $L$ is fixed}.
\medskip

As stated at the beginning of the introduction, computations with
continued fractions appear also when one passes from the canonical
plumbing structure on the boundary of a normal surface singularity to its
minimal JSJ structure. Using this, Neumann \cite{N 81} showed that the
topological type of the minimal good resolution of the germ is
determined by the topological type of the link. In fact all continued
fractions appearing in Neumann's work are the algebraic counterpart of
pairs $(L, \sigma)$ canonically determined by the topology of the
boundary. Using this remark, we prove the stronger statement (see
Theorem \ref{invplumb}):

\medskip
  \textit{The plumbing structure on the boundary of a normal surface
  singularity associated to the minimal normal crossings resolution is
  determined up to isotopy by the oriented ambient manifold. In
  particular, it is invariant up to isotopy under the group of
  orientation-preserving self-diffeomorphisms of the boundary.}
\medskip

In order to prove this theorem we have to treat separately the
boundaries of Hirzebruch-Jung and cusp singularities. In both cases,
we show that the oriented boundary determines naturally a pair $(L,
\sigma)$ as before. If one changes the orientation of the boundary,
one gets a supplementary cone. In this way, the involution defined
before on both sets of singularities is a manifestation of the
geometric duality between supplementary cones (see Propositions
\ref{orientchange} and \ref{dualcusp}). 
\medskip

For us, the moral of the story we tell in this paper is the following one:

\medskip
  \textit{If one meets computations with either Euclidean or
  Hirzebruch-Jung continued fractions in a geometrical problem, it
  means that somewhere behind is present a natural 2-dimensional
  lattice $L$ and a couple of lines in the associated real vector
  space. One has first to choose one of the
  two pairs of opposite cones determined by the four lines and
  secondly an ordering of the edges of those cones. These choices may
  be  dictated by choices of orientations of the manifolds which led to the
  construction of the lattice and the cones. So, in order to think
  geometrically at the computations with continued fractions, 
  recognize the lattice, the lines and the orientation choices}.
\medskip

Let us outline now the content of the paper. 

Someone who is interested only in the algebraic relations between the
Euclidean and the Hirzebruch-Jung continued fraction expansions of a
number can consult only section 2. If one is also interested in their
geometric interpretation, one can read sections 3 and 4. 

In section 5 we prove geometrically the relations
between the two kinds of continued fractions using the duality between
supplementary cones described before. We introduce also a new kind of
graphical representation which we call the \textit{zigzag diagram},
allowing to visualize at the same time the algebra and the geometry of
the continued fraction expansions of a number. 

In section 6 we give applications of zigzag diagrams to the algebraic
description of special curve and surface singularities, defined using
toric geometry. 

Sections 8 and 9 are dedicated to the study of topological aspects of
the links of normal surface singularities, after having recalled in
section 7 general facts about Seifert, graph, plumbing and
JSJ-structures on 3-manifolds. 

We think that the new results of the paper are Proposition
\ref{parallel}, Theorem \ref{canlens} and Theorem \ref{invplumb}, as
well as the very easy Proposition \ref{isom}, which is nevertheless
essential in order to understand the relation between dual cones in
terms of parallelism, using Proposition \ref{parallel}.

\medskip

We wrote this paper having in mind as a potential reader a graduate
student who wants to be initiated either to the algebra of surface
singularities or to their topology. That is why we tried to
communicate basic intuitions, often referring to the references for
complete proofs. 
\medskip

\textbf{Acknowledgments.} We are very grateful to Friedrich
Hirzebruch for the historical comments he sent us, as well as to Paolo
Lisca, Andras N{\'e}methi, Bernard Teissier and the anonymous referee for
their pertinent remarks and suggestions.

\medskip

\section{Algebraic comparison of Euclidean and Hirzebruch-Jung
  continued fractions} \label{algecomp}

\begin{definition} \label{continued}
If $x_1,...,x_n$ are variables, we consider two kinds of continued
fractions associated to them:

$$[x_1,...,x_n]^+:= x_1 + \cfrac{1}{x_2 +
                          \cfrac{1}{\cdots + \cfrac{1}{x_n}}}$$
$$[x_1,...,x_n]^-:= x_1 - \cfrac{1}{x_2 -
                          \cfrac{1}{\cdots - \cfrac{1}{x_n}}}$$

We call $[x_1,...,x_n]^+$ a \textbf{Euclidean continued fraction}
(abbreviated \textbf{E-conti\-nued fraction}) and $[x_1,...,x_n]^-$ a
\textbf{Hirzebruch-Jung continued fraction} (abbreviated
\textbf{HJ-continued fraction}). 
\end{definition}

The first name is motivated by the
fact that E-continued fractions are tightly related to the Euclidean
algorithm: if one applies this algorithm to a couple of positive
integers $(a,b)$ and the successive quotients are $q_1,...,q_n$, then
$a/b=[q_1,...,q_n]^+$.  See Hardy \& Wright \cite{HW 88}, Davenport \cite{D
  99} for an introduction to their arithmetics  and  Fowler 
\cite{F 87} for the relation with the Greek theories of
proportions. An extended bibliography on their applications can be found
in Brezinski \cite{B   91} and Shallit \cite{S 92}. 

The second name is motivated by the fact 
that HJ-continued fractions appear naturally in the Hirzebruch-Jung
method of resolution of singularities, originating in Jung \cite{J 08}
and Hirzebruch \cite{H 53}, as explained after Definition \ref{HJsing}
below.

Define two sequences $(Z^\pm (x_1,...,x_n))_{n \geq 1}$ of polynomials
with integer coefficients,  by the initial data
  $$ Z^\pm (\emptyset) =1, \: Z^\pm (x)= x$$
and the recurrence relations:
  \begin{equation} \label{recrel}
   Z^\pm (x_1,...,x_n)= x_1 Z^\pm (x_2,...,x_n) \pm Z^\pm
  (x_3,...,x_n), \: \forall \: n \geq 2.
  \end{equation}

Then one proves immediately by induction on $n$ the following equality
of rational fractions:
  \begin{equation} \label{exprfrac}
   [x_1,..., x_n]^\pm = \frac{Z^\pm (x_1,...,x_n)}{Z^\pm
    (x_2,...,x_n)}, \: \forall \: n \geq 1.
  \end{equation}

Also by induction on $n$, one proves the following twin of relation
(\ref{recrel}):

\begin{equation} \label{twin}
  Z^\pm (x_1,...,x_n)= Z^\pm (x_1,...,x_{n-1}) x_n \pm Z^\pm
  (x_1,...,x_{n-2}), \: \forall \: n \geq 2.
  \end{equation}
which, combined with (\ref{recrel}), proves the following symmetry
property:

\begin{equation} \label{sym}
  Z^\pm (x_1,...,x_n)= Z^\pm (x_n,...,x_1) , \: \forall \: n \geq 1.
  \end{equation}

If $(y_1,...,y_k)$ is a finite sequence of numbers or variables and $m \in
\mathbf{N}\cup \{ +\infty\}$, we denote by 
  $$(y_1,...,y_k)^m$$
the sequence obtained by repeating $m$ times the sequence
$(y_1,...,y_k)$. By convention, when $m=0$, the result is the
empty sequence.  

Each number $\lambda \in \R$ can be expanded as (possibly infinite)
Euclidean and Hirzebruch-Jung continued fractions:
 $$ \lambda= [a_1,a_2,...]^+ = [\alpha_1,\alpha_2,...]^-$$
with the conditions:
 \begin{equation} \label{restr1}
   a_1\in \Z, \: a_n \in \N-\{0\}, \: \forall \: n \geq 1
 \end{equation}

\begin{equation} \label{restr2}
   \alpha_1\in \Z, \: \alpha_n \in \N -\{0,1\}, \: \forall \: n \geq 1
 \end{equation}
Of course, we consider only indices $n$ effectively present. For an
infinite number of terms, 
these conditions ensure the existence of the limits $[a_1,
a_2,...]^+:= \displaystyle{\lim_{n\rightarrow
    +\infty}}[a_1,...,a_n]^+$ and $[\alpha_1, 
\alpha_2,...]^-:= \displaystyle{\lim_{n\rightarrow
    +\infty}}[\alpha_1,...,\alpha_n]^-$.

Any sequence $(a_n)_{n\geq 1}$ which verifies the restrictions
(\ref{restr1}) can appear  and the only ambiguity in the expansion of
a number as a E-continued fraction comes from the identity:
 \begin{equation} \label{id1}
   [a_1,...,a_n,1]^+= [a_1,...,a_{n-1}, a_n +1]^+
 \end{equation}

We deduce that any real number $\lambda \neq 1$ admits a unique expansion as a
E-continued fraction such that condition (\ref{restr1}) is satisfied
and in the case that the sequence $(a_n)_n$ is finite, its last term
is different from $1$. When we speak in the sequel about \textit{the
  E-continued fraction expansion} of a number $\lambda\neq 1$, it will
be about this one.  By analogy with the vocabulary of the Euclidean
algorithm, we say 
that the numbers $(a_n)_{n \geq 1}$ are the 
\textit{E-partial quotients} of $\lambda$.

Similarly, any sequence $(\alpha_n)_{n\geq 1}$ which verifies the
restrictions (\ref{restr2}) can appear and the only ambiguity in the
expansion of a number as a HJ-continued fraction comes from the identity:
 \begin{equation} \label{id2}
   [\alpha_1,...,\alpha_n,(2)^{\infty}]^-= [\alpha_1,...,\alpha_{n-1},
   \alpha_n -1]^- 
 \end{equation}

We see that any real number $\lambda$ admits a unique expansion as a
HJ-continued fraction such that condition (\ref{restr2}) is satisfied
and the sequence $(\alpha_n)_n$ is not infinite and ultimately
constant equal to $2$. When we speak in the sequel about \textit{the
  HJ-continued fraction expansion} of a number $\lambda$, it will
be about this one.  We call the numbers $(\alpha_n)_{n \geq 1}$ the
\textit{HJ-partial quotients} of $\lambda$.

\medskip

The following lemma (see Hirzebruch \cite[page 257]{H 73}) can be easily proved
by induction on the integer $b \geq 1$.  

\begin{lemma}
 If $a\in \Z,\: b \in \mathbf{N}-\{0\}$ and $x$ is a variable, then:
     $$[a,b,x]^+ = [a+1,(2)^{b-1},x+1]^-$$ 
\end{lemma}

Using this lemma one sees how to pass from the E-continued fraction
expansion of a real number $\lambda$ to its HJ-continued fraction
expansion:

\begin{proposition} \label{transf}
If $(a_n)_{n \geq 1}$ is a (finite or infinite)
  sequence of positive integers, then:  
  $$[a_1,...,a_{2n}]^+ = [a_1 +1, (2)^{a_2-1},a_3 +2, (2)^{a_4-1},..., 
                            (2)^{a_{2n}-1}]^-$$
  $$[a_1,...,a_{2n+1}]^+ = 
      [ a_1 +1, (2)^{a_2-1},a_3 +2, (2)^{a_4-1},...,
      (2)^{a_{2n}-1}, a_{2n+1} +1]^-$$ 
  $$[a_1, a_2, a_3, a_4,...]^+ =
      [ a_1 +1, (2)^{a_2-1},a_3 +2, (2)^{a_4-1}, a_5 +2, (2)^{a_6
      -1},...]^-$$ 
(recall that, by convention, $(2)^0$ denotes the empty sequence).
\end{proposition}

\begin{example} \label{extrans}
   $\frac{11}{7}= [(1)^3, 3]^+=   [2,3,(2)^2]^-$.
\end{example}

Notice that this procedure can be inverted. In particular, 
an immediate consequence of the previous proposition is that a number
has bounded E-partial quotients if and only if it has bounded
HJ-partial quotients. Similarly, it has ultimately periodic
E-continued fraction (which happens if and only if it is a quadratic
number, see Davenport \cite{D 99}) if and only if it has
ultimately periodic HJ-continued fraction. In this case, Proposition
\ref{transf} explains how to pass from its E-period to its HJ-period.

\medskip

The continued fraction expansions of two numbers which differ by an
integer are related in an evident and simple way. For this reason,
from now on we 
restrict our attention to real numbers $\lambda >1$. 
The map 
\begin{equation} \label{invol}
  \lambda \longrightarrow \frac{\lambda}{\lambda -1}
\end{equation}
is an involution of the interval $(1, +\infty)$ on itself. The
E-continued fraction expansions of the numbers in the same orbit of this
involution are related in a very simple way: 

\begin{lemma} \label{expos}
 If $\lambda \in (1, +\infty)$ and $\lambda=
 [a_1,a_2,...]^+$ is its expansion as a (finite or infinite)
 continued fraction, then:
  $$\frac{\lambda}{\lambda -1} = \left\{ 
        \begin{array}{ll}
            [1+a_2, a_3, a_4,...]^+, & \: \mbox{ if } \: a_1=1, \\ 
            
            [1,a_1 - 1, a_2, a_3,...]^+, &\: \mbox{ if } \: a_1\geq 2 
        \end{array} \right.$$
\end{lemma}

The proof is immediate, once one notices that $\frac{\lambda}{\lambda
  -1} = [1, \lambda -1]^+$. Notice also that the involutivity of the
map (\ref{invol}) shows that the first equality in the
previous lemma is equivalent to the second one.

\begin{example} \label{exinv} 
If $\lambda= \frac{11}{7}= [(1)^3, 3]^+$,
then $\frac{11}{4}= \frac{\lambda}{\lambda-1}= [2,1,3]^+$. 
\end{example}

By combining Proposition \ref{transf} and Lemma \ref{expos}, we get
the following relation between the HJ-continued fraction expansions of
the numbers in the same orbit of the involution (\ref{invol}):

\begin{proposition} \label{exneg}
  If $\lambda \in \mathbf{R}$ is greater than $1$ and 
     $$\lambda= [(2)^{m_1},n_1+3,(2)^{m_2}, n_2 +3,...]^-$$ 
is its expression as a
 (finite or infinite) continued fraction, with $m_i, n_i \in 
 \mathbf{N}, \: \forall \: i \geq 1$, then: 
    $$\frac{\lambda}{\lambda -1} = 
        [m_1+2, (2)^{n_1},m_2 +3,(2)^{n_2}, m_3+3,...]^-$$
\end{proposition}

For $\lambda$ rational, this was proved in a different way by Neumann
\cite[Lemma 7.2]{N 81}. It reads then:

$$ \begin{array}{ll} 
      \lambda= [(2)^{m_1},n_1+3,(2)^{m_2}, ..., n_{s} +3,(2)^{m_{s+1}}]^-
 \Longrightarrow \\ 
 \dfrac{\lambda}{\lambda -1} = 
        [m_1+2, (2)^{n_1},m_2 +3,..., (2)^{n_{s}}, m_{s+1}+2]^-
    \end{array}$$

The important point here is that even a value $m_{s+1}=0$ contributes to
the number of partial quotients in the HJ-continued fraction expansion of
$\frac{\lambda}{\lambda -1}$. 

The next proposition is equivalent to the previous one, as an easy
inspection shows. Its advantage is that it gives a graphical way to
pass from the HJ-continued fraction expansion of a number $\lambda >1$
to the analogous expansion of $\frac{\lambda}{\lambda -1}>1$. 

\begin{proposition}
 Consider a number $\lambda \in \mathbf{R}$ greater than $1$ and let  
     $$\lambda= [\alpha_1,\alpha_2,...]^-, \: \: 
      \frac{\lambda}{\lambda -1} = [\beta_1,\beta_2,...]^-$$
be  the expressions of $\lambda$ and $\frac{\lambda}{\lambda -1}$ as
(finite or infinite) HJ-continued fractions. Construct a diagram made
of points organized in lines and columns in the following way: 
    
 $\bullet$ its lines are numbered by the positive integers;
 
 $\bullet$ the line numbered $k \geq 1$ contains $\alpha_k -1$ points;
 
 $\bullet$ the first point in the line numbered $k+1$ is placed under the
 last point of the line numbered $k$.

Then the $k$-th column contains $\beta_k -1$ points.
\end{proposition}

This
graphical construction seems to have been first noticed by
Rie\-menschneider in \cite{R 74} when $\lambda \in \Q_+$. 
Nowadays one usually speaks about \textit{Riemenschneider's point
  diagram} or \textit{staircase diagram}.  

\medskip
\begin{example} \label{exdiag}
 If $\lambda= \frac{11}{7}=[2,3,(2)^2]^-$,
 the associated point diagram is:
 $$ \begin{array}{cc}
     \bullet & \\
       \bullet &  \bullet \\
          &  \bullet \\
          &  \bullet 
    \end{array}$$
 One deduces from it that $\frac{\lambda}{\lambda -1} =[3,4]^-$. 
\end{example}

\medskip

\section{Klein's geometric interpretation of Euclidean continued
  fractions} \label{Kleinint}

We let Klein \cite{K 04} himself speak about his interpretation, in order to
emphasize his poetical style:

\begin{quote}

{\small

Let us now enliven these considerations with geometric
pictures. Confining our attention to positive numbers, let us
\textit{mark all those points} in the positive quadrant of the $xy$
plane \textit{which have integral coordinates}, forming thus a
so-called \textit{point lattice}. Let us examine this lattice, I am
tempted to say this ``firmament" of points, with our point of view at
the origin. [...] Looking from
$0$, then, one sees points of the lattice \textit{in all rational
  directions and only in such directions}. The field of view is
everywhere ``densely" but not completely and continuously filled with
``stars". One might be inclined to compare this view with that of the
milky way. With the exception of $0$ itself there is \textit{not a
  single integral point lying upon an irrational ray
  $\frac{x}{y}=\omega$, where $\omega$ is irrational}, which is very
remarkable. If we recall Dedekind's definition of irrational number,
it becomes obvious that such a ray makes a \textit{cut in the field of
  integral points by separating the points into two point sets}, one
lying to the right of the ray and one to the left. If we inquire how
these point sets converge toward our ray $x/y=\omega$, we shall find a
very simple relation to the continued fraction for $\omega$. By
marking each point $(x=p_\nu, y=q_\nu)$, corresponding to the
convergent $p_\nu/q_\nu$, we see that the rays to these points
approximate to the ray $x/y=\omega$ better and better, alternately
from the left and from the right, just as the numbers $  p_\nu/q_\nu$
approximate to the number $\omega$. Moreover, if one makes use of the
known number-theoretic properties of $p_\nu, q_\nu$, one finds the
following theorem: \textit{Imagine pegs or needles affixed at all the
  integral points, and wrap a tightly drawn string about the sets of
  pegs to the right and to the left of the $\omega$-ray, then the
  vertices of the two convex string-polygons which bound our two point
  sets will be precisely the points $(p_\nu, q_\nu)$ whose coordinates
  are the numerators and denominators of the successive convergents to
  $\omega$, the left polygon having the even convergents, the right
  one the odd.} This gives a new and, one may well say, an extremely
graphical definition of a continued fraction.
}

\end{quote}

In the original article \cite{K 96}, one finds moreover the following
interpretation of the E-partial quotients:

\begin{quote}
{ \small
  Each edge of the polygons [...] may contain integral points. The
  number of parts in which the edge is decomposed by such points is
  exactly equal to a partial quotient.} 
\end{quote}

Before Klein, Smith expressed a related idea in \cite{S 76}:

\begin{quote}
  {\small If with a pair of rectangular axes in a plane we
  construct a system of unit points (\emph{i.e.} a system of points of
  which the coordinates are integral numbers), and draw the line $y=
  \theta x$, we learn from that theorem that if $(x,y)$ be a unit
  point lying nearer to that line than any other unit point having a
  less abscissa (or, which comes to the same thing, lying at a less
  distance from the origin), $\frac{y}{x}$ is a convergent to
  $\theta$; and, \emph{vice versa}, if $\frac{y}{x}$ is a convergent,
  $(x,y)$ is one of the `nearest points'. Thus the `nearest points'
  lie alternately on opposite sides of the line, and the double area
  of the triangle, formed by the origin and any two consecutive
  `nearest points', is unity.}
\end{quote}

Proofs of the preceding properties can be found in
Stark \cite{S 78}. Here we only sketch the reason of Klein's
interpretation. For explanations about our vocabulary, read next
section. 

Let $\lambda>1$ be a real number. In the first quadrant $\sigma_0$,
consider the half-line $L_\lambda$ of slope $\lambda$ (see Figure
1). It is defined 
by the equation $y=\lambda x$, which shows that $\lambda= \omega
^{-1}= \theta$, where $\omega$ is Klein's notation and $\theta$ is
Smith's. It subdivides the quadrant $\sigma_0$ into 
two closed cones with vertex the origin, $\sigma_x(\lambda)$ adjacent to the
axis of the variable $x$ and $\sigma_y(\lambda)$ adjacent to the
axis of the variable $y$. 

\begin{lemma} \label{firstedge}
  The segment which joins the lattice points of coordinates $(1,0)$
  and $(1, a_1)$ is a compact edge of the convex hull of the set of lattice
  points different from the origin contained in the cone
  $\sigma_x(\lambda)$, where $\lambda=[a_1,a_2,...]^+$ is the
  E-continued fraction expansion of $\lambda$. 
\end{lemma}

\textbf{Proof:} Indeed, the half-line starting from $(1,0)$
  and directed towards $(1, a_1)$ cuts the half-line $L_\lambda$
  inside the segment $[(1, [\lambda]), \: (1, [\lambda]+1))$, where
  $[\lambda]$ is the integral part of $\lambda$. But $[\lambda]=a_1$,
  which finishes the proof. 
  \hfill $\Box$
\medskip

{\tt    \setlength{\unitlength}{0.92pt}}
\begin{figure} \label{Kleininit}
   \epsfig{file=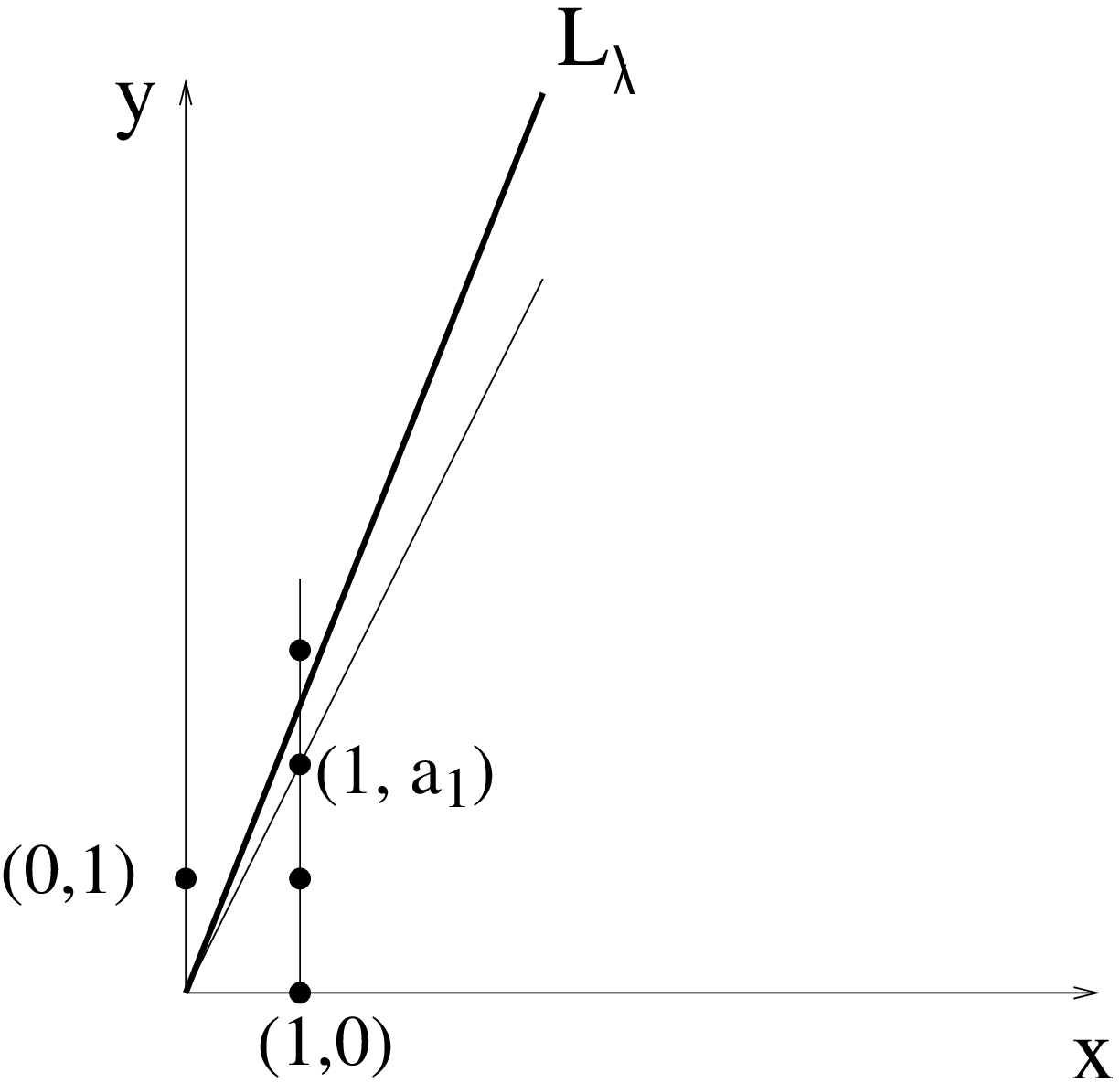, height= 50 mm}
   \caption{Figure illustrating the proof of Lemma \ref{firstedge}}
\end{figure}

Replace now the initial basis of the lattice  by $(0,1), \: (1,
a_1)$. With respect to this new basis, the 
slope of the half-line $L_\lambda$ is $(\lambda -a_1)^{-1}=
[a_2,a_3,...]^+$. This allows one to prove Klein's interpretation by
induction.
\medskip

If one considers all lattice points on the compact edges of the
boundaries of the two previous convex hulls instead of only the
vertices, and then one looks at the slopes of the lines which join
them to the origin, one obtains the so-called \textit{slow approximating
  sequence} of $\lambda$. This kind of sequence appears naturally when
one desingularizes germs of complex analytic plane curves by
successively blowing up points (see Enriques \&
Chisini \cite{EC 18}, Michel \& Weber \cite{MW 85} and L{\^e}, Michel \&
Weber \cite{LMW 89}). We leave as an exercise for the interested
reader to interpret this geometrically (first, read Section
\ref{moncurves}).

As explained by Klein himself in \cite{K 96}, his interpretation
suggests to generalize 
the notion of continued fraction to higher dimensions by taking the
boundaries of convex hulls of lattice points situated inside convex
cones. For references about recent research in this area, see Arnold
\cite{A 98} and Moussafir \cite{M 00}.

\section{Cohn's geometric interpretation of Hirzebruch-Jung
  continued fractions} \label{geomHJ} 

A geometric interpretation of HJ-continued fractions analogous to
Klein's interpretation of Euclidean ones was given by Cohn \cite{C
  73} (see the comment on his work in Hirzebruch \cite[2.3]{H 73}). It
seems to have soon become folklore among people doing toric 
geometry (see section \ref{toric}). Before describing this
interpretation, let us introduce some vocabulary in order to speak
with more precision about convex hulls of lattice points in the
plane. 

Let $L$ be a lattice of rank $2$, that is, a free abelian group of
rank $2$. It embeds canonically into the associated real vector space 
$L_{\R}= L \otimes_{\Z} \R$. When we picture the elements of $L$ as
points in the affine plane $L_\R$, we call them the \textit{integral points}
of the plane. When $A$ and $B$ are points of the affine plane $L_\R$,
we denote by $AB$ the element of the vector space $L_\R$ which
translates $A$ into $B$, by $[AB]$ the closed segment in $L_\R$ of
extremities $A,B$ and by $[AB$ the closed  half-line having $A$ as an
extremity  and directed towards $B$. 

If $(u,v)$ is an ordered basis of $L_{\R}$ and $l$ is a line of
$L_{\R}$, its \textit{slope} is the quotient $\beta/ \alpha \in \R \cup
\{\infty\}$, where $\alpha u + \beta v$ generates $l$. 

\begin{definition} \label{elem}
A (closed convex) triangle $ABC$ in $L_\R$ is called \textbf{elementary} if
its vertices are integral and they are the only intersections of the
triangle with the lattice $L$. 
\end{definition}

If the triangle $ABC$ is elementary, then each pair of vectors $(AB,
AC)$, $(BC, BA)$, $(CA, CB)$ is a basis of the lattice
$L$. Conversely, if one of 
these pairs is a basis of the lattice, then the triangle is
elementary. 

We call a line or a half-line in $L_\R$
\textit{rational} if it contains at least two integral points. If so,
then it contains an infinity of them. If $A$ and $B$ are two integral
points, the \textit{integral length} $l_\Z[AB]$ of the segment $[AB]$ is the
number of subsegments in which it is divided by the integral points it
contains. A vector $OA$ of $L$ is called \textit{primitive} if $l_{\Z}[OA]=1$.

Let $\sigma$ be a closed strictly convex 2-dimensional cone in the plane
$L_\R$, that is, the convex ``angle" (in the language of plane
elementary geometry) delimited by two non-opposing half-lines
originating from $0$. These half-lines are called the \textit{edges}
of $\sigma$. The cone $\sigma$ is called \textit{rational} if its
edges are rational. A cone is called \textit{regular} if its
edges contain points $A,B$ such that the triangle $OAB$ is
elementary. The name is motivated by the fact that the associated
toric surface $\mathcal{Z}(L, \sigma)$ is smooth (that is, all its
local rings are \textit{regular}) if and only if $\sigma$ is regular
(see section \ref{elemtor}).

Let $K(\sigma)$ be the convex hull of the set of 
lattice points situated inside $\sigma$, with the exception of the
origin, that is:
  $$K(\sigma):= \mathrm{Conv}(\sigma \cap (L - \{0\})).$$
The closed convex set $K(\sigma)$ is unbounded. 
Denote by $P(\sigma)$ its boundary: it is a connected polygonal line. It 
has two ends (in the topological sense), each one being asymptotic to
(or contained inside) an edge of 
$\sigma$ (see Figure 2).  An edge of $\sigma$ intersects  $P(\sigma)$
if and only if it is rational. 

Denote by $\mathcal{V}(\sigma)$ the set of vertices of $P(\sigma)$ and by
$\mathcal{E}(\sigma)$ the set of its (closed) edges.   For example, in
Figure 3 the vertices are the points $A_0, A_2, A_5$ and the edges are
the segments $[A_0 A_2], [A_2 A_5]$ and two half-lines contained in
$l_-, l_+$, starting from $A_0$, respectively $A_5$.  

Now  \textit{order} arbitrarily the edges of $\sigma$. Denote by $l_-$ the
first one and by $l_+$ the second one. This orients the plane $L_\R$,
by deciding to turn from $l_-$ towards $l_+$ \textit{inside}
$\sigma$. 
If we orient $P(\sigma)$ from the end which is asymptotic to $l_-$ towards the
end which is asymptotic to $l_+$, we get induced orientations of its
edges.

{\tt    \setlength{\unitlength}{0.92pt}}
\begin{figure} \label{Convex}
   \epsfig{file=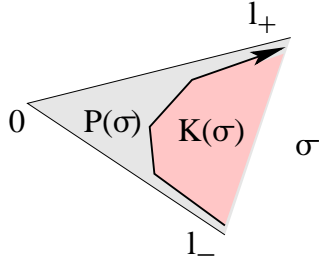, height= 35 mm}
   \caption{The polygonal line associated to a convex cone }
\end{figure}
\medskip

\textit{Suppose now that
  the edge $l_-$ of $\sigma$ is rational}. Denote then by $A_-\neq 0$
  the integral point of the half-line $l_-$ 
which lies nearest to $0$, and by $V_-\neq A_-$ the vertex of $P(\sigma)$
which lies nearest to $A_-$. Define in the same way $A_+$ and $V_+$
whenever $l_+$ is rational.  
Denote  by $(A_n)_{n \geq 0}$ the sequence of
integral points on 
$P(\sigma)$, enumerated as they appear when one travels on this
polygonal line in the positive direction, starting from $A_0=A_-$. If
$l_+$ is a rational half-line, then we stop this sequence when we
arrive at the point $A_+$. If $l_+$ is irrational, then this sequence
is infinite. Define $r \geq 0$ such that $A_{r+1}=A_+$. So,
$r=+\infty$ if and only if $l_+$ is irrational. 

\begin{example} \label{exallpoints} 
 We consider the lattice $\Z^2 \subset \R^2$
and the cone $\sigma$ with rational edges, generated by the vectors $(1,0)$ and
$(4,11)$ (see Figure 3). The small dots represent 
integral points in the plane and the bigger ones represent integral
points on the polygonal lines $P(\sigma)$.  In this example we have
$V_+=V_-=A_2$.  
\end{example}

{\tt    \setlength{\unitlength}{0.92pt}}
\begin{figure} \label{Ex411}
  \epsfig{file=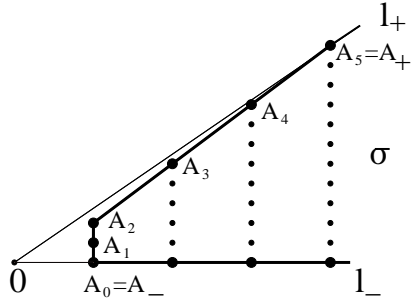, height= 40 mm}
  \caption{An illustration of Example \ref{exallpoints} }
\end{figure}

Each triangle $OA_{n}A_{n+1}$ is elementary, by the construction of the
convex hull $K(\sigma)$, which implies that all the couples $(OA_n, 
OA_{n+1})$ are bases of $L$. This shows that for any $n\in
\{1,...,r\}$, one has a relation of the type:
  \begin{equation} \label{eq3}
     OA_{n+1} + OA_{n-1} = \alpha_n \cdot OA_n 
  \end{equation}
with $\alpha_n \in \Z$, and the convexity of $K(\sigma)$ shows that:
  \begin{equation} \label{ineq}
      \alpha_n \geq 2, \: \forall \: n \in \{1,...,r\}
  \end{equation}

Conversely: 

\begin{proposition} \label{concrete}
 Suppose that $(OA_n)_{n \geq 0}$ is a (finite or infinite)
  sequence of   primitive vectors of $L$, related by relations of the form
  (\ref{eq3}). Then  we have
    $$OA_n=Z^-(\alpha_1,...,\alpha_{n-1})OA_1 - 
    Z^-(\alpha_2,...,\alpha_{n-1})OA_0, \: \forall \: n\geq 1$$ and
    the slope of the half-line 
    $l_+= \lim_{n \rightarrow \infty}[OA_n)$ in the base 
    $(-OA_0, OA_1)$ is equal to $[\alpha_1, \alpha_2,...]^-$. 
\end{proposition} 

\textbf{Proof:} Recall that the polynomials $Z^-$ were defined by the
recursion formula (\ref{recrel}). The first assertion can be easily
proved by induction, 
using the relations (\ref{eq3}). The second one is a consequence of
formula (\ref{exprfrac}), which shows that the slope of the half-line
$[OA_n$ in the base $(-OA_0, OA_1)$ is equal to $[\alpha_1,...,
\alpha_{n-1}]^-$. \hfill $\Box$

\begin{proposition} \label{converse}
  Let $\sigma$ be the closure of the convex hull of the
  union of the half-lines $([OA_n)_{n \geq 0}$. Then $\sigma$ is
  strictly convex and the points $\{A_n\}_{n \geq 1}$ are precisely
  the integral 
  points on the compact edges of the polygonal line $P(\sigma)$ if and
  only if the conditions (\ref{ineq}) are satisfied and the sequence
  $(\alpha_n)_{n \geq 1}$ is not infinite and ultimately constant
  equal to $2$.
\end{proposition}

\textbf{Proof:} $\bullet$ What remains to be proved about the
\textit{necessity} is that if the sequence $ (\alpha_n)_{n \geq 1}$ is
infinite, then it cannot be ultimately constant equal to $2$. If this
was the case, by relation (\ref{id2}) we would deduce that
$[\alpha_1,\alpha_2,...]^-$ is rational, and Proposition
\ref{concrete} would imply that $l_+$ is rational. Then $P(\sigma)$
would contain a finite number of integral points on its compact edges,
which would contradict the infinity of the sequence $(\alpha_n)_{n
  \geq 1}$.

$\bullet$ Let us prove now the \textit{sufficiency}. As $\alpha_n\geq
2$, $\forall \: n \in \{1,...,r\}$, we see that the triangles
$(OA_n A_{n+1})_{n \geq 0}$ turn in the same sense. Moreover,
Proposition \ref{concrete} shows that $\sigma$ is a strictly convex
cone. The vertices of the polygonal line $P=A_0A_1A_2...$ are
precisely those points $A_n$ for which $\alpha_n \geq 3$. As all the
triangles $OA_nA_{n+1}$ are elementary, we see that the origin $O$ is
the only integral point of the connected component of $\sigma-P$ which
contains it. Moreover, conditions  (\ref{ineq}) show that the other
component is convex. So, $P\subset P(\sigma)$. 

The proposition is
proved. \hfill $\Box$

\medskip

\section{Geometric comparison of Euclidean and HJ-continued fractions}
\label{geomcomp}  

In  section \ref{comp1} we relate the two preceding interpretations, by
describing  a duality between two supplementary cones in the plane, an
underlying lattice being fixed (see Proposition \ref{parallel}). In
section \ref{zigdiag} we introduce a so-called \textit{zigzag diagram} based
on this duality, which makes it very easy to visualize the
various relations between continued fractions proved algebraically in
section \ref{algecomp}. In section \ref{reldual} we give a proof of
the isomorphism between the supplementary cone $(L, \sigma')$ and the
dual cone $(\check{L}, \check{\sigma})$ of a given cone $(L, \sigma)$.
\medskip

\subsection{A geometric duality between supplementary cones} \label{comp1}
$\:$
\medskip

Suppose again that $\sigma$ is \textit{any} strictly convex
cone in $L_\R$, whose edge $l_-$ is not necessarily rational. 
Let $l_-'$ be the half-line opposite to $l_-$ and $\sigma'$
be the closed convex cone bounded by $l_+$ and $l_-'$. So, $\sigma$
and $\sigma'$ are \textit{supplementary cones}:

\begin{definition} \label{suppl}
  Two strictly convex cones in a real plane are called
  \textbf{supplementary} if they have a
common edge and if their union is a half-plane.
\end{definition}

By analogy with
what we did in the previous section for $\sigma$, orient the polygonal
line $P(\sigma')$ from 
$l_-'$ towards $l_+$. If $l_-$ is rational, define the point $A_-'$
and the sequence $(A_n')_{n\geq 0}$, with $A_0'=A_-'$. They are the
analogs for $\sigma'$ of the points $A_-$ and $(A_n)_{n \geq 0}$ for
$\sigma$. In particular, $OA_- + OA_-'=0$. 

\begin{example} \label{exsym}
 Consider the same cone as in Example \ref{exallpoints}. Then the
 polygonal lines $P(\sigma)$ and $P(\sigma')$ are  represented in
 Figure 4 using heavy segments. 
\end{example}

{\tt    \setlength{\unitlength}{0.92pt}}
\begin{figure}
  \epsfig{file=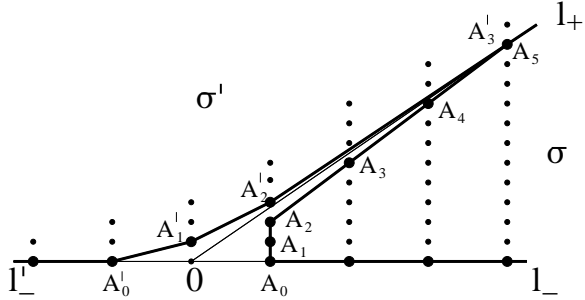, height= 40 mm}
  \caption{An illustration of Example \ref{exsym}}

\end{figure}

\medskip
The basis for our
geometric comparison of Euclidean and Hirzebruch-Jung continued
fractions is the observation that {\em the polygonal line $P(\sigma')$ can
be constructed in a very simple way once one knows
$P(\sigma)$}. Namely, starting from the origin, one draws the
half-lines parallel to the oriented edges of $P(\sigma)$. On each
half-line, one considers the integer point which is nearest to the
origin. Then the polygonal line which joins those points is the union
of the compact edges of $P(\sigma')$.

\medskip

Now we describe this with more precision. If
$e\in \mathcal{E}(\sigma)$ is an edge of $P(\sigma)$, denote by
$\mathcal{I}(e)\in L$ the integral point  
 such that $O\mathcal{I}(e)$ is a primitive vector of $L$ positively
parallel to $e$ (where $e$ is oriented according to the chosen
orientation of $P(\sigma)$). Then it is an easy exercise to see that 
 $\mathcal{I}(e) \in \sigma'$ (use the fact that the line containing $e$
 intersects $l_-$ and $l_+$ in interior points).
We can define a map:
  \begin{equation} \label{mapi}
  \begin{array}{llcl}
     \mathcal{I}: & \mathcal{E}(\sigma) & \longrightarrow & \sigma'\cap L\\
        & e & \longrightarrow & \mathcal{I}(e)
    \end{array}
  \end{equation}
As the edges of $P(\sigma)$ always turn in the same direction, one
sees that the map $\mathcal{I}$ is injective.

\begin{proposition} \label{parallel}
  The map $\mathcal{I}$ respects the orientations and the image of
  $\mathcal{I}$ verifies the double inclusion  
   $$\mathcal{V}(\sigma') \subset \mathrm{Im} (\mathcal{I}) \subset
   P(\sigma')\cap L.$$  
The difference $\mathrm{Im} (\mathcal{I})- \mathcal{V}(\sigma')$
   contains at most the 
points  $\mathcal{I}[A_- V_-]$, $\mathcal{I}[V_+A_+]$. Such a point is
   a vertex of 
$P(\sigma')$ if and only if the integral length of the corresponding
edge of $P(\sigma)$ is $\geq 2$.  
In particular,  one has the equality $\mathcal{V}(\sigma') =
\mathrm{Im} (\mathcal{I})$ if 
and only if $l_\Z[A_- V_-]\geq 2$ and $l_\Z[V_+A_+]\geq 2$, whenever
these segments exist. 
\end{proposition}

\textbf{Proof: }
 Denote by $(V_j)_{j \in J}$ the vertices of $P(\sigma)$, enumerated
 in the positive direction. The 
 indices form a set of consecutive integers, well-defined only up to
 translations. 

For any $j \in J$, denote by $V_j^-$ and $V_j^+$ respectively the
integral points of  
$P(\sigma)$ which precede and follow $V_j$. If $V_j$ is
an interior point of $\sigma$, denote by
$W_j\in L$ the point such that $OW_j = OV_j^- +OV_j^+$, and by
$W_j^-$ its nearest integral point in the interior of the segment
$[OW_j]$ (see Figure 5). 

{\tt    \setlength{\unitlength}{0.92pt}}
\begin{figure}
   \epsfig{file=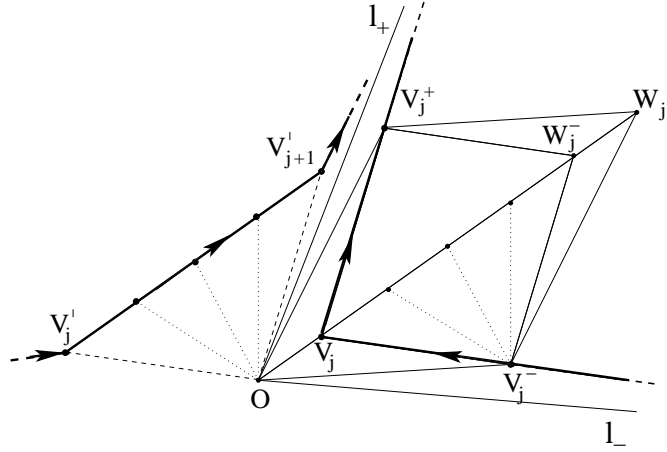, height= 60 mm}
   \caption{The first illustration for the proof of Proposition
     \ref{parallel}} 
\end{figure}

As $OV_j ^- V_j$ and $OV_jV_j^+$ are elementary triangles, it implies
that both $(OV_j^- ,OV_j)$ and $(OV_j, OV_j^+)$ are bases of $L$. So,
there exists an integer $n_j$  such that 
  \begin{equation} \label{defweight}
    OV_j^- + OV_j^+ = (n_j+3)OV_j.
  \end{equation} 
As $V_j$ is a vertex of $P(\sigma)$, we see that $n_j \geq 0$. 
We deduce that the points $O, V_j, W_j^-, W_j$ are aligned in this
order, that $V_jV_j^- + V_j V_j^+= V_j W_j^-$ and that
$l_{\Z}[V_jW_j^-]= n_j +1$.  

Let us join each one of the $n_j$ interior
points of $[V_jW_j^-]$ to $V_j^-$. This gives a decomposition of 
the triangle $V_j^-V_jW_j^-$  into $(n_j +1)$
triangles. These are necessarily elementary, because the triangle 
$OV_j^-V_j$ is. Denote 
  $$V_j'= \mathcal{I}[V_{j-1}V_j] \:\:\mathrm{and}\:\:
  V_{j+1}'=\mathcal{I}[V_jV_{j+1}].$$  
By the definition of the map $\mathcal{I}$, we see that $OV_j'= V_j^-V_j$  and
$OV_{j+1}' =V_jV_j^+ = V_j^-W_j^-.$
 This implies that \textit{the triangle $OV_j'V_{j+1}'$ is the translated image
by the vector  $V_j^-O$ of the triangle $V_j^-V_jW_j^-$}. The preceding
     arguments show that its only integral points are its vertices and
     $n_j$ other points in the interior of the segment
     $[V_j'V_{j+1}']$. Indeed:
\begin{equation} \label{proport} 
      V_j'V_{j+1}'= V_jW_j^- =(n_j+1)OV_j
   \end{equation}
Moreover, the triangle  $OV_j'V_{j+1}'$is included in the cone
   $\sigma'$ and the couple of 
   vectors $(OV_j', OV_{j+1}')$ has the same orientation as $(l_-',l_+)$. 

This shows that the triangles $(OV_j'V_{j+1}')_{j \in J}$ are pairwise
   disjoint and that their union does not contain integral points in
   its interior. 

$\bullet$ \textit{If both edges of $\sigma$ are irrational}, then the
   closure of 
   the union of the cones $\R_+ OV_j' + \R_+ OV_{j+1}'$ is the cone
   $\sigma'$, as the edges $l_-$ and $l_+$ are asymptotic to
   $P(\sigma)$.  We deduce from relation (\ref{proport}) that the
   sequence $(\lambda_j)_{j \in J}$  
of slopes of the vectors $(V_j'V_{j+1}')_{j \in J}$, expressed in a
base $(u_-, u_+)$ of $L_\R$ which verifies $l_\pm =\R_+ u_\pm$ is
\textit{strictly} increasing, and that $\lim_{j \rightarrow -\infty}
\lambda_j =0$, $\lim_{j \rightarrow +\infty}\lambda_j =+\infty$. This
shows that the closure of the connected component of $\sigma'-
\bigcup_{j \in J} [V_j'V_{j+1}']$ which does not contain the origin is
\textit{convex}. As a consequence,
  $$\bigcup_{j \in J} [V_j'V_{j+1}']= P(\sigma').$$
Moreover, as $n_j \geq 0$, the strict monotonicity of the sequence
$(\lambda_j)_{j \in J}$ implies that the
points $(V_j')_{j \in J}$ are precisely the vertices of
$P(\sigma')$. The proposition is proved in this case.

$\bullet$ \textit{Suppose now that $l_-$ is rational.} Then choose the
index set 
$J$ such that $V_0= A_-$ and $V_1= V_-$. By the construction of the
map $\mathcal{I}$, \textit{the triangle $OV_0'V_1'$ is the translated image of
$V_0OV_0^+$ by the vector $V_0 O$} (see Figure 6).      

{\tt    \setlength{\unitlength}{0.92pt}}
\begin{figure}
   \epsfig{file=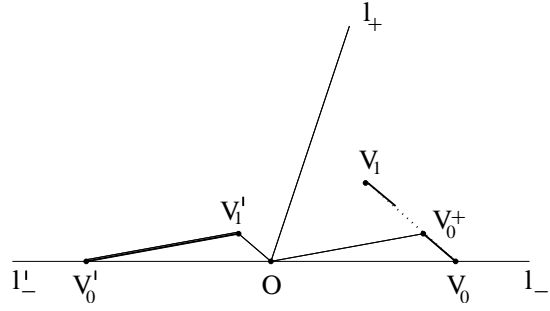, height= 40 mm}
   \caption{The second illustration for the proof of Proposition
     \ref{parallel}} 
\end{figure}

In particular, $V_0'V_1'=OV_0^+$. But $V_1'V_2'=(n_1 +1)OV_1$ by  relation
(\ref{proport}), which shows that the vectors $V_0'V_1'$ and
$V_1'V_2'$ are proportional if and only if $V_0^+ =V_1$, which is
equivalent to $l_\Z[A_-V_-] =1$. Moreover, the property of
monotonicity for the slopes of the vectors $(V_j'V_{j+1}')_{j \in J}$
is true as before, if one starts from $j =0$.

$\bullet$ An analogous reasoning is valid for $l_+$ \textit{if this
  edge of $\sigma$ is 
rational}. By combining all this, the proposition is proved.  
\hfill $\Box$

\medskip

The previous proposition explains a geometric \textit{duality} between the
supplementary cones $\sigma, \sigma'$ with respect to the lattice
$L$. We see that, with possible exceptions for the compact edges which
intersect the edges of $\sigma$ and $\sigma'$, the
compact edges of $P(\sigma)$ correspond to the vertices of
$P(\sigma')$ interior to $\sigma'$ and conversely (by permuting the
roles of $\sigma$ and $\sigma'$), which is a kind of
point-line polarity relation. 

The next corollary shows that the involution (\ref{invol}) studied
algebraically in section 
\ref{algecomp} is closely related to the previous duality.

\medskip

\begin{corollary} \label{slopes}
Suppose that $l_-$ is rational and that $\sigma$ is not regular. If
$(OA_0', U)$ 
is a basis of $L$ with respect to which the slope of $l_+$ is greater
than $1$, then $U= OA_1$. 
   If $\lambda >1$ denotes the slope of the half-line  $l_+$ in the
   base $(OA_0', OA_1)$, then $\frac{\lambda}{\lambda -1}$ is its slope
   in the base $(OA_0, OA_1')$.
\end{corollary}

\textbf{Proof:} 
We leave the first affirmation to the reader (look at
Figure 6). 

As the triangles $OA_0A_1$ and $OA_0'A_1'$ are elementary, we see that $(OA_0,
OA_1)$ and $(OA_0', OA_1')$ are indeed two bases of the lattice $L$. 
Proposition \ref{parallel} shows that $OA_0'=A_0A_1$, which allows us
to relate the two bases:

\begin{equation} \label{change}
 \left\{ \begin{array}{l}
            OA_0'=-OA_0\\
            OA_1'= OA_1 -OA_0
          \end{array}   \right.
\end{equation}

Let $v\in L_\R$ be a vector which generates the half-line $l_+$. We
want to express it in these two bases. As $l_+$ lies between the
half-lines $[OA_0'$ and $[OA_1$, we see that:
\begin{equation} \label{eq1}
   v = -q \:OA_0 + p \: OA_1, \: \mbox{ with } p, q \in \R_+^*
\end{equation}
The equations (\ref{change}) imply then that:
\begin{equation} \label{eq2}
   v = -(p-q) \: OA_0' + p \: OA_1'
\end{equation}
which shows that $p-q >0$, as $l_+$ lies between the half-lines
$[OA_1'$ and $[OA_0$. This implies that $\lambda:= \frac{p}{q}>1$. We
then deduce the corollary from equation (\ref{eq2}). 
  \hfill $\Box$
\medskip

 The previous corollary shows that the number $\lambda >1$ can be
 canonically attached to the pair $(L, \sigma)$, once a rational edge
 of $\sigma$ is chosen as the first edge $l_-$. This motivates the
 following definition: 

\begin{definition} \label{typecone}
  Suppose that $l_-$ is rational and that the cone $\sigma$ is not
  regular. 
  We say that the pair $(L, \sigma)$ with the chosen ordering of
  sides is \textbf{of type } $\lambda >1$ if $\lambda$ is the slope of
  the half-line $l_+$ in the base $(OA_0', OA_1)$. 
\end{definition}

Proposition \ref{concrete} shows that, if $(L, \sigma)$ is of type
$\lambda >1$, then $\lambda=[\alpha_1, \alpha_2,...]^-$, where the
sequence $(\alpha_n)_{n \geq 1}$ was defined using relation (\ref{eq3}).

Suppose now that \textit{both edges of $\sigma$ are rational}. Then one can
choose $p,q \in \N^*$ with $gcd(p,q)=1$ in relation (\ref{eq1}),
condition which determines them uniquely. So,
$\lambda=\frac{p}{q}$. The following proposition describes the type of
$(L, \sigma)$ after changing the ordering of the sides. 

\begin{proposition} \label{changeorder}
  If $(L, \sigma)$ is of type $\dfrac{p}{q}$ with respect to the
  ordering $l_-, l_+$, then it is of
  type $\dfrac{p}{\overline{q}}$ with 
  respect to the ordering $l_+, l_-$, where 
  $q\overline{q}\equiv 1  (\mbox{mod } p)$. 
\end{proposition}

\textbf{Proof:} By relation (\ref{eq1}), we have $OA_+= -q OA_- + p
OA_1$. Multiply both sides by $\overline{q}$. By the definition of
$\overline{q}$, there exists $k \in \N$ such that $q\overline{q}= 1 +
kp$. We deduce that $OA_-= -\overline{q}OA_+ + p(\overline{q} OA_1 -k
OA_-)$. So, $(-OA_+, \overline{q} OA_1 -k OA_-)$ is a base of $L$ in
which the slope of $l_-$ is $\frac{p}{\overline{q}} >1$. By the first
affirmation of Corollary \ref{slopes}, the proposition is proved.
 \hfill $\Box$
\medskip

By combining the previous proposition with Proposition \ref{concrete},
we deduce the following classical fact (see \cite[section III.5]{BHPV 04}):

\begin{corollary}
  If $\dfrac{p}{q}= [\alpha_1, \alpha_2,..., \alpha _r]^-$, then
  $\dfrac{p}{\overline{q}} = [\alpha_r,\alpha_{r-1},..., \alpha_1]^-$. 
\end{corollary}

Another immediate consequence of Corollary \ref{slopes} is:

\begin{proposition} \label{dualtype}
  If $(L, \sigma)$ is of type $\dfrac{p}{q}$ with respect to the
  ordering $l_-, l_+$, then $(L, \sigma')$ is of
  type $\dfrac{p}{p-q}$ with 
  respect to the ordering $l_-', l_+$.
\end{proposition}

The previous proposition describes the relation between the types of
two supplementary cones. In section \ref{zigdiag}, we describe more
precisely the relation between numerical invariants attached to the
edges and the vertices of $P(\sigma)$ and $P(\sigma')$. 

\medskip

\subsection{A diagram relating Euclidean and HJ-continued fractions}
\label{zigdiag} 
$\: $
\medskip

We introduce now a diagram which allows one to ``see" the duality between
$P(\sigma)$ and $P(\sigma')$, as well as the relations between the
various  numerical 
invariants attached to these polygonal lines. 

$\bullet$ \textit{Suppose first that both $l_-$ and $l_+$ are
  irrational.} Consider two consecutive vertices $V_j, V_{j+1}$ of
$P(\sigma)$.  Let us {\em attach the weight $n_j +3$ to the vertex}
$V_j$, where $n_j \geq 0$ was defined by relation
  (\ref{defweight}). Introduce also the integer $m_{j+1}\geq 0$ such
  that 
$l_\Z[V_jV_{j+1}]=m_{j+1} +1$. The relation (\ref{proport}) shows that
$l_\Z[V_j'V_{j+1}']=n_j +1$. By reversing the roles of the polygonal
lines $P(\sigma')$ and $P(\sigma)$, we deduce that the weight of 
the vertex $V_{j+1}'$ of $P(\sigma')$ is $m_{j+1}
+3$. 

We can visualize the relations between the vertices $V_j,
V_{j+1}, V_j', V_{j+1}'$ as well as the numbers associated to them and
to the segments $ [V_jV_{j+1}]$, $ [V_j'V_{j+1}']$ by using a 
diagram, in which the heavy lines represent the polygonal lines $P(\sigma),
P(\sigma')$, and each vertex $V_j$ is joined to $V_j'$ and
  $V_{j+1}'$ (see Figure 7). In this way, the region contained between
  the two curves 
  representing $P(\sigma)$ and $P(\sigma')$ is subdivided into
  triangles. Each edge $E$ of $P(\sigma)$, $P(\sigma')$ is contained in
  only one of those triangles. Look at its opposite vertex. We 
  say that $E$ is the \textit{opposite edge} of that vertex in the
  zigzag diagram. We see that
  \textit{the weight of a vertex is equal to the length 
  of the opposite edge augmented by $2$}.

{\tt    \setlength{\unitlength}{0.92pt}}
\begin{figure}
  \epsfig{file=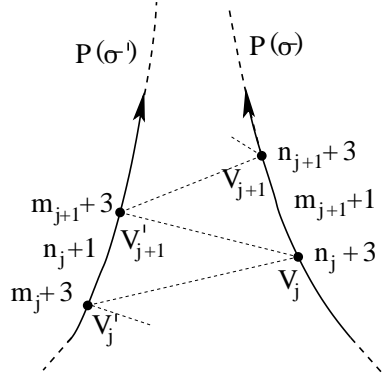, height= 50 mm}
 \caption{Local aspect of the zigzag diagram}
\end{figure}

As an edge and its opposite vertex are dual through the morphism
$\mathcal{I}$ (see Proposition \ref{parallel}) and its analog
$\mathcal{I}'$ attached to the cone $\sigma'$, the triangles appearing
in the zigzag diagram are a convenient graphical representation of
the duality explained in section \ref{comp1}.

$\bullet$ \textit{When $l_-$ is rational and $l_+$ is irrational}, we
draw a little differently the diagram (see Figure 8). The curves representing
$P(\sigma)$ and $P(\sigma')$ start from points $V_0$ and $V_0'$ of a
horizontal line representing the line which contains $l_-$. We 
represent the integral point $V_1'$ differently from the points $V_2',
V_3',...$, because it may not
be a vertex of $P(\sigma')$, as explained in Proposition
\ref{parallel}. The length of $[V_0'V_1']$ is always $1$. The relation
between the length of an edge and the weight of the opposite vertex is
the same as before, with the exception of the triangle $V_1'V_0V_1$,
where the weight of $V_1'$ is equal to $l_{\Z}[V_0V_1]+1$.

{\tt    \setlength{\unitlength}{0.92pt}}
\begin{figure}
   \epsfig{file=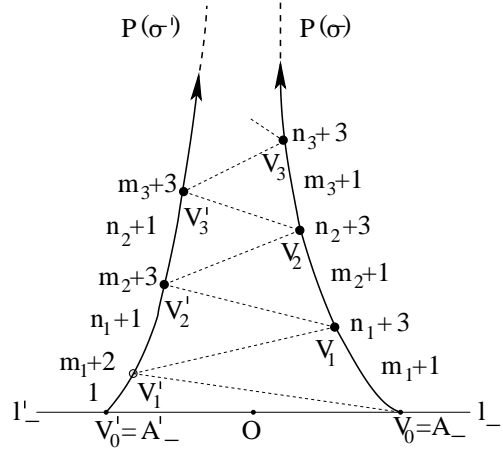, height= 60 mm}
   \caption{The zigzag diagram when $l_-$ is rational}
\end{figure}

$\bullet$  \textit{When both $l_-$ and $l_+$ are rational and there is
at least one vertex on $P(\sigma)$ lying strictly between $A_-$ and $A_+$
(that is, $s\geq 1$)}, the curves representing $P(\sigma)$ and
$P(\sigma')$  start again  from a horizontal line, but now they join
in a point $A_+$ (see Figure 9).

{\tt    \setlength{\unitlength}{0.92pt}}
\begin{figure}
   \epsfig{file=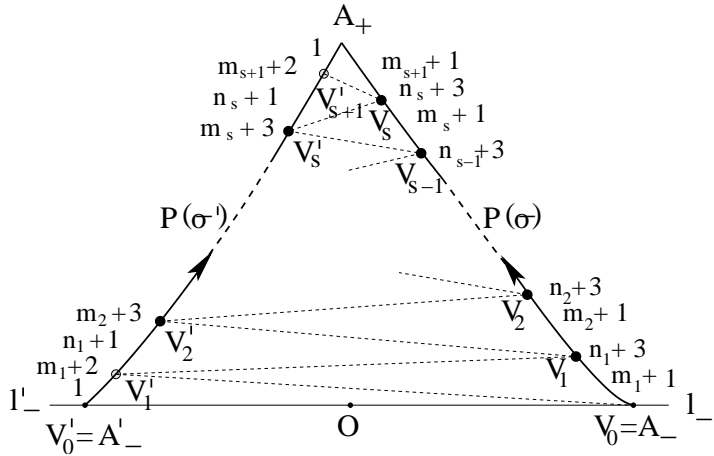, height= 60 mm}
   \caption{The zigzag diagram when both $l_-$ and $l_+$ are rational}
\end{figure}

$\bullet$  \textit{When both $l_-$ and $l_+$ are rational and
  $[A_-A_+]$ is an edge of $P(\sigma)$ (that is, $s=0$)}, the diagram
  is represented in Figure 10.

{\tt    \setlength{\unitlength}{0.92pt}}
\begin{figure}
  \epsfig{file=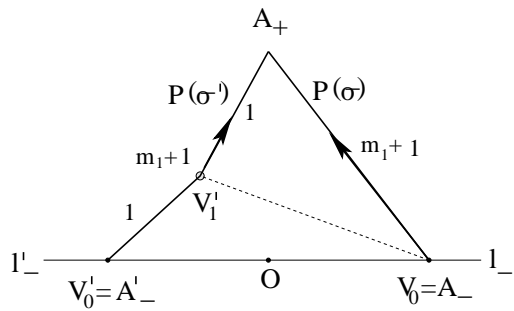, height= 40 mm}
\caption{The zigzag diagram when $P(\sigma)$ has only one compact edge}
\end{figure}

To summarize, we have the following procedure for constructing and
decorating the diagram when $l_-$ is rational: 

\textbf{Procedure: }\textit{Suppose that $l_-$ is rational. Then draw a
  horizontal line  with three marked points $V_0'=A_-', O, V_0=A_-$
  in this order, $V_0'$ on the left and $V_0$ on the right. Starting
  from $V_0'$ and $V_0$, draw in the upper 
  half-plane two curves $P(\sigma')$, respectively
  $P(\sigma)$, concave towards $0$ and coming closer and closer from
  one another. If $l_+$ is rational, join them in a point $A_+$. Draw
  a zigzag line 
  starting from $V_0$ and going alternatively from $P(\sigma)$ to
  $P(\sigma')$. Denote its successive vertices by $V_1', V_1,
  V_2',...$ and stop at the point $V_{s+1}'$. Decorate the edges
  $V_0'V_1'$ and $V_{s+1}'A_+$ by 
  $1$. The other edges and vertices will be decorated using the
  initial data (discussed in the sequel), by respecting the following
  rule:}

\textbf{Rule: } \textit{The weight of a vertex is equal to the length
  of the opposite edge augmented by the number of its vertices
  distinct from the points $A_-,A_-',A_+$.}

\textbf{Initial data: } \textit{If $\sigma$ is of type $\lambda$,
  write the HJ-continued fraction expansion of $\lambda$ in the form:
\begin{equation} \label{CF1}  
   \lambda= [(2)^{m_1}, n_1 +3, (2)^{m_2}, n_2 +3,...]^-
 \end{equation}
Then decorate the edges of $P(\sigma)$ with the numbers $m_1 +1, m_2
 +1,...$ and the vertices with the numbers $n_1 +3, n_2 +3,...$. }

\begin{definition} \label{zigzag}
  We call the previous diagram the \textbf{zigzag diagram} associated
  to the pair $(L, \sigma)$ and to the chosen ordering of the edges of
  $\sigma$, or to the number $\lambda >1$, where $(L, \sigma)$ is of
  type $\lambda$ with respect to this ordering. We denote it by
  $ZZ(\lambda)$.  
\end{definition}

\medskip

The zigzag diagrams allow one  to visualize the relations between
Euclidean and Hirzebruch-Jung continued fractions, proved
algebraically in section \ref{algecomp}. Indeed, one can read the
HJ-continued fraction expansion of $\lambda >1$ on the right-hand
curved line of $ZZ(\lambda)$. By Corollary \ref{slopes}, we can read
the HJ-continued fraction expansion of $\frac{\lambda}{\lambda -1}$ on
the left-hand curved line $P(\sigma)$ of $ZZ(\lambda)$. So, by looking
at Figure 
9, which can be easily constructed from the initial data by respecting
the rule, we get:
 \begin{equation} \label{CF2}  
   \frac{\lambda}{\lambda -1}= [m_1 +2, (2)^{n_1},m_2 +3,(2)^{n_2},
    m_3 +3,...]^-
 \end{equation}
which gives a geometric proof of Proposition \ref{exneg}.

Now, by Klein's geometric interpretation of E-continued fractions (see
section \ref{Kleinint}),
we see that the E-continued fraction expansion of
$\frac{\lambda}{\lambda -1}$ can be obtained by writing alternatively
the integral lengths of the edges of the polygonal lines $P(\sigma)$
and $P(\sigma') -[V_0'V_1']$ (indeed, $\frac{\lambda}{\lambda -1}$ is
the slope of $l_+$ in the base $(OV_0, OV_1')$):
  \begin{equation} \label{CF3}  
    \frac{\lambda}{\lambda -1}= [m_1 +1, n_1 +1, m_2 +1, n_2 +1, m_3
    +1,...]^+.
  \end{equation}
This proves geometrically Proposition \ref{transf}. 

In order to read the E-continued fraction expansion of $\lambda$ on
the diagram, one has to look at $ZZ(\lambda)$ from left to right
instead of from right to left and draw a new zigzag line starting
from $V_0'$. The important point here is that one has to discuss
according to the alternative $m_1=0$ or $m_1 >0$. In the first case,
the zigzag line joins $V_0'$ to $V_1$ and $V_1$ to $V_2'$. In the
second case, it joins $V_0'$ to a new point representing $A_1$ and
$A_1$ to $V_1'$. Compare this with Lemma \ref{expos}.

\begin{example} \label{allinone}
  Take $\lambda=\frac{11}{7}$. After computing $\lambda=[2,3,2,2]^-$,
  we can construct the associated zigzag diagram
  $ZZ(\frac{11}{7})$. We see that the extreme points $V_1', V_2'$ are
  vertices of $P(\sigma')$. 
One can read on it the results of the Examples \ref{extrans},
\ref{exinv}, \ref{exdiag}. 

{\tt    \setlength{\unitlength}{0.92pt}}
\begin{figure}
 \epsfig{file=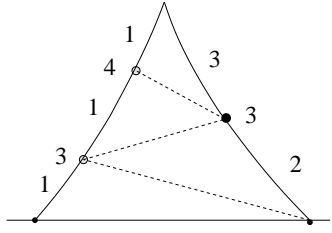, height= 30 mm}
\caption{The first illustration for Example \ref{allinone}: $ZZ(\frac{11}{7})$}
\end{figure}

If one had starts instead from
$\lambda= \frac{11}{4}=[3,4]^-$, the corresponding
diagram would be $ZZ(\frac{11}{4})$. In this case the extreme
points are not vertices of $P(\sigma')$, because their weights are
equal to $2$. 

{\tt    \setlength{\unitlength}{0.92pt}}
\begin{figure}
  \epsfig{file=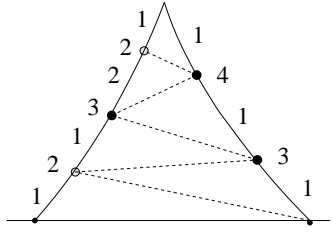, height= 30 mm} 
  \caption{The second illustration for Example \ref{allinone}:
    $ZZ(\frac{11}{4})$} 
\end{figure}

\end{example}

\medskip

\subsection{Relation with the dual cone} \label{reldual}
$ \:$
\medskip
   
Denote by $\check{L}:= \mathrm{Hom}(L, \mathbf{Z})$ the \textit{dual
  lattice} of 
$L$. Inside the associated vector space $\check{L}_\R$ lives the
  \textit{dual cone} $\check{\sigma}$ of  $\sigma$, defined by:
  $$\check{\sigma}:= \{ \check{u} \in \check{L}_\R \: | \: \check{u} .u
  \geq 0, \: \forall   \: u \in \sigma\}.$$

Let $\omega$ be the volume form on $L_\R$
which verifies $\omega(u_1, u_2) =1$ for any basis $(u_1, u_2)$ of $L$
defining the \textit{opposite} orientation to $(l_-, l_+)$. It is a
symplectic form, that is, a non-degenerate alternating bilinear form
on $L_\R$. But we prefer to look at it as a morphism (obtained by
making interior products with the elements of $L$):
   $$\omega: L \longrightarrow \check{L}.$$ 

\begin{proposition} \label{isom}
   The mapping $\omega$ realizes an isomorphism between the pairs 
 $(L, \sigma')$ and $(\check{L}, \check{\sigma})$. 
\end{proposition}

\textbf{Proof:} Indeed we have:
  $$
       \omega^{-1}(\check{\sigma}) 
          = \{ u\in L \: | \: \omega(u) \in \check{\sigma}\}=
       \{ u\in L \: | \: \omega(u,v)\geq 0,\: \forall \: v \in L\}=
       \sigma'. $$ 
  While writing the last equality, we used our convention on the
 orientation of $\omega$. Notice that the dual cone $\check{\sigma}$
 can be defined without 
 the help of any orientation, in contrast with the morphism
          $\omega$. \hfill $\Box$\medskip

The previous proposition shows that the construction of the polygonal
line $P(\sigma')$ explained in Proposition \ref{parallel} describes
also the polygonal line $P(\check{\sigma})$. This observation is
crucial when one wants to use zigzag diagrams for understanding computations
with invariants of toric surfaces (see next section). 

It also helps to
understand geometrically the duality between the convex polygons
$K(\sigma)$ and 
$K(\check{\sigma})$ explained in Gonzalez-Sprinberg \cite{GS 77} and
in Oda \cite[pages 27-29]{O 88}. As 
Dimitrios Dais kindly informed us after seeing a version of this paper on
ArXiv, a better 
algebraic understanding of that duality is explained in Dais, Haus
\& Henk \cite[section 3]{DHH 98}. In particular, modulo Proposition
\ref{isom}, the Theorem 3.16 in the previous reference leads easily to an 
algebraic proof of our Proposition \ref{parallel}. 
\medskip

(Added in proof) Emmanuel Giroux has informed us that he had realized
the existence of a duality between supplementary cones (see
\cite[section 1.G]{G 00}).

\medskip

\section{Relations with toric geometry} \label{toric}

First we introduce elementary notions of toric geometry (see section
\ref{elemtor}). 
In section \ref{torsurf} we explain how to get combinatorially
various invariants of 
a normal affine toric surface and of the corresponding Hirzebruch-Jung
analytic surface singularities. In
Section \ref{moncurves} we explain how to read the combinatorics
of the minimal embedded resolution of a plane monomial curve on
an associated zigzag diagram.

The basics about resolutions of surface singularities needed in order
to understand this section are recalled in section \ref{resdualgr}.
\medskip

\subsection{Elementary notions of toric geometry} \label{elemtor}
$\: $
\medskip

For details about toric geometry, general
references are the books of Oda \cite{O 88} and Fulton \cite{F
  93}, as well as the first survey of it by Kempf,
Knudson, Mumford \& St. Donat \cite{KKMD 73}. 

In the previous section, our fundamental object of study was a pair
$(L, \sigma)$, where $L$ is a lattice of rank 2 and $\sigma$ is
a strictly convex cone in the 2-dimensional vector space
$L_\R$. 

Suppose now that the lattice $L$ has arbitrary finite rank $d\geq 1$ and that
$\sigma$ is a strictly convex {\em rational} cone in $L_{\R}$. The
pair $(L, \sigma)$ gives rise 
canonically to an affine algebraic variety: 

$$ \mathcal{Z}(L, \sigma):= \mbox{Spec} \: \C[\check{\sigma} \cap
\check{L}].$$

This means that the algebra of regular functions on $\mathcal{Z}(L,
\sigma)$ is generated by the monomials whose exponents are elements of
the semigroup $\check{\sigma} \cap
\check{L}$ of integral points in the dual cone of $\sigma$. If $v \in
\check{\sigma} \cap 
\check{L}$, we formally write such a monomial as $X^v$. One can show
that the variety $\mathcal{Z}(L,\sigma)$ is \textit{normal} (see the
definition at 
the beginning of section \ref{resdualgr}).

The closed points of $\mathcal{Z}(L, \sigma)$ are the morphisms of
semigroups $(\check{\sigma}\cap 
\check{L}, +) \rightarrow (\C, \cdot)$. Among them, those whose image is
contained in $\C^*$ form a $d$-dimensional \textit{algebraic torus}
$\mathcal{T}_L= 
\mbox{Spec}  \: \C[\check{L}]$, that is, a complex algebraic group
isomorphic to $(\C^*)^d$. The elements of $L$ correspond to the
\textit{1-parameter subgroups} of $\mathcal{T}_L$, that is, the group
morphisms $(\C^*,\cdot )\rightarrow (\mathcal{T}_L, \cdot)$. The action of
$\mathcal{T}_L$ on itself by 
multiplication extends canonically to an algebraic action on
$\mathcal{Z}(L, \sigma)$, such that $\mathcal{T}_L$ is the unique open
orbit. If $(\overline{L}, \overline{\sigma})$ is a second pair and
$\phi:\overline{L}\rightarrow L$ is a morphism such that
$\phi(\overline{\sigma})\subset \sigma$, one gets an associated {\em toric
morphism}:
 $$\phi_*: \mathcal{Z}(\overline{L}, \overline{\sigma})\rightarrow
 \mathcal{Z}(L, \sigma)$$
It is birational if and only if $\phi$ realizes an isomorphism between
 $\overline{L}$ and $L$. In this case $\phi_*$ identifies the tori contained
 inside $\mathcal{Z}(\overline{L}, \overline{\sigma})$ and
 $\mathcal{Z}(L, \sigma)$. 

In general:

\begin{definition} \label{torvar}
  Given an algebraic torus $\mathcal{T}$, a \textbf{toric variety}
  $\mathcal{Z}$ is an algebraic variety containing $\mathcal{T}$ as a
  dense Zariski open set and endowed with an action \linebreak $\mathcal{T}
  \times \mathcal{Z} \rightarrow \mathcal{Z}$ which extends the group
  multiplication of $\mathcal{T}$. 
\end{definition}

Oda \cite{O 88} and Fulton \cite{F 93} 
study mainly the \textit{normal} toric varieties. For an introduction
to the study of non-necessarily 
normal toric varieties, one can consult Sturmfels \cite{S 96} and
Gonz{\'a}lez P{\'e}rez \& Teissier \cite{GPT 02}. 

A normal toric variety can be
described combinatorially using $\textit{fans}$, that is finite
families of rational strictly convex cones, closed under the
operations of taking faces or intersections. If $L$ is a lattice and
$\mathcal{F}$ is a fan in $L_{\R}$, we denote by $\mathcal{Z}(L,
\mathcal{F})$ the associated normal toric variety. It is obtained by
glueing the various affine toric varieties $\mathcal{Z}(L, \sigma)$
when $\sigma$ varies among the cones of the fan $\mathcal{F}$. As
glueing maps, one uses the toric birational maps
$\mathcal{Z}(\overline{L}, \overline{\sigma})\rightarrow 
 \mathcal{Z}(L, \sigma)$ induced by the inclusion morphisms $(L,
 \overline{\sigma})\rightarrow (L, \sigma)$, for each pair
 $\overline{\sigma}\subset \sigma$ of cones of $\mathcal{F}$. 

The variety $\mathcal{Z}(L, \mathcal{F})$ is smooth if and only if
each cone of the fan $\mathcal{F}$ is {\em regular}, that is,
generated by a subset of a basis of the lattice $L$. 

\medskip

\subsection{Toric surfaces} \label{torsurf}
$\: $
\medskip

We restrict now to the case of surfaces.  Consider a 2-dimensional normal
toric surface $\mathcal{Z}(L, \sigma)$, where $\sigma$ is a strictly
convex cone with non-empty interior.  There is a unique
0-dimensional orbit $O$, whose maximal ideal is generated by the
monomials with exponents in the semigroup $\check{\sigma}\cap
  \check{L} -O$. The surface is smooth outside $O$, and $O$ is a
  smooth point of it if and only if $\sigma$ is a regular
  cone. Supposing that $\sigma$ is not regular, we explain how to
  describe combinatorially the minimal resolution morphism of
  $\mathcal{Z}(L, \sigma)$ and the effect of blowing-up the point
  $O$. We also give a formula for the embedding dimension of the germ
  $(\mathcal{Z}(L, \sigma),O)$, which is a so-called {\em Hirzebruch-Jung
  singularity}.  
\medskip

With the notations of section 4, let us subdivide $\sigma$ by drawing
the half-lines starting from 
$O$ and passing through the points $A_k, \: \forall \: k \in
\{1,...,r\}$. In this way we decompose $\sigma$ in a finite number of
regular subcones. They form \textit{the minimal regular subdivision}
of $\sigma$, in the sense that any subdivision of $\sigma$ by regular
cones is necessarily a refinement of the preceding one. 

The family consisting of the 2-dimensional cones in the subdivision,
of their edges and of the origin form a fan
$\mathcal{F}(\sigma)$. For each such subcone $\sigma'$ of $\sigma$,
there is a canonical birational morphism $\mathcal{Z}(L, \sigma')
\rightarrow \mathcal{Z}(L, \sigma)$, which realizes an isomorphism of
the tori. Using these morphisms, one can glue canonically the tori
contained in the surfaces $\mathcal{Z}(L, \sigma')$ when $\sigma'$
varies,  and obtain a new
toric surface $  \mathcal{Z}(L, \mathcal{F}(\sigma))$, endowed with a
morphism:
  $$ \mathcal{Z}(L, \mathcal{F}(\sigma))
  \stackrel{p_\sigma}{\longrightarrow} \mathcal{Z}(L, \sigma)$$

\begin{proposition} \label{resmintor}
The morphism $p_\sigma$ is the
minimal resolution of singularities of the surface $\mathcal{Z}(L,
\sigma)$. Moreover, its exceptional locus $E_\sigma$
is a normal crossings 
divisor and the dual graph of $E_\sigma$ is topologically  a segment. 
\end{proposition}

\textbf{Proof:} For details, see \cite{F 93}. Here we outline only the
  main steps. The morphism $p_\sigma$ is proper, birational and realizes an
  isomorphism over 
  $\mathcal{Z}(L, \sigma)-O$. As $ \mathcal{Z}(L,
  \mathcal{F}(\sigma))$ is smooth, $p_\sigma$ is a a resolution of
  singularities of $\mathcal{Z}(L, \sigma)$ (see Definition
  \ref{resdef}). There is a canonical
  bijection between the irreducible
  components $E_k$ of the exceptional divisor
  $E_\sigma=p_\sigma^{-1}(0)$ and the half-lines $[OA_k$, for $k\in
  \{1,...,r\}$. Moreover, $E_k$ is a smooth compact rational curve and

\begin{equation} \label{selfint}
   E_k^2 =-\alpha_k, \: \forall \: k \in \{1,...,r\}
\end{equation}
where the numbers $\alpha_k$ were introduced in relation (\ref{eq3}). 

Using the inequality (\ref{ineq}), we deduce that no component of
$E_{\sigma}$ is exceptional of the first kind (see the comments which
follow Definition \ref{resdef}). This implies that $p_{\sigma}$ is
the minimal resolution of singularities of $\mathcal{Z}(L,
\sigma)$. The proposition is proved.    \hfill $\Box$

\medskip

Notice that relation (\ref{selfint}) gives an intersection-theoretical  
interpretation of the weights attached through relation (\ref{eq3}) to
the integral points situated on $P(\sigma)$ which are interior to $\sigma$.

Conversely (see \cite{BHPV 04} and \cite{PP 05}):

\begin{proposition} \label{convHJ}
 Suppose that  a smooth surface $\mathcal{R}$ contains a compact normal
crossings divisor $E$ whose components are  smooth rational curves of
self-intersection $\leq -2$ and
whose dual graph is topologically a segment. Denote by $\alpha_1,...,
\alpha_r$ the self-intersection numbers read orderly along the
segment.  Then $E$ can be
contracted by a map $p: (\mathcal{R} ,E) \rightarrow (\mathcal{S}, 0)$ to
a normal surface $\mathcal{S}$ and the germ $(\mathcal{S}, 0)$ is
analytically isomorphic to a germ of the form $(\mathcal{Z}(L,
\sigma), O)$, where $\sigma$ is of type $\lambda := [\alpha_1,...,
\alpha_r]^-$. 
\end{proposition}

 This motivates:

\begin{definition} \label{HJsing}
  A normal surface singularity $(\mathcal{S}, 0)$ isomorphic 
  to a  germ of the form $(\mathcal{Z}(L, \sigma), O)$ is called a
  \textbf{Hirzebruch-Jung singularity}. 
\end{definition}

Hirzebruch-Jung singularities can also be defined as \textit{cyclic
  quotient singularities} (see \cite{BHPV 04} and \cite{PP 05}). 
They  appear naturally in the so-called \textit{Hirzebruch-Jung
method} of studying an arbitrary surface singularity. Namely, one projects the
given singularity by a finite morphism on a smooth surface, then one
makes an embedded resolution of the discriminant curve and takes the
pull-back of the initial surface by this morphism. In this case, the
normalization of the new surface has only Hirzebruch-Jung
singularities (see Laufer \cite{L 71},
Lipman \cite{L 75}, Brieskorn \cite{B 00} for details  and
  Popescu-Pampu \cite{PP 05} for a
generalization to higher dimensions). 

The proof of Proposition \ref{resmintor} shows that the germs
$(\mathcal{Z}(L, \sigma), O)$ and $(\mathcal{Z}(\overline{L},
\overline{\sigma}), O)$ are analytically isomorphic if and only if
there exists an isomorphism of the lattices $L$ and $\overline{L}$
sending $\sigma$ onto $\overline{\sigma}$. The same is true for
strictly convex cones in arbitrary dimensions, as proved by Gonz{\'a}lez
P{\'e}rez \& Gonzalez-Sprinberg \cite{GPGS 04}. Previously we had proved
this for simplicial cones in \cite{PP 05}. 

A Hirzebruch-Jung singularity isomorphic to $(\mathcal{Z}(L, \sigma),
O)$ is said to be \textit{of type }
$\mathcal{A}_{p,q}$, with $1 \leq q < p$ and $\mbox{gcd}(p,q)=1$ if
(using Definition \ref{typecone}) the pair $(L, \sigma)$ is of type
$\frac{p}{q}$ with respect to one of the orderings of the sides of
$\sigma$. Then, by Proposition \ref{concrete}, we have 
$\frac{p}{q}=[\alpha_1,...,\alpha_r]^-$. By Proposition
\ref{changeorder}, one has $\mathcal{A}_{p,q}\simeq \mathcal{A}_{p',
  q'}$ if and only if $p=p'$ and $q'\in \{q, \overline{q}\}$, where
$q\overline{q}\equiv 1 \: (mod \: p)$.

The singularities of type
$\mathcal{A}_{n+1,n}$ are also called \textit{of type} $\mathbf{A}_n$. They are
those for which the polygonal line $P(\sigma)$ has only one compact edge, as
$\frac{n+1}{n}= [(2)^n]^-$ (a case emphasized in Section
\ref{zigdiag}), and also the only Hirzebruch-Jung 
singularities of embedding dimension $3$ (more precisely, they can be
defined by the equation $z^{n+1}=xy$). Indeed:

\begin{proposition} \label{embdim}
  If $\dfrac{p}{q}= [\alpha_1,..., \alpha_r]^-= [(2)^{m_1}, n_1
  +3,..., n_{s}+3, (2)^{m_{s+1}}]^-$, 
  then:
    $$ embdim(\mathcal{A}_{p,q}) = 3 + \sum_{i=1}^r (\alpha_i -2)=\:  3 + s +
    \sum_{k=1}^{s} n_k.$$ 
\end{proposition}

\textbf{Proof:} 
  If $S$ is a generating 
  system of the semigroup $\check{L} \cap \check{\sigma} -O$, then the
  monomials $(X^v)_{v \in S}$ form a
  generating system of the Zariski cotangent space
  $\mathcal{M}/\mathcal{M}^2$ of the germ at the singular point, where
  $\mathcal{M}$ is the maximal ideal of the local algebra of the
  singularity $\mathcal{A}_{p,q}$. By
  taking a minimal generating system, one gets a basis of
  this cotangent space. But such a minimal generating system is
  unique, and consists precisely of the integral points of 
  $P(\check{\sigma})$ interior to $\check{\sigma}$. By Propositions 
  \ref{isom}  and \ref{exneg}, we see that this number is as given in
  the Proposition.\hfill $\Box$

\medskip

Hirzebruch-Jung singularities are particular cases of \textit{rational
  singularities}, introduced by M. Artin \cite{A 62}, \cite{A 66} in
  the 60's (see also \cite{BHPV 
  04}). In \cite{T 68}, Tjurina proved that the blow-up of a rational
  surface singularity is a normal surface which has again only
  rational singularities (see also the comments of L{\^e} \cite[
  4.1]{L 00}). As any surface can be 
  desingularized by a sequence of blow-ups of its
  singular points followed by normalizations (Zariski \cite{Z 39},
  see also Cossart \cite{C 00} and the references therein), this shows that a 
  rational singularity can be desingularized by a sequence of blow-ups
  of closed points. 
In particular this is true for a Hirzebruch-Jung singularity. As the
operation of blow-up is analytically invariant, we can describe the
blow-up of $O$ in the model surface $\mathcal{Z}(L,\sigma)$. We use
  notations introduced at the beginning of the proof of Proposition
  \ref{parallel}.

\begin{proposition} \label{blowup}
   Suppose that the cone $\sigma$ is not regular. 
   Subdivide it by drawing the half-lines starting from
   $O$ and passing through the points $A_1, V_1,$ $V_2...,V_s, A_r$. Denote
   by $\mathcal{F}_0(\sigma)$ the fan obtained in this way. Then the
   natural toric morphism 
     $\mathcal{Z}(L,\mathcal{F}_0(\sigma)) \stackrel{p_0}{\longrightarrow}
   \mathcal{Z}(L,\sigma)$
   is the blow-up of $O$ in $\mathcal{Z}(L,\sigma)$. 
\end{proposition}

\textbf{Proof:} A proof is sketched by Lipman in \cite{L 75}.  
Here we give more details. 

Let $(\mathcal{S},0)$ be any germ of normal surface. Consider its
minimal resolution $p_{min}: (\mathcal{R}_{min}, E_{min})
\rightarrow (\mathcal{S},0)$ and its exceptional divisor
$E_{min}=\sum_{k=1}^r E_k$. The divisors $Z \in \sum_{k=1}^r \Z E_k$
which satisfy $Z\cdot E_k \leq 0, \: \forall \: k \in \{1,...,r\}$
form an additive semigroup with a unique minimal element $Z_{top}$,
called the \textit{fundamental cycle} of the singularity. It verifies
\begin{equation} \label{fundcyc}
 Z_{top} \geq \sum_{k=1}^{r}E_k
\end{equation} 
for the componentwise order on the set of cycles with integral
coefficients. \textit{In the case of a  rational singularity}, Tjurina
\cite{T 68} showed that the divisors $E_k$ which appear in the blow-up
of $0$ on $\mathcal{S}$ can be characterized using the fundamental
cycle: they are precisely those for which $Z_{top} \cdot E_k <0$. 

In our case, where $(\mathcal{S},0)= (\mathcal{Z}(L,\sigma),O)$,
Proposition \ref{resmintor} shows that $p_{min}=
p_{\sigma}$. Using the relations (\ref{selfint}) and
(\ref{fundcyc}), we see that $Z_{top}= \sum_{k=1}^{r}E_k$. Again using
relation (\ref{selfint}), we get:
 $$Z_{top} \cdot E_k <0 \: \Longleftrightarrow \: \mbox{either } k \in
 \{1,r\} \mbox{ or } \alpha_k \geq 3.$$
This shows that the components of $E_\sigma$ which appear when one
 blows-up the origin, are precisely those which correspond to the
 half-lines $[OA_1, [OV_1,$ $[OV_2,...,[OV_s,$ $[OA_r$. But the surface
 obtained by blowing-up the origin is again normal, by Tjurina's
 theorem, which shows that it coincides with $
 \mathcal{Z}(L,\sigma)$. \hfill $\Box$

\medskip

One sees that after the first blow-up, the new surface has only
singularities of type $\mathbf{A}_n$, where $n$ varies in a finite set
of positive numbers. The singular points are
contained in the set of $0$-dimensional orbits of the toric surface
$\mathcal{Z}(L,\mathcal{F}_0(\sigma))$, which in turn correspond
bijectively to the 2-dimensional cones of the fan
$\mathcal{F}_0(\sigma)$. The germs of the surface at those points are
Hirzebruch-Jung singularities of types $\A_{n_0},..., \A_{n_s}$, where 
$n_0=l_\Z[A_1 V_1], n_1=l_\Z[V_1 V_2],..., n_s= l_\Z[V_s A_r]$. 
\medskip

We have spoken until now of algebraic aspects of Hirzebruch-Jung
singularities. We discuss their topology in section \ref{topchar}. 
\medskip

\subsection{Monomial plane curves} \label{moncurves}
$\:$
\medskip

Suppose that $(\mathcal{S},0)$ is a germ of \textit{smooth} surface
and that $(\mathcal{C},0)\subset (\mathcal{S},0)$ is a germ of reduced
curve. A proper birational morphism $p:\mathcal{R} \rightarrow \mathcal{S}$
is called an \textit{embedded resolution} of the germ
$(\mathcal{C},0)$ if $\mathcal{R}$ is smooth, $p$ is an isomorphism above
$\mathcal{S}-0$ and 
the \textit{total transform} $p^{-1}(\mathcal{C})$ of $\mathcal{C}$ is
a divisor with 
normal crossings on $\mathcal{R}$ in a neighborhood of the
\textit{exceptional divisor} $E:= p^{-1}(0)$. The difference
$p^{-1}(\mathcal{C})- p^{-1}(0)$ is called the \textit{strict
  transform} of $\mathcal{C}$ 
by the morphism $p$. 

It is known since the XIX-th century that any germ of plane curve can
be resolved in an embedded way by a sequence of blow-ups of points
(see Enriques \& Chisini \cite{EC 18}, Laufer \cite{L 71}, Brieskorn
\& Kn{\"o}rrer \cite{BK 86}). The combinatorics of the exceptional
divisor of the resolution can be determined starting from the
Newton-Puiseux exponents of the irreducible components of the curve
and from their intersection numbers using E-continued fraction 
expansions. We explain here how to read the sequence of
self-intersection numbers of the components of the exceptional divisor
of the minimal embedded resolution of a monomial plane curve by using
a zigzag diagram, instead of just 
doing blindly computations with continued fractions. 
\medskip

If $p,q \in \N^*, \: 1\leq q <p$ and $gcd(p,q)=1$, consider the plane
curve $C_{p/q}$ defined by the equation:
\begin{equation} \label{eqmon}
  x^p-y^q=0
\end{equation}
It can be parametrized by:
\begin{equation} \label{param}
    \left\{ \begin{array}{l}
                x=t^q\\
                y=t^p
    \end{array}\right.
\end{equation}

As $p$ and $q$ are relatively prime, one sees that (\ref{param})
describes the  normalization morphism for $C_{p/q}$ (see its
definition at the beginning of section \ref{resdualgr}). As $t^p$ and $t^q$
are monomials, one says that $C_{p/q}$ is a 
\textit{monomial curve}. There is a natural generalization  to higher
dimensions (see Teissier  \cite{T 04}). 

If one identifies the plane $\C^2$ of coordinates $(x,y)$ with the toric
surface $\mathcal{Z}(L_0, \sigma_0)$, where $L_0= \Z^2$ and $\sigma_0$
is the first quadrant, then it is easy to see (look at equation
(\ref{param})) that $C_{p/q}$ is the closure in $\C^2$ of the image of
the 
1-parameter subgroup of the complex torus $\mathcal{T}_{L_0}=(\C^*)^2$
corresponding to the point $(q,p)$. 

Consider again the notations introduced before Lemma 
\ref{firstedge}. Let $l_-:= [O \: (1,0)$ and $l_+:= [O \: (q,p)$ be the
edges of the cone $\sigma_x(\frac{p}{q})$. We leave to the reader the
proof of the following lemma, which is very similar to the proof of Lemma
\ref{firstedge}. Recall that the {\em type} of a cone was introduced in
Definition \ref{typecone}.

\begin{lemma} \label{identype}
  With respect to the chosen ordering of its edges, the cone
  $\sigma_x(\frac{p}{q})$ is of type $\frac{p}{p-q}$. Moreover, with
  the notations of section \ref{geomcomp}, $A_1=(1,1)$, 
  $A_1'=(0,1)$ and $A_+= (q,p)$. 
\end{lemma}

Even if the proof is very easy, it
  is important to be conscious of this result, as it allows to apply
  the study done in section \ref{geomcomp}  to our context.

Given the pair $(p,q)$, we want to describe the process of embedded
resolution of the curve $C_{p/q}$ by blow-ups, as well as the final
exceptional divisor, the self-intersections of its components and
their orders of appearance during the process. 

\begin{lemma} \label{ordap}
  The blow-up $\pi_0 : \mathcal{R}_0 \rightarrow \C^2$ of $0$ in
  $\C^2$ is a toric morphism corresponding to the subdivision of
  $\sigma_0$ obtained by joining $O$ to $A_1= (1,1)$. The strict
  transform of $C_{p/q}$ passes through the $0$-dimensional orbit of
  $\mathcal{R}_0$ associated to the cone $\R_+ OA_1 + \R_+ OA_1'$. 
\end{lemma}

\textbf{Proof: } With the notations of Section \ref{Kleinint}, we
consider the fan $\mathcal{F}_0$ subdividing $\sigma_0$ which consists
of the cones $\sigma_x(1), \sigma_y(1)$, their edges and the
origin. Let $\pi_{\mathcal{F}_0}: \mathcal{Z}(L, \mathcal{F}_0)
\rightarrow \mathcal{Z}(L, \sigma_0)$ be the associated toric
morphism. It is obtained by gluing the maps $\pi_x : \mathcal{Z}(L,
\sigma_x(1))\rightarrow \mathcal{Z}(L, \sigma_0)$  and $\pi_y : \mathcal{Z}(L,
\sigma_y(1))\rightarrow \mathcal{Z}(L, \sigma_0)$  over
$(\C^*)^2$. With respect to the coordinates given by the monomials
associated to the primitive vectors of $L$ situated on the edges of
the cones $\sigma_0, \sigma_x(1), \sigma_y(1)$, the maps $\pi_x$ and
$\pi_y$ are respectively described by:
 $$\left\{ \begin{array}{l}
              x= x_1 y_1\\
              y= y_1
           \end{array}  \right. \mathrm{and} \:  
    \left\{ \begin{array}{l}
              x= x_2\\
              y= x_2 y_2
           \end{array}    \right. $$
One recognizes the blow-up of $0$ in $\C^2$. Now, in order to compute
the strict transform of $C_{p/q}$, one has to make the previous
changes of variables in equation (19). The lemma follows immediately.
\hfill $\Box$
\medskip

Starting from Lemma \ref{identype} and using the previous lemma as an
induction step, we get:

\begin{proposition} \label{algres}
 The following
procedure constructs the dual graph of the total transform of
$C_{p/q}$ by the minimal embedded resolution morphism, starting from
the zigzag diagram $ZZ(\frac{p}{p-q})$:

\medskip

$\bullet$ On each edge of integral length $l \geq 1$, add $(l-1)$ vertices of
weight $2$. Then erase the weights of the edges (that is, their integral
length).

$\bullet$  Attach the weight $1$ to the vertex $A_+$. Then change the
signs of all the weights of the vertices.

$\bullet$ Label the vertices by the symbols $E_1, E_2, E_3,...$
starting from $A_1$ on $P(\sigma)$ till arriving at $V_1$, continuing
from the first vertex which follows $V_1'$ on $P(\sigma')$ till
arriving at $V_2'$, coming then back to $P(\sigma)$ at the first
vertex which follows $V_1$ and so on, till labelling the vertex $A_+$.

$\bullet$ Erase the horizontal line, the zigzag line and the curved
segment between $V_0'$ and the first vertex which follows $V_1'$. 

$\bullet$ Add an arrow to the vertex $A_+$ and keep only the
weights of the vertices and their labels $E_n$.

The arrowhead vertex represents the strict transform of the curve
$C_{p/q}$ and the indices of 
  the components $E_i$ correspond to the orders of appearance during
  the process of blow-ups. 
\end{proposition}

It is essential to remark that in the previous construction one starts
from $ZZ(\frac{p}{p-q})$ and 
not from $ZZ(\frac{p}{q})$ (look again at Lemma \ref{identype}).

\medskip

\begin{example} \label{dualgraph}
  Consider the curve $x^{11}-y^4=0$. Then $\lambda= \frac{11}{11-4}=
  \frac{11}{7}$. Its zigzag diagram $ZZ( \frac{11}{7})$ was
  constructed in Example \ref{allinone}. So, the dual graph of the
  total transform of $C_{11/4}$ by the minimal embedded resolution
  morphism has 6 vertices, of easy computable weights (see Figure
  13). 

{\tt    \setlength{\unitlength}{0.92pt}}
\begin{figure}
  \epsfig{file=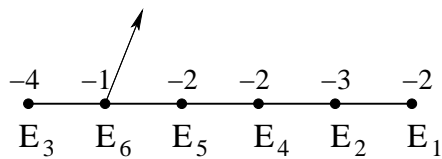, height= 20 mm} 
  \caption{The dual graph of the
  total transform of $C_{11,4}$}
\end{figure}

\end{example}

Proposition \ref{algres} endows us with an easy way of
remembering the following classical description of the minimal
embedded resolution of a monomial plane curve (see Jurkiewicz 
\cite{J 85}, who attributes it to Hirzebruch; Spivakovsky \cite{S 90}
extends it to the case of monomial-type valuations on function-fields
of surfaces):

\begin{proposition} \label{dualgen}
  If $\dfrac{p}{q}= [m_1 +1, n_1 +1, m_2 +1,..., n_s +1, m_{s+1}+1]^+$,
  then the dual graph of the total transform of the monomial curve
  $C_{p/q}$ is the one which appears in Figure 14.
\end{proposition}

\textbf{Proof: } Combine formulae (\ref{CF3}) and (\ref{CF1}) with
Figure 9 and Proposition \ref{algres}. 
  \hfill $\Box$
\medskip

{\tt    \setlength{\unitlength}{0.92pt}}
\begin{figure}
  \epsfig{file=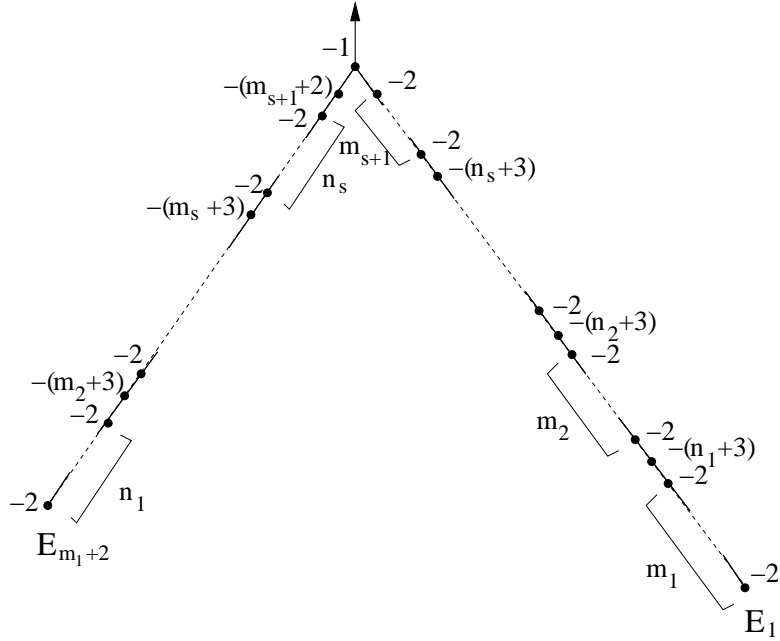, height= 85 mm} 
  \caption{The dual graph of the
  total transform of $C_{p/q}$}
\end{figure}

In Figure 14 we have indicated only the orders of appearance of the
components of 
the exceptional divisor corresponding to the extremities of the
graph. We leave as an exercise for the reader to complete the diagram
with the sequence $(E_k)_{k \geq 1}$. 

Notice that in the E-continued fraction expansion of $\frac{p}{q}$
used in the previous proposition, there is the possibility  that
$m_{s+1}=0$. In this case, the canonical expansion is obtained using
relation (\ref{id1}). But in order to express in a unified form the
result of the application of the algorithm, it was important for us to
use an expansion of $\frac{p}{q}$ with an \textit{odd} number of
partial quotients (which is always possible, precisely according to
formula (\ref{id1})). 

One can use the combinatorics of the embedded resolution of monomial
plane curves as building blocks for the description of the
combinatorics of the resolution of any germ of plane curve. A
detailed description of the passage between \textit{the Eggers tree},
which encodes the Newton-Puiseux exponents of the components of the
curve, and the dual graph of the total transform of the curve by its
embedded resolution morphism can be found in Garc{\'\i}a Barroso
\cite{GB 96} (see also Brieskorn \& Kn{\"o}rrer \cite[section 8.4]{BK 86} and Wall
\cite{W 04}). A topological interpretation 
of the trees appearing in these two encodings was given in
Popescu-Pampu \cite[chapter 4]{PP 01}. 

In higher dimensions, Gonz{\'a}lez P{\'e}rez \cite{GP 03} used toric
geometry in order to describe embedded resolutions of quasi-ordinary
hypersurface singularities. Again, the building blocks are monomial
varieties. A prototype for his study is the method of resolution of an 
irreducible germ of plane curve by only one toric morphism, developed
by Goldin \& Teissier \cite{GT 00}. 

In the classical treatise of Enriques \& Chisini \cite{EC 18},
resolutions of curves by blow-ups of points are not studied using
combinatorics of divisors, but instead using the \textit{infinitely near
  points} through which the strict transforms of the curve pass during
the process of blowing ups. Those combinatorics were also encoded in a
diagram, called nowadays \textit{Enriques diagram} (see Casas-Alvero
\cite{CA 00}). Enriques diagrams are very easily constructed using the
knowledge of the orders of appearance of the divisors during the
process of blowing ups. For this reason, zigzag diagrams combined with
Proposition \ref{algres} give an easy way to draw them for a
monomial plane curve. We leave the
details to the interested reader. Then one uses this again as building
blocks for the 
analysis of general plane curve singularities (see \cite{CA 00}).

\medskip

\section{Graph structures and plumbing structures on 3-manifolds}
\label{graplumb} 

This section contains preparatory material for the topological study of the
3-manifolds appearing as abstract boundaries of normal surface
singularities, done in sections \ref{general} and \ref{canplumb}.

We recall 
general facts about Seifert, graph and plumbing structures on
3-manifolds, as well as about JSJ theory. We also 
define  particular classes of plumbing structures on thick tori and
solid tori, starting from naturally arising pairs $(L,
\sigma)$, where $L$ is a 2-dimensional lattice and $\sigma$ is a
rational strictly convex cone in $L_{\R}$. Namely, given a pair of
essential curves on 
the boundary of a thick torus $M$, their classes generate two lines in
the lattice $L:=H_1(M, \Z)$. A choice of orientations of these lines
distinguishes one of the four cones in which the lines divide
the plane... 

\medskip

\subsection{Generalities on manifolds and their splittings} 
$\:$
\medskip

 We denote by $\I$ the interval $[0,1]$, by $\D$
the closed disc of 
dimension $2$ and by $\textbf{S}^n$ the sphere of dimension $n$. An
\textit{annulus} is a surface diffeomorphic to $I \times \mathbf{S}^1$. 

A simple closed 
curve on a 2-dimensional torus is called \textit{essential} if it is
non-contractible. It is classical that an oriented essential curve on a torus
$T$ is determined up to isotopy by its image in $H_1(T,
\Z)$ (see \cite[section 2.3]{FM 97}). Moreover, the vectors of $H_1(T,
\Z)$ which are homology 
classes of essential curves are precisely the primitive ones.

We say that a manifold is \textit{closed} if it is compact and without
boundary. If $M$ is a
manifold with boundary, we denote by $\stackrel{\circ}{M}$ its
interior and by $\partial M$ its boundary. If moreover $M$ is
oriented, we orient 
$\partial M$ in such a way that at a point of $\partial M$, an outward
pointing tangent vector to 
$M$, followed by a basis of the tangent
space to $\partial M$, gives a basis of the tangent space to $M$ (this
is the convention which makes Stokes' theorem $\int_M d \omega = 
\int_{\partial M} \omega$ true). We say then that $\partial M$ is
oriented \textit{compatibly with} $M$. 

If $M$ is an oriented manifold, we denote by $-M$ the same manifold
with reversed orientation. If $M$ is a closed oriented surface, then
$-M$ is orientation-preserving diffeomorphic to $M$. This fact is no
longer true in dimension 3, that is why it is important to describe 
carefully the choice of orientation. In this sense, see Theorem
\ref{changeor}, as well as Propositions \ref{orientchange} and
\ref{dualcusp}. 

We denote by $\operatorname{Diff}(M)$ the group of
self-diffeomorphisms of $M$, by 
$\operatorname{Diff}^{\circ}(M)$ the subgroup of self-diffeomorphisms which are
isotopic to the identity and by $\operatorname{Diff}^+(M)$ the subgroup of
diffeomorphisms which preserve the orientation of $M$ (when $M$ is
orientable).

\begin{definition} \label{thickfill}
Let $M$ be  a 3-manifold with boundary. We say that $M$ is a
\textbf{thick torus} if it is 
diffeomorphic to  $\mathbf{S}^1\times\mathbf{S}^1 \times \I$. We say
that $M$ is a 
\textbf{solid torus} if it is diffeomorphic to $\D \times
\mathbf{S}^1$. We say that $M$ is a \textbf{thick Klein bottle} if it
is diffeomorphic to a unit tangent circle bundle to the M{\"o}bius
band. 
\end{definition}

In the definition of a thick Klein bottle $M$ we use an arbitrary
riemannian metric on a M{\"o}bius band. The manifold obtained like this is
independent of the choices up to diffeomorphism. Moreover, it is
orientable, because any tangent bundle is orientable and the manifold
we define appears as
the boundary of a unit tangent disc bundle. The preimage of a central circle
of the M{\"o}bius band by the fibration map is a Klein bottle, and the
manifold $M$ appears then as a tubular neighborhood of it, which
explains the name. For details, see \cite[section 3]{W 67} and
\cite[section 10.11]{FM 97}.

On the boundary of a solid torus $M$ there exists an essential curve
which is 
contractible in $M$. Such a curve, which is unique up to
isotopy (see \cite{FM 97}), is called a \textit{meridian} of $M$. A
3-manifold $M$ is called 
\textit{irreducible} if any embedded sphere bounds a ball. A surface
embedded in $M$ is called \textit{incompressible} if its $\pi_1$
injects in $\pi_1(M)$. Two tori embedded in
$M$ are called \textit{parallel} if they are disjoint and they cobound
a thick torus embedded in $M$. The manifold $M$ is called
\textit{atoroidal} if any embedded incompressible torus is parallel to
a component of $\partial M$.

\medskip

\begin{definition} \label{split}
Let $M$ be an orientable manifold and $S$ be
an orientable closed (not necessarily connected) hypersurface of 
$M$. 
 A manifold with boundary $M_S$ endowed with a map $M_S
 \stackrel{r_{M,S}}{\longrightarrow} M$ is called \textbf{a splitting
 of $M$ along $S$}  if:

$\bullet$ $r_{M,S}$ is a local embedding;

$\bullet$ $\partial M_S = (r_{M,S})^{-1}(S)$ and the restriction $r_{M,S}
|_{\partial M_S}$ is a trivial double covering of $S$;

$\bullet$ the restriction $(r_{M,S})|_{\stackrel{\circ}{M}_S}:
\stackrel{\circ}{M}_S \longrightarrow M-S$ is a diffeomorphism.

If this is the case, the map $r_{M,S}$ is called the
\textbf{reconstruction map} associated to the splitting. We say that
$S$ \textbf{splits} $M$ \textbf{into} $M_S$ and that 
the connected components of $M_S$ are the \textbf{pieces} of the
splitting. If $N$ is a piece of $M_S$ and $P \subset M$ is a set, we
say that $P$ \textbf{contains} $N$ if $r_{M,S}(N) \subset P$.  
\end{definition} 

It can be shown easily that splittings of $M$ along $S$ exist and are unique up
to unique isomorphism. The idea is very intuitive, one simply thinks at
$M$ being split open along each connected component of $S$. A way to
realize this is to take the complement of an open tubular neighborhood
of $S$ in $M$ and to deform the inclusion mapping in an arbitrarily
small neighborhood of the boundary in order to push it towards $S$
(see Waldhausen \cite{W 68} and Jaco \cite{J 80}).

If $\phi \in \operatorname{Diff}^+(M)$, one can also canonically split
$\phi$ and get a 
diffeomorphism $\phi_S$ of manifolds with boundary (we leave the
axiomatic definition of $\phi_S$ to the reader):
  $$ \phi_S : M_S \longrightarrow
  M_{\phi(S)}$$  

Among closed 3-manifolds, two particular classes will be especially 
important for us, the lens spaces and the torus fibrations. The reason
why we treat them simultaneously will appear clearly in section
\ref{topchar}.

\begin{definition} \label{lenstor}
  Let $M$ be an orientable 3-manifold. We say that $M$ is a
  \textbf{lens space} if it contains an embedded torus $T$ such that
  $M_T$ is the disjoint union of two solid tori whose meridians have
  non-isotopic images on $T$. We say that $M$ is a
  \textbf{torus fibration} if it contains an embedded torus $T$
  such that $M_T$ is a thick torus. 
\end{definition}

Lens spaces can also be defined as quotients of $\mathbf{S}^3$ by linear free
cyclic actions or - and this explains the name - as  manifolds 
obtained by gluing in a special way the faces of a lens-shaped
polyhedron (see  \cite{ST 80} or \cite[section 4.3]{FM 97}). We impose the
condition on the meridians in order to avoid the manifold $\mathbf{S}^1 \times
\mathbf{S}^2$, which can also be split into two solid tori, but whose
universal cover is not the 3-dimensional sphere, a difference which
makes it to be excluded from the set of lens spaces by most
authors. There exists a classical encoding of oriented lens spaces by
positive integers. We recall it at the end of section \ref{caselens}
(see Proposition \ref{isolens}).

If $M$ is a torus fibration and $T\subset M$ splits it into a thick
torus, then a trivial foliation of $M_T$ by tori parallel to the boundary
components is projected by $r_{M,T}$ onto a foliation by pairwise parallel
tori. The space of leaves is topologically a circle and the projection
$\pi:M \rightarrow \mathbf{S}^1$ is a locally trivial fibre bundle
whose fibres are tori, which explains the name. 

\begin{definition} \label{algmon}
Let $\pi:M \rightarrow \mathbf{S}^1$ be a locally trivial fibre bundle
whose fibres are tori. 
  Fix a fibre of $\pi$ (for
example the initial torus $T$) and also an orientation of the base
space $\mathbf{S}^1$. The \textbf{algebraic monodromy operator} $m$ 
is by definition the first return map of the natural parallel
transport on the first homology fibration over $\mathbf{S}^1$, when
one travels in the positive direction.
\end{definition}

The map $m$ is a well-defined linear automorphism $m \in SL(H_1(T,
\Z))$, once an orientation of $\mathbf{S}^1$ was chosen.  Its conjugacy class
in $SL(2, \Z)$ is independent of the 
choice of the fibre. If one changes the orientation of $\mathbf{S}^1$, then
$m$ is replaced by $m^{-1}$. This shows that the trace of $m$ is independent of
the choice of $T$ and of the orientation of $\mathbf{S}^1$. Remark
that no choice of orientation of $M$ is needed in order to define it. 

For more information about torus fibrations, see Neumann \cite{N 81}
and Hatcher \cite{H 00}. We come back to them in Section
\ref{torfibr}, with special emphasis on subtleties related to their
orientations. 

\medskip

\subsection{Seifert structures}
$\:$
\medskip

Seifert manifolds are special 3-manifolds whose study can be reduced
in some way to the study of lower-dimensional spaces.

\begin{definition} \label{seifert}
   A \textbf{Seifert structure} on a 3-manifold $M$ is a foliation
   by circles such that any leaf has a compact orientable 
   saturated neighborhood. A leaf with trivial holonomy is called a
   \textbf{regular fibre}. A leaf which is not regular is called an 
   \textbf{exceptional fibre}. The space of leaves is called the
   \textbf{base} of the Seifert structure. We say that a Seifert structure is
   \textbf{orientable} if there is a continuous orientation of all the
   leaves of the foliation. If such an orientation is fixed, one says
   that the Seifert structure is \textbf{oriented}. If there exists a
   Seifert structure on 
   $M$, we say that $M$ is a \textbf{Seifert manifold}. 
\end{definition}

The condition on the leaves to have compact saturated neighborhoods is
superfluous if the ambient manifold $M$ is compact, it is enough then
to ask that any leaf be orientation-preserving, as was shown by
Epstein \cite{E 72}. This is no longer true on non-compact manifolds,
as was shown  by Vogt \cite{V 89}. 

The initial definition of Seifert \cite{S 33} was slightly
different: 

a) He did not speak of ``foliation", but of ``fibration". 
 
b) He gave models for the possible neighborhoods of the leaves.  

In what concerns point a),  Seifert's definition is one of the
historical sources of the concept of fibration and fibre bundle. For
him a fibration is a decomposition of a manifold into ``fibres"; only
in a second phase can one try to construct the associated ``orbit
space", or the ``base" with our vocabulary. This shows that his
definition is closer to the present notion of foliation; in fact his
``fibration" is a foliation, but this can be seen only by 
using the required condition on model neighborhoods. We prefer to
speak about ``Seifert structure" and not ``Seifert fibration" precisely
because what is important to us is to see the structure as living
\textit{inside} the manifold, which makes possible to speak about
isotopies. For details about the historical development of different notions of
fibrations, see Zisman \cite{Z 99}. 

In what concerns point b), the possible orientable saturated
neighborhoods of foliations by circles coincide up to a
leaf-preserving diffeomorphism with Seifert's model neighborhoods. If
one drops the orientability condition, appears a new model which was not
considered by Seifert, but which is very useful in the classification
of \textit{non-orientable} 3-manifolds (see Scott \cite{S 83}, Bonahon \cite{B
  02}). Some general references about Seifert manifolds are Orlik
\cite{O 72}, Neumann \& Raymond \cite{NR 78} (where the base was
defined as an orbifold), Scott \cite{S 83}, Fomenko \& Matveev
\cite{FM 97} and Bonahon \cite{B 02}. 

In the sequel, we are interested in Seifert structures only \textit{up
  to isotopy}.  

\begin{definition} \label{isoseif}
  Two Seifert structures $\mathcal{F}_1$ and $\mathcal{F}_2$ on $M$
  are called \textbf{isotopic} if there exists $\phi \in
  \operatorname{Diff}^{\circ}(M)$ such that $\phi(\mathcal{F}_1)=
  \mathcal{F}_2$.  
\end{definition}

The following proposition is proved in Jaco \cite{J 80} and Fomenko \&
Matveev \cite{FM 97}. 

\begin{proposition} \label{seifbound}
  The only orientable compact connected 3-manifolds with non-empty
  boundary which admit more than one Seifert structure up to isotopy
  are the thick torus, the solid torus and the thick Klein bottle. 

 a) {\em If $M$ is a thick torus}, any essential curve on one of its
 boundary components is the fibre of a Seifert structure on $M$,
 unique up to isotopy, and devoid of exceptional fibres. Moreover, $M$
 appears like this as the total space of a trivial circle bundle over
 an annulus.

 b) {\em If $M$ is a solid torus} and $\gamma$ is a meridian of it, an
 essential curve $c$ on its boundary is 
 a fibre of a Seifert structure on $M$ if and only if their homological
 intersection 
 number $[c]\cdot [\gamma]$ (once they are arbitrarily oriented) is
 non-zero. In  this case, the 
 associated structure is unique up to isotopy and has at most one
 exceptional fibre. All fibres are regular if and only if $[c] \cdot
 [\gamma]=\pm 1$. In this last case, $M$ appears as the total space of
 a trivial  circle bundle over a disc. 

 c) {\em If $M$ is a thick Klein bottle}, it admits up to isotopy two
 Seifert structures. One of them is devoid of exceptional fibres and
 its space of orbits is a M{\"o}bius band. The other one has two
 exceptional fibres with holonomy of order $2$ and its space of orbits
 is topologically a disc. 
\end{proposition}

The closed orientable 3-manifolds which admit more than one Seifert
structure up to isotopy are also classified (see Bonahon \cite{B 02}
and the references therein). In this paper we need only the following
less general result, which can be deduced by combining \cite{B 02}
with \cite{N 81} (see Definition \ref{bound}):

\begin{proposition} \label{nonun}
  The only 3-manifolds which are diffeomorphic to abstract boundaries
  of normal surface singularities and which admit non-isotopic Seifert
  structures are the lens spaces.
\end{proposition}

\medskip
\subsection{Graph structures and JSJ decomposition theory} \label{graphJSJ}
$\:$
\medskip

If one glues various Seifert manifolds along components of their
boundaries, one obtains so-called \textit{graph-manifolds}:

\begin{definition} \label{graph}
   A \textbf{graph structure} on a 3-manifold $M$ is a pair 
   $(\mathcal{T}, \mathcal{F})$, where $\mathcal{T}$ is an embedded 
   surface in $M$ whose connected components are tori and where
   $\mathcal{F}$ is a Seifert structure on $M_{\mathcal{T}}$ (see
   Definition \ref{split}). We say
   that a graph structure is \textbf{orientable} if 
   $\mathcal{F}$ is an orientable Seifert structure on
   $M_{\mathcal{T}}$. If there exists a graph structure on $M$, we 
   say that $M$ is a \textbf{graph  manifold}. 
\end{definition}

Notice that no particular graph structure is specified when one speaks
about  a graph manifold. One only supposes that there exists one. In
the sequel we are interested in graph structures on a given manifold only
\textit{up to isotopy}:

\begin{definition} \label{isot}
  Two graph structures $(\mathcal{T}_1, \mathcal{F}_1)$, $(\mathcal{T}_2,
  \mathcal{F}_2)$ on $M$ are called \textbf{isotopic} if there exists
  $\phi\in \operatorname{Diff}^{\circ}(M)$ such that
  $\phi(\mathcal{T}_1)= \mathcal{T}_2$ 
  and $\phi_{\mathcal{T}_1}(\mathcal{F}_1)$ is isotopic to
  $\mathcal{F}_2$.
\end{definition}

Graph manifolds were introduced by Waldhausen \cite{W 67},
generalizing von Randow's \textit{tree manifolds} (see their
definition in the next paragraph) studied in \cite{R
  62}. Following Mumford \cite{M 61} who proved Poincar{\'e} conjecture
for the abstract boundaries of normal surface singularities (see
Definition \ref{bound}), von
Randow proved it for tree manifolds; his proof contained a gap which
was later filled by Scharf \cite{S 75}. 

Waldhausen's 
definition was different from Definition \ref{graph}. On one side
 he did not allow 
exceptional fibres in the Seifert structure on $M_{\mathcal{T}}$ and
on another side he did not fix (up to isotopy)  a precise fibration by
circles, but only supposed that such a fibration existed. He
represented 
a graph structure  by a finite graph with decorated vertices and edges
(corresponding respectively to the pieces of $M_{\mathcal{T}}$ and to
the components of $\mathcal{T}$), which explains the name. Tree
manifolds are then the graph manifolds which admit a graph
structure $(\mathcal{T}, \mathcal{F})$ such that the corresponding
graph is a tree and the base of 
the Seifert structure on $\mathcal{F}$ has genus $0$. With our
definition,  graph structures can
also be encoded by graphs. One has only to add more decorations to the
vertices, in order to keep in memory the exceptional fibres of the
corresponding Seifert fibred pieces.

With
his definition, Waldhausen solved the homeomorphism problem for
graph-manifolds, by giving normal forms for the graph structures on a
given manifold and by
showing that with exceptions in a finite explicit list, any
irreducible graph-manifold has a graph-structure in normal form which is unique
\textit{up to isotopy}. 

Later, Jaco \& Shalen \cite{JS 79} and
Johannson \cite{J 79} showed that there remains no exception in the
classification up to isotopy if one modifies the notion of
graph-structure by allowing exceptional fibres, that is, when one
works with Definition \ref{graph}. More generally, they proved:

\begin{theorem} \label{JSJthm}
  Let $M$ be a compact, connected, orientable and irreducible
  $3$-manifold (with possible non-empty boundary). Then $M$ contains
  an embedded surface $\mathcal{T}$ whose connected components are
  incompressible tori and such that any piece of $M_{\mathcal{T}}$ is
  either a Seifert manifold or is atoroidal. Moreover, if
  $\mathcal{T}$ is minimal for the inclusion among surfaces with this
  property, then it is well-defined up to isotopy.
\end{theorem}

We say that a minimal family $\mathcal{T}$ as in the previous
theorem is a 
 \textit{JSJ family of tori} in $M$. 

A variant
of the previous theorem considers also embedded annuli. These various theorems
of canonical decomposition  are 
called nowadays \textit{Jaco-Shalen-Johannson (JSJ) decomposition
  theory}, and were the 
starting point of Thurston's geometrization program, as well as of 
the theory of JSJ decompositions for groups. For details about
JSJ decompositions, in addition to the previously quoted books one can
consult Jaco \cite{J 80},  Neumann \& Swarup
\cite{NS 97},  Hatcher \cite{H 00} and Bonahon \cite{B 02}. In
\cite{PP 01} and \cite{PP 02}, we showed that also knot theory inside
an irreducible 3-manifold reflects the ambient JSJ decomposition. 

We define now a notion of \textit{minimality} for graph structures on
a given manifold. 

\begin{definition} \label{mingraph}
  Suppose that $(\mathcal{T}, \mathcal{F})$ is a graph structure on
  $M$. We say that it is \textbf{minimal} if the following conditions
  are verified:
    
 $\bullet$ No piece of $M_{\mathcal{T}}$ is a thick torus or a solid
 torus.

 $\bullet$ One cannot find a Seifert structure $\mathcal{F}'$ on
 $M_{\mathcal{T}}$ such that the images of its leaves by the
 reconstruction mapping 
 $r_{M, \mathcal{T}}$ coincide on a component of $\mathcal{T}$.
\end{definition}

As a corollary of Theorem \ref{JSJthm}, if $(\mathcal{T},
\mathcal{F})$ is a minimal graph structure on $M$, then $\mathcal{T}$
is the minimal JSJ system of tori in $M$. But one can prove more:

\begin{theorem} \label{minthm}
  Each closed orientable irreducible graph manifold which is not a
  torus fibration 
  with $|tr \: m|\geq 3$ admits a minimal graph 
  structure. Moreover, the family $\mathcal{T}$ of tori associated to
  a minimal graph structure coincides
  with the JSJ family of tori. In particular, it is unique up to an isotopy.
\end{theorem}
  
Suppose that $(\mathcal{T}, \mathcal{F})$ is a given graph structure
without thick tori and solid tori among its pieces. In view of
Proposition \ref{seifbound}, its only pieces which can have
non-isotopic Seifert structures are the thick Klein bottles. This
shows that, in
order to check whether $(\mathcal{T}, \mathcal{F})$ is minimal or not, one
has only to consider the possible choices of Seifert structures on
them up to isotopy (that is $2^n$ possibilities, where $n$ is the
number of such pieces).

Suppose that $M$ is a graph manifold which is neither a torus
fibration with $|tr \: m|\geq 3$, nor
a Seifert manifold which admits non-isotopic Seifert 
structures. Then, if $\mathcal{T}$ is a family of tori associated to a
minimal graph structure, there is a unique Seifert structure on
$M_{\mathcal{T}}$ up to isotopy, such that each piece which is a thick
Klein bottle has an orientable base.

\begin{definition}  \label{cangraph}
 Suppose that $M$ is an orientable graph manifold which is neither a torus
 fibration with $|tr \: m|\geq 3$ 
 nor a Seifert manifold which admits non-isotopic Seifert 
structures. We say that a minimal graph structure  is the
 \textbf{canonical graph structure} on $M$ if each piece which is a
 thick Klein bottle has an orientable base. 
\end{definition}

\medskip
\subsection{Plumbing structures} \label{plumbdef}
$\:$
\medskip

Plumbing structures are special types of graph structures:

\begin{definition} \label{plumbing}
  A \textbf{plumbing structure} on a 3-manifold $M$ is a graph
  structure without exceptional fibres $(\mathcal{T}, \mathcal{F})$ on
  $M$, such that for any component 
  $T$ of $\mathcal{T}$,  the homological intersection number on $T$
  of two fibres of $\mathcal{F}$ coming from opposite sides is equal
  to $\pm 1$. 
\end{definition}

Plumbing structures are the ancestors of graph structures. They were
introduced by Mumford \cite{M 61} in the study of singularities of
complex analytic surfaces (see Hirzebruch \cite{H 62},
Hirzebruch, Neumann \& Koh \cite{HNK 71}, as well as our
explanations in section \ref{plumcros}). In fact Mumford does not
speak about 
``plumbing structure". Instead, he describes a way to construct the
abstract boundary of a normal surface singularity (see Definition
\ref{bound}) by gluing total spaces of
circle-bundles over real surfaces using ``plumbing fixtures". 

Later on,
``plumbing" was more used as a verb than as a noun. That is, one
concentrated more on the operations needed to construct a new object
from elementary pieces, than on the structure obtained on the manifold
resulting from the construction. The fact that we are  interested
precisely in this structure up to isotopy and not on the graph which
encodes it, is a  difference  with Neumann \cite{N 81} for example. 

In \cite{N 81}, Neumann describes an algorithm for deciding if two
manifolds obtained by plumbing are diffeomorphic. He uses as an
essential ingredient Waldhausen's classification theorem of graph
manifolds (according to the definition which does not allow
exceptional fibres, see the comments made in  section \ref{graphJSJ}). In
fact, by using the 
uniqueness up to isotopy of the JSJ-tori, we can deduce the uniqueness
up to isotopy for special plumbing structures on singularity
boundaries. This is 
the subject of section \ref{canplumb}. 

Even if Definition \ref{plumbing} seems to suggest the opposite, the
class of graph manifolds is the same as the class of manifolds which
admit a plumbing structure. A way to see this is to use the
construction of plumbing structures on thick tori and solid tori
described in section \ref{HJplumb}. 
For a detailed comparison of graph structures and plumbing structures,
as well as for a study of the elementary operations on them, one can
consult Popescu-Pampu \cite[chapter 4]{PP 01}.

\subsection{Hirzebruch-Jung plumbing structures on thick tori and
  solid tori}
  \label{HJplumb} $\:$
\medskip

In this section we define special classes of plumbing structures on
thick tori and solid tori, which will be used in section 9. The
starting point is in both cases a pair 
$(L, \sigma)$ of a 2-dimensional lattice and a rational strictly
convex cone $\sigma \subset L_\R$, naturally attached to essential
curves on the boundary of the 3-manifold. 

 $\bullet$ \textit{Suppose first that $M$ is an \emph{oriented} thick
 torus}. 

On each  component of its boundary, we consider an essential  
 curve. Denote by $\gamma, \delta$ these curves. We
 suppose that their homology classes (once they are arbitrarily
 oriented) in $H_1(M, \R)\simeq \R^2$ are
 non-proportional. So, we are in presence of a 2-dimensional lattice
 $L= H_1(M, \Z)$ and of two distinct rational lines in it, generated
 by the homology classes $[\gamma], [\delta]$.

 Orient $\partial M$ compatibly with $M$. Then \textit{order} in an
 arbitrary way the components of $\partial M$: call the first one $T_-$ 
 and the second one $T_+$.  Denote by $\gamma_-$ the simple closed curve
 drawn on $T_-$ and by $\gamma_+$ the one drawn on $T_+$. Then
 \textit{orient} $\gamma_-$ and $\gamma_+$. By hypothesis, their
 homology classes 
 $[\gamma_-], [\gamma_+]$ are non-proportional primitive vectors in the
 2-dimensional lattice $L= H_1(M, \Z)$. This shows that  $([\gamma_-],
 [\gamma_+])$ is 
 a basis of $L_\R = H_1(M, \R)$ which induces an orientation of this
 vector space. As $T_+$ is a deformation retract of $M$, one has
 canonically $H_1(T_+, \Z)=L$, and so the ordered pair $(\gamma_-,
 \gamma_+)$ induces  an orientation of $T_+$. 

\begin{definition} \label{compat1}
 We say that $\gamma_-$ and $\gamma_+$ are oriented
 \textbf{compatibly with the orientation of } $M$ if, when taken in
 the order $( \gamma_-, \gamma_+)$, they induce on $T_+$ an
 orientation which coincides with 
 its orientation as a component of $\partial M$.
\end{definition}

Of course, a priori there is no reason for choosing this notion of
compatibility rather than the opposite one. Our choice was done in
order to get a more pleasant formulation for Lemma \ref{orientconv}.

Let $\sigma$ be the cone generated by $[\gamma_-]$ and
$[\gamma_+]$ in $L_\R$. As these homology classes were supposed
non-proportional, the cone $\sigma$ is strictly convex and has
non-empty interior. Denote by
$l_\pm$ the edge of $\sigma$ which contains the integral point
$[\gamma_\pm]$. Then, with the notations of section \ref{geomHJ},
$A_\pm= [\gamma_\pm]$.  Indeed, as $\gamma_{\pm}$ is an essential
curve of $T_{\pm}$, its homology class is a primitive vector of $L$. 

Let $(A_n)_{0 \leq n \leq r+1}$ be
the integral points on the compact edges of $P(\sigma)$, defined 
in section \ref{geomHJ}.  So, $OA_0= [\gamma_-]$ and $OA_{r+1}=
[\gamma_+]$. Let $(T_n)_{0 \leq n \leq r+2}$ be a sequence of pairwise
parallel
tori in $M$, such that $T_0=T_-$ and $T_{r+2}=T_+$. Moreover, we
number them in the order in which they appear between $T_-$ 
and $T_+$. Denote $\mathcal{T}:= \bigsqcup_{n=1}^{r+1} T_n$. If $M_n$
denotes the piece of $M_{\mathcal{T}}$ whose boundary components are
$T_n$ and $T_{n+1}$, where $n \in \{0,...,r+1\}$, we consider on it a
Seifert structure such that the homology class of its
fibres in $L$ is $OA_n$.

{\tt    \setlength{\unitlength}{0.92pt}}
\begin{figure}
  \epsfig{file=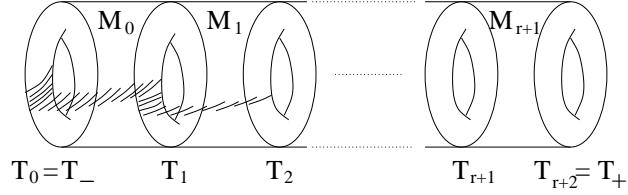, height= 25 mm} 
  \caption{Hirzebruch-Jung plumbing structures on thick tori}
\end{figure}

We get like this a plumbing structure on $M$, well-defined up to
isotopy, and depending only on the triple $(M, \gamma_-, \gamma_+)$. We
see that the simultaneous change of the orientations of $\gamma_-$ and
$\gamma_+$ or the change of their ordering (in order to respect the
compatibility condition of Definition \ref{compat1}) leads to the same
(unoriented) plumbing structure. 

\begin{definition} \label{HJPS1}
  We say that the previous unoriented plumbing structure on the oriented thick
  torus $M$ is \textbf{the Hirzebruch-Jung plumbing structure}
  associated to $(\gamma, \delta)$ and we denote it by
  $\mathcal{P}(M, \gamma, \delta)$.
\end{definition} 

\medskip
\pagebreak

$\bullet$ \textit{Suppose now that $M$ is an \emph{oriented} solid 
 torus}. 

{\tt    \setlength{\unitlength}{0.92pt}}
\begin{figure}
  \epsfig{file=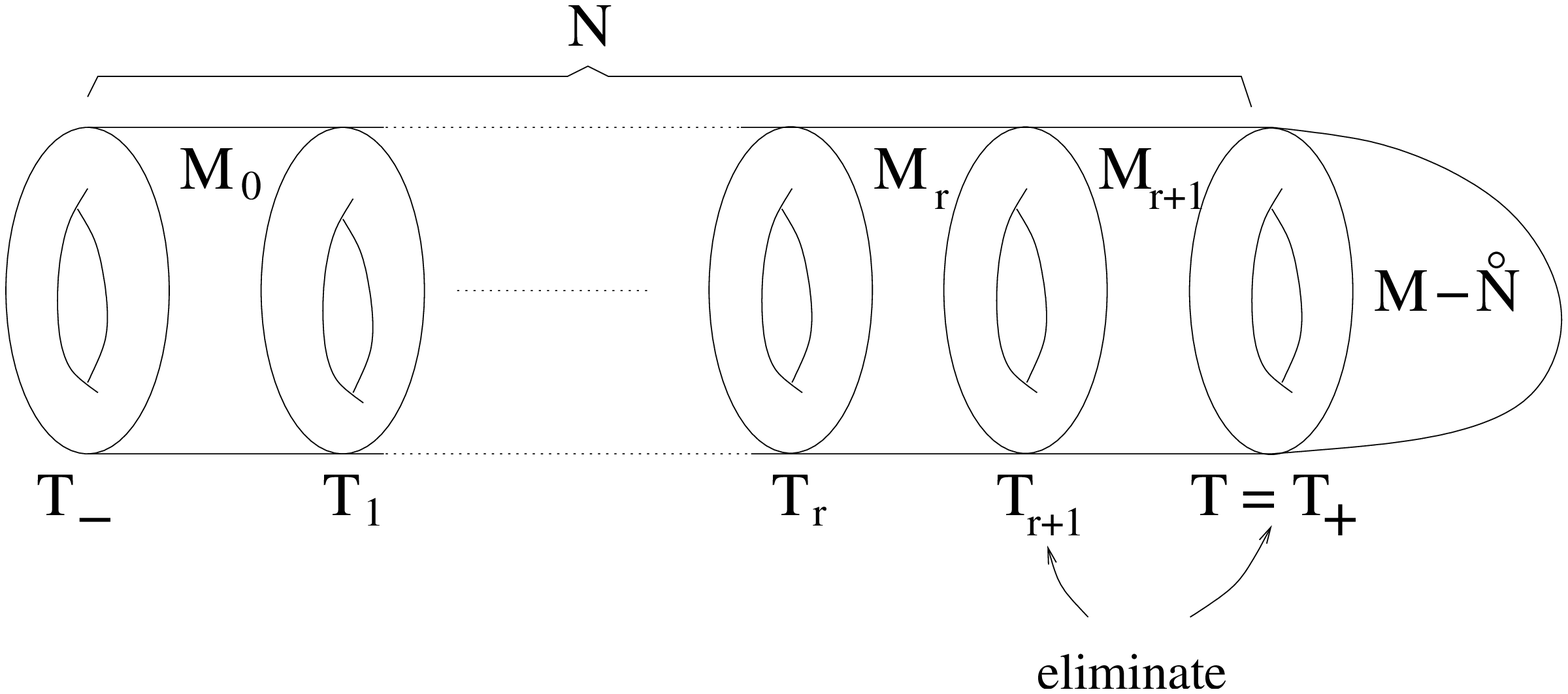, height= 35 mm} 
  \caption{Hirzebruch-Jung plumbing structures on solid tori}
\end{figure}

We consider an essential curve $\gamma$ on $\partial M$ \textit{which
  is not a meridian}. Take a torus
$T$ embedded in $\stackrel{\circ}{M}$ and parallel to 
$\partial M$. Denote by $N$ the thick torus contained between
$\partial M$ and $T$. Put $T_-= \partial M, \: T_+=T, \:
\gamma_-=\gamma$ and let $\gamma_+$ be an essential curve on $T_+$
which is a meridian of the solid torus
$M-\stackrel{\circ}{N}$ (see Figure 16). Consider the
Hirzebruch-Jung plumbing 
structure $\mathcal{P}(N, \gamma_-, \gamma_+)$. With the notations of
the construction done for thick tori, denote $\mathcal{T}(M, \gamma):=
\bigsqcup_{n=1}^r T_n$. Then the pieces of $M_{\mathcal{T}(M,
  \gamma)}$ are the thick tori $M_0, M_1,..., M_{r-1}$ and a solid
torus which is the ``union" of $M_r, M_{r+1}$ and
$M-\stackrel{\circ}{N}$. On $M_0,..., M_{r-1}$ we keep the Seifert
structure of $\mathcal{P}(N, \gamma_-, \gamma_+)$. On the solid torus
we extend the Seifert structure of $M_r$. By Proposition
\ref{seifbound} b), we see that this Seifert structure has no
exceptional fibres. This shows that we have constructed a
plumbing structure on $M$. It is obviously well-defined up to
isotopy, once the isotopy class of $\gamma$ is fixed. 

\begin{definition} \label{HJPS2}
  We say that the previous unoriented plumbing structure on the oriented solid
  torus $M$ is \textbf{the Hirzebruch-Jung plumbing structure}
  associated to $\gamma$ and we denote it by
  $\mathcal{P}(M, \gamma)$.
\end{definition} 
\medskip

\section{Generalities on the topology of surface singularities} \label{general}
 
In this section  we look at the boundaries $M(\mathcal{S})$ of
normal surface singularities $(\mathcal{S},0)$. We explain how to
associate to any normal crossings resolution $p$ of $(\mathcal{S},0)$
a plumbing structure on $M(\mathcal{S})$. Then we explain how to pass
from the plumbing structure associated to the minimal normal crossings
resolution of $(\mathcal{S},0)$ to the canonical graph structure on
$M(\mathcal{S})$ (see Definition \ref{cangraph}).

We recommend the survey articles of N{\'e}methi \cite{N 04} and Wall
\cite{W 00} for an introduction to the classification of normal
surface singularities. 
\medskip

\medskip
\subsection{Resolutions of normal surface singularities and their dual
  graphs} \label{resdualgr} 
$\:$
\medskip

First we recall basic facts about normal analytic spaces. 
Let $\mathcal{V}$ be a reduced analytic space. It is called
\textit{normal} if for any point $P\in \mathcal{V}$, the germ
$(\mathcal{V},P)$ is irreducible and its local algebra is integrally
closed in its field of fractions. If $\mathcal{V}$ is not normal, then
there exists a finite map $\nu:\tilde{\mathcal{V}}\rightarrow
\mathcal{V}$ which is an isomorphism over a dense open set of
$\mathcal{V}$ and such that $\tilde{\mathcal{V}}$ is normal. Such a
map, which is unique up to unique isomorphism, is called a
\textit{normalization} map of $\mathcal{V}$. 

A reduced analytic curve is normal if and only if it is smooth. If a
germ $(\mathcal{S},0)$ of reduced surface is normal, then there exists
a representative of it, which we keep calling $\mathcal{S}$, such that
$\mathcal{S}-0$ is smooth. The converse is not true. 

\medskip

Let $(\mathcal{S},0)$ be a germ of normal complex analytic
surface. We say also that $(\mathcal{S},0)$ is a \textit{normal surface
singularity} (even if the point 0 is regular on $\mathcal{S}$). In the
sequel, we use the same notation $(\mathcal{S},0)$ for the germ and
for a
sufficiently small representative of it. If $e: (\mathcal{S},0) \rightarrow
(\C^N,0)$ is any local embedding, denote by $\mathcal{S}_{e,r}$ the
intersection of $\mathcal{S}$ with a euclidean ball of $\C^N$ of
radius $r \ll 1$ and by $M_{e,r}(\mathcal{S})$ the boundary of
$\mathcal{S}_{e,r}$. 

By general transversality theorems due to Whitney, when $r >0$ is
small enough, $M_{e,r}(\mathcal{S})$ is a smooth manifold, {\em naturally
oriented} as the boundary of the complex manifold
$\mathcal{S}_{e,r}$. It does not depend on the choices of embedding
$e$ and radius $r \ll 1$ made to define it (see Durfee \cite{D 83}). 

\begin{definition} \label{bound} 
 An oriented 3-manifold $M(\mathcal{S})$ orientation-preserving
 diffeomorphic with the manifolds $M_{e,r}(\mathcal{S})$, where $r >0$ is
 small enough, is called the \textbf{(abstract) boundary} or the
 \textbf{link} of the singularity $(\mathcal{S},0)$.
\end{definition}

It is important to keep in mind that {\em in the sequel $M(\mathcal{S})$
is supposed naturally oriented as explained before}. In
order to understand better this remark, look at
Theorem \ref{changeor}.
\medskip

The easiest way to describe the topological type of the manifold
$M(\mathcal{S})$ is (as first done by Mumford \cite{M 61}) by
retracting it to the exceptional divisor of a 
resolution of $(\mathcal{S},0)$. Let us first define this last notion.

\begin{definition} \label{resdef}
An analytic map  $p: (\mathcal{R}, E) \rightarrow (\mathcal{S},0)$ is
called \textbf{a resolution} of the singularity $(\mathcal{S},0)$ with
\textbf{exceptional divisor} $E= p^{-1}(0)$ if the
following conditions are simultaneously satisfied: 

$\bullet$ $\: \mathcal{R}$ is a smooth surface;

$\bullet$ $\: p$ is a proper morphism;

$\bullet$ the restriction of $p$ to $\mathcal{R} -E=\mathcal{R} -f^{-1}(0)$ is
an isomorphism onto $\mathcal{S}-0$.

We say that $p: (\mathcal{R}, E) \rightarrow (\mathcal{S},0)$ is \textbf{a
  normal crossings resolution} if one has moreover: 

$\bullet$ $E$ is a divisor with normal crossings.
\end{definition}

Recall that, by definition, a divisor on a smooth complex surface has
  \textit{normal crossings} if in the neighborhood of any of its
  points, its support is either smooth,  or the union of transverse
  smooth curves.

Normal crossings resolutions always exist (see Laufer \cite{L 71} and
Lipman \cite{L 75} for a careful presentation of the Hirzebruch-Jung
method of resolution, as well as Cossart \cite{C 00} for Zariski's
method of resolution by normalized blow-ups). 

There is a unique 
\textit{minimal} resolution, which we  denote 
$p_{min}:(\mathcal{R}_{min}, E_{min})\rightarrow (\mathcal{S},0).$  
The minimality property means that any other resolution
$p:(\mathcal{R},E)\rightarrow (\mathcal{S},0)$ can be factorized as
$p=p_{min}\circ q$, where $q:\mathcal{R}\rightarrow \mathcal{R}_{min}$
is a proper bimeromorphic map. 
The minimal resolution $p_{min}$ is characterized by the fact 
that $E_{min}$ contains no component $E_i$ which is smooth, rational and of
self-intersection $-1$ (classically called an \textit{exceptional
  curve of the first kind}).  

Analogously, there is a unique resolution which is minimal among
normal crossings ones. We denote it: 
$$p_{mnc}:(\mathcal{R}_{mnc}, E_{mnc})\rightarrow (\mathcal{S},0)$$ 
It is characterized by the fact that $E_{mnc}$ has normal crossings
and each component $E_i$ of $E_{mnc}$ which is
 an exceptional
  curve of the first kind contains at least 3 points which are
singular on $E_{mnc}$.

If a
normal crossings resolution has moreover only smooth components, one
says usually that the resolution is \textit{good}; there exists also a
unique minimal good resolution, but in this paper we don't consider it. 

The following criterion allows one to recognize the divisors which are
exceptional with respect to some resolution of a normal surface
singularity.

\begin{theorem} \label{critcontr} 
  Let $E$ be a reduced compact connected divisor in a smooth surface
  $\mathcal{R}$. Denote by $(E_i)_{1 \leq i \leq n}$ its components. Then $E$
  is the exceptional divisor of a resolution of a normal surface
  singularity if and only if the intersection matrix $(E_i \cdot
  E_j)_{i,j}$ is negative definite.
\end{theorem}

The necessity is classical (see \cite[section 9]{HNK 71}, where is
presented Mumford's proof of \cite{M 61} and where the oldest
reference is to Du Val \cite{V 44}). The sufficiency was proved by
Grauert \cite{G 62} (see also Laufer \cite{L 71}). If $E$ verifies the
conditions which are stated to be equivalent in the theorem, one also
says that $E$ \textit{can be contracted} on $\mathcal{R}$.

\medskip

\textit{From now on we suppose that $p:(\mathcal{R},E)\rightarrow
  (\mathcal{S},0)$ is a normal crossings resolution of 
$(\mathcal{S},0)$.} 

Denote 
by  $\Gamma(p)$  its \textit{weighted dual graph}. Its set of vertices
$\mathcal{V}(p)$ is in 
bijection with the irreducible components of $E$. Depending on the
context, we  think about $E_i$ as a curve on $\mathcal{R}$ or a
vertex of $\Gamma(p)$. The vertices
which represent the components $E_i$ and $E_j$ are joined by as many
edges as $E_i$ and $E_j$ have intersection points on $\mathcal{R}$. In
particular, 
there are as many \textit{loops} based at the vertex $E_i$ as singular
points (that is, self-intersections) on the curve $E_i$ (see Figure
17). Each vertex 
$E_i$ is decorated by two weights, the \textit{geometric genus} $g_i$
of the curve $E_i$ 
(that is, the genus of its normalization) and its self-intersection
number $e_i \leq -1$ in $\mathcal{R}$. Denote 
also by $\delta_i$ the \textit{valency} of the vertex $E_i$, that is, 
the number of edges starting from it (where each loop counts
for 2). For example, in Figure 17 one has $\delta_1=9, \delta_2=5,$ etc.

{\tt    \setlength{\unitlength}{0.92pt}}
\begin{figure}
  \epsfig{file=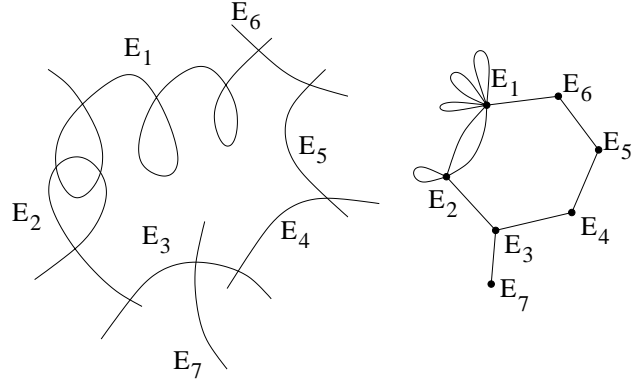, height= 50 mm} 
  \caption{A normal crossings divisor and its dual graph}
\end{figure}
\medskip


\subsection{The plumbing structure associated to a normal crossings
  resolution} \label{plumcros}
$\:$
\medskip

By Definition \ref{bound}, $M(\mathcal{S})$ is diffeomorphic to
$M_{e,r}(\mathcal{S})$, where $e: (\mathcal{S},0) \rightarrow
(\C^N,0)$ is an embedding and $r \ll 1$. But $M_{e,r}(\mathcal{S})$ is
the level-set at level $r$ of the function $\rho_e : (\mathcal{S},0)
\rightarrow 
(\R_+,0)$, the restriction to $e(\mathcal{S})$ of the distance-function
to the origin in $\C^N$. 

As the resolution $p$ realizes by definition an isomorphism between
$\mathcal{R}- E$ 
and $\mathcal{S}-0$, it means that
$M_{e,r}(\mathcal{S})=\rho_e^{-1}(r)$ is diffeomorphic to
$\psi_e^{-1}(r)$, where $\psi_e:= \rho_e \circ p$. The advantage of
this changed viewpoint on $M(\mathcal{S})$ is that it appears now
orientation-preserving diffeomorphic to the boundary of a ``tubular
neighborhood" of the curve $E$ in the smooth manifold $\mathcal{R}$. 
As in general $E$ has singularities, one has to discuss the precise
meaning of the notion of tubular neighborhood. We quote Mumford \cite[pages
230-231]{M 61}:

\begin{quote}

{\small Now the general problem, given a complex $K \subset E^n$,
  Euclidean $n$-space, to define a tubular neighborhood, has been
  attacked by topologists in several ways although it does not appear
  to have been treated definitively as yet. J.H.C. Whitehead \cite{W
  39}, when $K$ is a subcomplex in a triangulation of $E^n$, has
  defined it as the boundary of the star of $K$ in the second
  barycentric subdivision of the given triangulation. I am informed
  that Thom \cite{T 59} has considered it more from our point of view:
  for a suitably restricted class of positive $C^{\infty}$ fcns. $f$
  such that $f(P)=0$ if and only if $P \in K$, define the tubular
  neighborhood of $K$ to be the level manifolds $f=\epsilon$, small
  $\epsilon$. The catch is how to suitably restrict $f$; here the
  archtype for $f^{-1}$ may be thought of as the potential
  distribution due to a uniform charge on $K$.}
\end{quote}

Let us come back to the normal crossings divisor $E$ in the smooth
surface $\mathcal{R}$.

{\em If $E$ is smooth}, then one can construct a diffeomorphism between a
tubular neighborhood $U(E)$ of $E$ in $\mathcal{R}$ and of $E$ in the
total space $N_{\mathcal{R}}E$ of its normal bundle in
$\mathcal{R}$. As $N_{\mathcal{R}}E$ is naturally fibred by discs,
this is also true for $U(E)$. The fibration of $U(E)$ can be chosen in such a
way that the levels $\psi_e^{-1}(r)$ are transversal to the fibres for
$r \ll 1$. In this way one gets a Seifert structure without singular
fibres on $\psi_e^{-1}(r)\simeq M(\mathcal{S})$.

Suppose now that $E$ {\em is not smooth, but that its irreducible
components are so}. One can also define in this situation a notion of
tubular neighborhood $U(E)$ of $E$ in $\mathcal{R}$. One way to do it
is to take the union of conveniently chosen tubular neighborhoods
$U(E_i)$ of $E$'s components $E_i$. Abstractly, one has to glue 
the 4-manifolds with boundary $U(E_i)$ by identifying
well-chosen neighborhoods of the points which get identified on
$E$. This procedure is what is called the ``plumbing" of disc-bundles
over surfaces (see Hirzebruch \cite{H 62}, Hirzebruch \& Neumann \&
Koh \cite{HNK 71}, Brieskorn \cite{B 00}). Its effect on the
boundaries $\partial U(E_i)$ is to take out saturated filled tori and
to identify their boundaries, by a diffeomorphism which permutes
fibres and meridians in an orientation-preserving way. This is the
3-dimensional ``plumbing" operation introduced by Mumford \cite{M 61},
alluded to in section \ref{plumbdef}.

In order to understand what happens near a singular point of $E$, it
is convenient to choose local coordinates $(x,y)$ on $E$ in the
neighborhood of the singular point, such that $E$
is defined by the equation $xy=0$. So, $y=0$ defines locally an
irreducible component $E_i$ of $E$ and similarly $x=0$ defines
$E_j$. It is possible that $E_i=E_j$, a situation excluded in the
previous paragraph for pedagogical reasons. If this equality is true,
then the same plumbing procedure can be applied, this time by
identifying well-chosen neighborhoods of points of the same 4-manifold
with boundary $U(E_i)$. 

\textit{At this point appears a subtlety}: the 4-manifold $U(E_i)$ to be
considered is no longer a tubular neighborhood of $E_i$ in
$\mathcal{R}$, but instead of the normalization $\tilde{E}_i$ of $E_i$
inside the modified normal bundle
$\nu_i^*T\mathcal{R}/T\tilde{E}_i$. Here $\nu_i:\tilde{E}_i
\rightarrow \mathcal{R}$ denotes the normalization map of $E_i$ and
$T\mathcal{R}$, respectively $T\tilde{E}_i$ denote the holomorphic
tangent bundles to the smooth complex manifolds $\mathcal{R}$ and
$\tilde{E}_i$. As a real differentiable bundle of rank 2, this vector
bundle over $ \tilde{E}_i$  is characterized by its Euler number
$\tilde{e}_i$, which is equal to the self-intersection number of
$\tilde{E}_i$ inside the total space of the bundle. This number is
related to the self-intersection of $E_i$ inside $\mathcal{R}$ in the
following way (see Neumann \cite[page 333]{N 81}):

\begin{lemma} \label{compeul}
  If $\tilde{e}_i$ is the Euler number of the real bundle
  $\nu_i^*T\mathcal{R}/T\tilde{E}_i$ over $ \tilde{E}_i$, where
  $\nu_i:\tilde{E}_i \rightarrow \mathcal{R}$ is the normalization map
  of $E_i$, then $\tilde{e}_i= e_i-\delta_i$. 
\end{lemma}

\textbf{Proof:} In order to understand this formula, just think at the
effect of a 
small isotopy of $E_i$ inside $\mathcal{R}$. Near each self-crossing
point of $E_i$, the intersection point of one branch of $E_i$ with the
image of the other branch after the isotopy is not counted when one
computes $\tilde{e}_i$. \hfill $\Box$
\medskip

Notice that Theorem \ref{critcontr} is {\em true} if one takes as diagonal
entries of the matrix the numbers $e_i=E_i^2$, but is {\em false}  if one
takes instead the numbers $\tilde{e}_i$. The easiest example is given
by an irreducible divisor $E=E_1$, with $e_1=1 >0$ and $\delta_1=2$
which, by Lemma \ref{compeul} implies that $\tilde{e}_1 =-1 <0$.

In Figure 18 we represent in two ways the local situation near
the chosen singular point of $E$. On the left we simply draw the union
of the two neighborhoods $U(E_i)$ and $U(E_j)$. On the right, ``the
corners are smoothed". This is precisely what happens when we look at
the levels of the function $\psi_e$. Moreover, we represent by
interrupted lines the real analytic set defined by the equation
$|x|=|y|$. Its intersection with $\psi_e^{-1}(r)\simeq  \partial U(E)
\simeq M(\mathcal{S})$ 
is a two-dimensional torus $T$. {\em This is the way in which such tori
appear naturally} as structural elements of the 3-manifolds
$M(\mathcal{S})$. One also sees how the complement of $T$ in $\partial
U(E)$ is fibred by boundaries of discs transversal to $E_i$ or $E_j$. 

By considering model neighborhoods of the singular points of $E$
structured as in
the right-hand side of Figure 18 and conveniently extending them to a
tubular neighborhood of all of $E$, one gets a retraction 
  $$\Phi: U(E) \rightarrow E$$
which restricts to a locally trivial disc-fibration over the smooth
locus of $E$ and whose fibre over each singular point of $E$ is a cone
over a 2-dimensional torus. By considering the restriction
$\Phi|_{\partial U(E)}$, 
we see that the fibres over the singular points of $E$ are embedded
tori, and that their complement gets fibred by circles. 

As $\partial U(E)$ is orientation-preserving diffeomorphic to
$M(\mathcal{S})$, we see that $M(\mathcal{S})$ gets endowed with a 
graph structure $(\mathcal{T}(p), \mathcal{F}(p))$ well-defined up to isotopy. 
It is a good test of the understanding of the complexifications of
Figure 18 to show that $(\mathcal{T}(p), \mathcal{F}(p))$ \textit{is
  in fact a plumbing structure} (see Definition \ref{plumbing}).

{\tt    \setlength{\unitlength}{0.92pt}}
\begin{figure}
  \epsfig{file=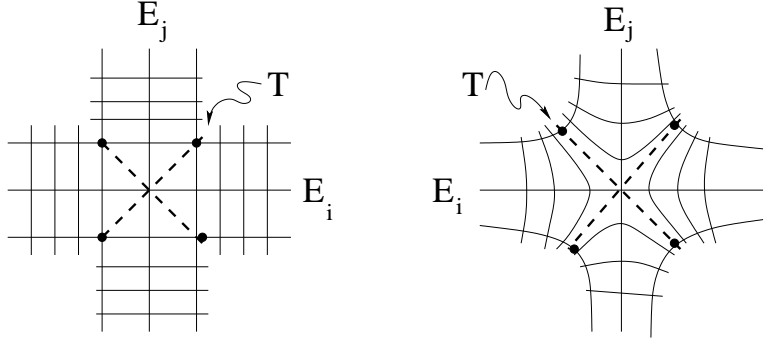, height= 45 mm} 
  \caption{The local configuration which leads to plumbing}
\end{figure}

The pieces of
$M(\mathcal{S})_{ \mathcal{T}(p)}$ correspond to the irreducible
components of $E$, that is to the vertices of $\Gamma(p)$. Denote by
$M(E_i)$ the piece which corresponds to $E_i$. The fibres 
of $M(E_i)$ are obtained up to isotopy by
cutting the boundary of the 
chosen sufficiently small tubular neighborhood of $E$ with  smooth
holomorphic curves transversal to $E$ at smooth points of $E_i$. So,
the plumbing structure $(\mathcal{T}(p),
\mathcal{F}(p))$ is naturally oriented. 

\begin{lemma} \label{orientconv}
  With their natural orientations, the fibres on both sides of any
  component of $\mathcal{T}(p)$ are oriented compatibly with the
  orientation of $M(\mathcal{S})$.
\end{lemma}

\textbf{Proof:} The notion of compatibility we speak about is the one
of Definition \ref{compat1}. We mean that, if we take an arbitrary
component $T$ of $\mathcal{T}(p)$, and a tubular neighborhood
$N(T)$ such that its preimage in $M(\mathcal{S})_{\mathcal{T}(p)}$
is saturated by the leaves of the foliation $\mathcal{F}(p)$, then two 
fibres, one in each boundary component of $N(T)$, are oriented
compatibly with the orientation of $N(T)$. Now, this is an
instructive exercise on the geometrical understanding of the relations
between the orientations of various objects in the neighborhood of a
normal crossing on a smooth surface. Just think of the
complexification of Figure 18. \hfill $\Box$

\begin{corollary} \label{simult}
  The orientation of the fibres of $(\mathcal{T}(p),
\mathcal{F}(p))$ is determined by the associated unoriented plumbing
structure up to a simultaneous change of orientation of all the fibres.
\end{corollary}

\textbf{Proof:} Consider the unoriented plumbing structure. Start from
an arbitrary piece $M(E_i)$, and choose one of the two continuous
orientations of its fibres. Then propagate this orientations farther
and farther through the components of $\mathcal{T}(p)$, by respecting
the compatibility condition on the neighboring orientations. As
$M(\mathcal{S})$ is connected, we know that after a finite number of
steps one has oriented the fibres of all the pieces. As one
orientation exists which is compatible in the neighborhood of all the
tori, we see that our process cannot arrive at a contradiction
(that is, a non-trivial monodromy around
a loop of $\Gamma(p)$ in the choice of orientations).  
\hfill $\Box$
\medskip

The following lemma is a particular case of the study done in Mumford
\cite[page 11]{M 61} and Hirzebruch \cite[page 250-03]{H 62}.

\begin{lemma} \label{reltor}
  Suppose that $E_i$ is a component of $E$ which is smooth, rational and
  whose valency in the graph $\Gamma(p)$ is $2$. In the thick torus
  $M(E_i)$ which corresponds to it in the plumbing structure
  $(\mathcal{T}(p), \mathcal{F}(p))$, consider an oriented fibre $f$
  of $M(E_i)$, as well as oriented fibres $f'$, $f''$ of the two
  (possibly coinciding) adjacent pieces. Then one has the following
  relation in the homology group $H_1(M(E_i), \Z)$:
    $$[f'] + [f''] = |e_i|\cdot [f].$$ 
\end{lemma}

\medskip

\subsection{The topological characterization of HJ and cusp singularities}
\label{topchar} $\:$
\medskip

We want now to understand how to pass from the plumbing structure
$(\mathcal{T}(p), \mathcal{F}(p))$ on $M(\mathcal{S})$ to the
canonical graph structure on it (see Definition \ref{cangraph}). We
see that the pieces of $M(\mathcal{S})_{\mathcal{T}(p)}$ which are
thick tori correspond to components $E_i$ which are smooth and
rational with $\delta_i=2$, and those which are solid tori correspond
to components $E_i$ which are smooth and rational with
$\delta_i=1$. It is then natural to introduce the following:

\begin{definition} \label{road}
  We say that a vertex $E_i$ of $\Gamma(p)$ is a \textbf{chain vertex}
  if $E_i$ is smooth, $g_i=0$ and $\delta_i \leq 2$. If moreover
  $\delta_i =2$, we call it an \textbf{interior chain vertex},
  otherwise we call it a \textbf{terminal chain vertex}. We say that a
  vertex of $\Gamma(p)$ is a \textbf{node} if it is not a chain
  vertex. 
\end{definition}

In \cite{LMW 89}, L{\^e}, Michel \& Weber used the name ``rupture vertex"
for a node in the dual graph associated to the minimal embedded
resolution of a plane curve singularity. In their situation, where all
the vertices represent smooth rational curves, nodes are
simply those of valency $\geq 3$. In our case this is no longer true,
as one can have also vertices of valency $\leq 2$, if they correspond
to curves $E_i$ which are either not smooth or of genus $g_i \geq 1$. 

Denote by $\mathcal{N}(p)$ the set of nodes of
$\Gamma(p)$. It is an empty set if and only if $\Gamma(p)$ is
topologically a segment or a circle and all the components $E_i$ are
smooth rational curves. The first situation occurs precisely for the
Hirzebruch-Jung singularities, defined in Section  \ref{torsurf} (see
Proposition \ref{resmintor}), and
the second one for cusp singularities,  introduced by Hirzebruch
\cite{H 73} in the number-theoretical context of the study of Hilbert
modular surfaces.

\begin{definition} \label{cusp}
  A germ $(\mathcal{S},0)$ of normal surface singularity is called a
  \textbf{cusp singularity} if it has a
resolution $p$ such that $\Gamma(p)$ is topologically a circle and $
   \mathcal{N}(p)=\emptyset$.
\end{definition}

For other definitions and details about them, see  Hirzebruch \cite{H
   73}, Laufer \cite{L 77} (where they appear as special cases of
   \textit{minimally elliptic singularities}), Ebeling \& Wall
   \cite{EW 85} (where they appear as special 
   cases of \textit{Kodaira singularities}), Oda \cite{O 88}, Wall
   \cite{W 00} and N{\'e}methi  \cite{N 04}. They were generalized to
   higher dimensions by Tsuchihashi (see Oda \cite[Chapter 4]{O 88}).

In the previous definition it is not possible to replace the
resolution $p$ by the minimal normal crossings one. Indeed:

\begin{lemma} \label{whynot}
  If $(\mathcal{S},0)$ is a cusp singularity, then 
  $\Gamma(p_{mnc})$ is topologically a circle and either
  $\mathcal{N}(p_{mnc})=\emptyset$, or $E_{mnc}$ is irreducible,
  rational, with one singular point where it has normal crossings.
\end{lemma}

\textbf{Proof:} One passes from $p$ to $p_{mnc}$ by successively
contracting components $F$ which are smooth, rational and verify
$F^2=-1$ (that is, exceptional curves of the first kind, by a remark
which follows Definition \ref{resdef}). The new exceptional divisor
verifies the same hypothesis as the one of $p$, except when one passes from a
divisor with 2 components to a divisor with one component. In this
last situation, this second irreducible divisor is rational, as its
strict transform $F$ is so. Moreover, it has one singular point with normal
crossing branches passing through it, as by hypothesis $F$ cuts
transversely the other component of the first divisor in exactly two
points. \hfill $\Box$
\medskip

We would like to emphasize the following theorem due to
Neumann \cite[Theorem 3]{N 81}, which characterizes Hirzebruch-Jung
and cusp singularities among normal surface singularities. 

\begin{theorem} \label{changeor}
  Let $(\mathcal{S},0)$ be a normal surface singularity. 
  The manifold $-M(\mathcal{S})$ is orientation-preserving
  diffeomorphic to the abstract boundary of a normal surface
  singularity if and only if $(\mathcal{S}, 0)$ is either a
  Hirzebruch-Jung singularity or a cusp-singularity. 
\end{theorem}

Recall that $-M(\mathcal{S})$ denotes the manifold $M(\mathcal{S})$
with reversed orientation. 

We will bring more light on this theorem with Propositions
\ref{orientchange} and \ref{dualcusp}, which show that for both
Hirzebruch-Jung and cusp singularities, the involutions
$M(\mathcal{S})\rightsquigarrow -M(\mathcal{S})$ are manifestations of the
duality described in section \ref{geomcomp}. 

As Hirzebruch-Jung singularities, cusp singularities
can also be defined using toric geometry (see Oda \cite[Chapter 4]{O 88}). 
In the same spirit, as a particular case of Laufer's \cite{L 73}
classification of \textit{taut} singularities, we have:

\begin{theorem} \label{taut}
  Hirzebruch-Jung and cusp singularities are taut, that is, their
  analytical type is determined by their topological type.
\end{theorem}

For this reason, it is natural to ask which 3-manifolds are obtained
as abstract boundaries of Hirzebruch-Jung singularities and cusp
singularities. This question is answered by:

\begin{proposition} \label{HJcusp}
  1) $(\mathcal{S},0)$ is a Hirzebruch-Jung singularity if and only if 
 $M(\mathcal{S})$ is a lens space.
Moreover, each oriented lens space appears like this. 

2) $(\mathcal{S},0)$ is a cusp singularity if and only if 
 $M(\mathcal{S})$ is a torus fibration with algebraic
 monodromy of trace $\geq 3$. 
Moreover, each oriented torus fibration of this type appears like this. 
\end{proposition}

\textbf{Proof:} This proposition is a particular case of  Neumann
\cite[Corollary 8.3]{N 81}. Here we sketch the proofs of
the necessities, in order to develop tools for sections \ref{caselens}
and \ref{torfibr}. 

Let $p:(\mathcal{R},E) \rightarrow (\mathcal{S},0)$ be the minimal
normal crossings resolution of $(\mathcal{S},0)$ (for notational
convenience, we drop the index ``mnc"). Denote by $U(E)$ a (closed)
tubular neighborhood of $E$ in $\mathcal{R}$ and by $\Phi: U(E)
\rightarrow E$ a preferred retraction, as defined in section
\ref{plumcros}. Denote also by $$\Psi: \partial U(E) \rightarrow E$$ the
restriction of $\Phi$ to $\partial U(E) \simeq M(\mathcal{S})$. 
\medskip

1) {\em Suppose that $(\mathcal{S},0)$ is a Hirzebruch-Jung
   singularity.} 

Orient the segment $\Gamma(p)$. Denote then by
   $E_1,...,E_r$ the components of $E$ in the order in which they
   appear along $\Gamma(p)$ in the positive direction. For each $i \in
   \{1,...,r-1\}$, denote by $A_{i, i+1}$ the intersection point of
   $E_i$ and $E_{i+1}$. Consider also two other points $A_{0,1} \in
   E_1, \: A_{r,r+1}\in E_r$ which are smooth points of $E$. Then consider on
   each component $E_i$ a Morse function 
     $$\Pi_i : E_i \rightarrow [\frac{i-1}{r}, \frac{i}{r}]$$
having as its only critical points $A_{i-1,i}$ (where $\Pi_i$ attains
its  minimum) and $A_{i, i+1}$ (where $\Pi_i$ attains
its maximum). As $\Pi_i(A_{i,i+1})=
\Pi_{i+1}(A_{i,i+1})$ for all $i \in \{1,...,r-1\}$, we see that the
maps $\Pi_i$ can be glued together in a continuous map 
   $$\Pi: E \rightarrow [0,1].$$

Consider the composed continuous map $\Pi \circ \Psi : M(\mathcal{S})
\rightarrow [0,1]$ (see Figure 19).

{\tt    \setlength{\unitlength}{0.92pt}}
\begin{figure}
  \epsfig{file=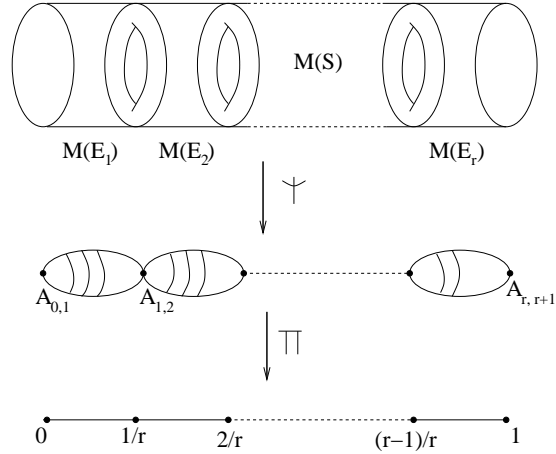, height= 60 mm} 
  \caption{The maps $\Pi$ and $\Psi$ for a Hirzebruch-Jung singularity}
\end{figure}

Our construction shows that its fibres over $0$ and $1$ are circles
and that those over interior points of $[0,1]$ are tori. Moreover,
each such torus splits $M$ into two solid tori. By Definition
\ref{lenstor}, we see that \textit{$M$ is a lens space}. 
\medskip

It remains now to prove that each oriented lens space appears like
this.

Denote $L:= H_1(M(\mathcal{S})- (\Pi \circ \Psi)^{-1}\{0,1\}, \Z)$. As
$M(\mathcal{S})- (\Pi \circ \Psi)^{-1}\{0,1\}$ is the interior of a
thick torus foliated by the tori $(\Pi \circ \Psi)^{-1}(c)$, where $c
\in (0,1)$, we 
see that $L$ is a 2-dimensional lattice. With the notations of section
\ref{plumcros}, let $f_i$ be an oriented fibre in the piece $M(E_i)$ of the
plumbing structure $(\mathcal{T}(p), \mathcal{F}(p))$ which
corresponds to $E_i$. Consider also $f_0$ and $f_{r+1}$, canonically
oriented meridians on the boundaries of tubular neighborhoods of $(\Pi
\circ \Psi)^{-1}(0)$, respectively $(\Pi \circ \Psi)^{-1}(1)$. 

For each $i \in \{0,...,r+1\}$, denote by $v_i := [f_i]\in L$ the
homology class of $f_i$. Recall that $e_i:= E_i^2$. By Lemma
\ref{reltor}, we see that 
\begin{equation} \label{relens} 
  v_{i+1}= 
  \mid e_i \mid \cdot v_i - v_{i-1}, \: \forall \: i \in \{ 0,...,r\}. 
\end{equation}
By Proposition \ref{resmintor}, $p$ is also the minimal resolution of
$(\mathcal{S},0)$, which shows that $|e_i|\geq 2, \: \forall i \in
\{1,...,r\}$. Now apply Proposition \ref{converse}. We deduce that the
numbers $e_i$ are determined by the oriented topological type of the
lens space $M(\mathcal{S})$, once the isotopy class of the tori $(\Pi
\circ \Psi)^{-1}(c)$ is fixed. 

This shows that, starting from any oriented lens space $M$ and torus
$T\subset M$ which splits $M$ into two solid tori, one can construct a
Hirzebruch-Jung singularity $(\mathcal{S},0)$ such that
$M(\mathcal{S})\simeq M$ only by looking at the classes of the
meridians of the two solid tori in the lattice $L=H_1(T, \Z)$. One has
only to be careful to orient them compatibly with the orientation of
$M$ (as explained at the beginning of the proof of Lemma
\ref{orientconv}). 
\medskip

2) {\em Suppose that $(\mathcal{S},0)$ is a cusp singularity.} 

$\bullet$ {\em Consider first the case where $r \geq 2$}. Orient
   the circle $\Gamma(p)$ and choose one of its vertices. Denote then
   by $E_1,...,E_r$ the components of $E$ in the order in which they
   appear along $\Gamma(p)$ in the positive direction, starting from
   $E_1$. For each $i \in \{1,...,r\}$, denote by $A_{i, i+1}$ the
   intersection point of $E_i$ and $E_{i+1}$, where $E_{r+1}=
   E_1$. Consider then functions $\Pi_i : E_i \rightarrow
   [\frac{i-1}{r}, \frac{i}{r}]$  with the same properties as in the
   case of Hirzebruch-Jung singularities. By passing to the quotient
   $\R \rightarrow \R/\Z$, we can glue the previous maps into a
   continuous map:
      $$\Pi: E \rightarrow \R/\Z.$$
Consider then the map $\Pi \circ \Psi : M(\mathcal{S})\rightarrow
\R/\Z$ (see Figure 20).

{\tt    \setlength{\unitlength}{0.92pt}}
\begin{figure}
  \epsfig{file=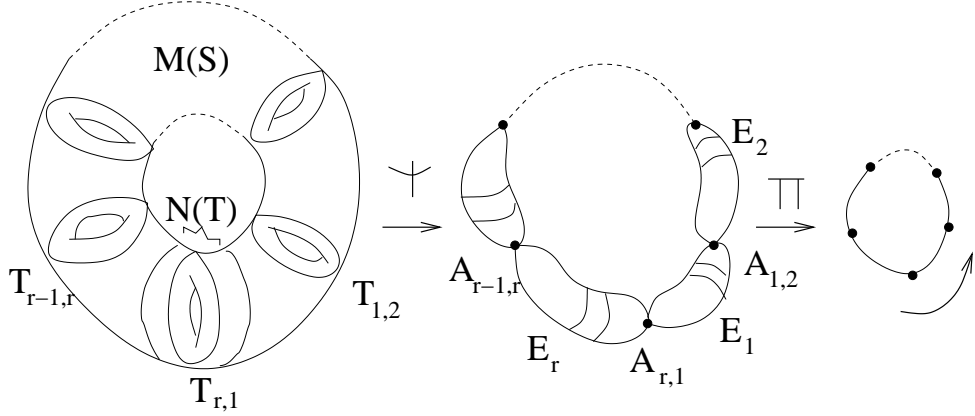, height= 55 mm} 
  \caption{The maps $\Pi$ and $\Psi$ for a cusp singularity}
\end{figure}

Our construction shows that $\Pi$ realizes $M(\mathcal{S})$ as the
total space of a torus fibration over $\R/\Z$.

Denote by $T_{i, i+1}:= \Psi^{-1}(A_{i, i+1})$ the torus of
$\mathcal{T}(p)$ which corresponds to the intersection point of $E_i$
and $E_{i+1}$. Denote $T:= T_{r,1}$ and let $N(T)$ be a (closed)
tubular neighborhood of $T$, which does not intersect any other torus
$T_{i, i+1}$, for $i \in \{1,..., r-1\}$ (see Figure 20). 

Denote $L:= H_1(M(\mathcal{S})-N(T), \Z)$. As $M(\mathcal{S})- N(T)$
is the interior of a thick torus, we see that $L$ is a 2-dimensional
lattice. With the notations of section  \ref{plumcros}, let $f_i$ be
an oriented 
fibre in the piece $M(E_i)$. We suppose moreover that $f_1$ and $f_r$
are situated on the boundary of $N(T)$. Consider two other circles
$f_0$ and $f_{r+1}$ on $\partial N(T)$, such that $f_0, f_r$ are
isotopic \textit{inside $N(T)$} and situated on distinct boundary
components and such that the same is true for the pair $f_1,
f_{r+1}$. 

For each $i \in \{0,..., r+1\}$, denote by $v_i:= [f_i]\in L$ the
homology class of $f_i$. By Lemma \ref{reltor}, we see that: 
\begin{equation} \label{recfibr}
 v_{i+1}=
\mid e_i\mid \cdot v_i - v_{i-1}= -e_i \cdot v_i -v_{i-1}, \: \forall
\: i \in \{0,...,r\},
\end{equation}
where $E_0:=E_r$. 

Denote by $n \in GL(L)$ the automorphism which sends the basis $(v_0,
v_1)$ of $L$ into the basis $(v_r, v_{r+1})$. The relations
(\ref{recfibr}) 
show that its matrix in the basis $(v_0, v_1)$ is:
  $$\left( \begin{array}{cc}
              0 & -1\\
              1 & e_1
           \end{array} \right)
\left( \begin{array}{cc}
              0 & -1\\
              1 & e_2
           \end{array} \right)\cdots
\left( \begin{array}{cc}
              0 & -1\\
              1 & e_r
           \end{array} \right)$$
A little thinking shows that $n$ is the inverse of the algebraic
monodromy $m \in GL(L)$ in the positive direction along $\R/\Z$. So,
the matrix of $m$ in the basis $(v_0, v_1)$ is:
  $$\left( \begin{array}{cc}
              e_r & 1\\
              -1 & 0
           \end{array} \right)
\left( \begin{array}{cc}
              e_{r-1} & 1\\
              -1 & 0
           \end{array} \right)\cdots
\left( \begin{array}{cc}
              e_1 & 1\\
              -1 & 0
           \end{array} \right)$$
We have reproved like this Theorem 6.1 IV in Neumann \cite{N
  81}. We deduce by induction the following
  expression for its trace, where the polynomials $Z^-$ were
  defined by formula (\ref{recrel}):
\begin{equation} \label{formtrace}
 tr \: m = Z^- (| e_1|,..., |e_r|)- Z^- (|e_2|,..., |e_{r-1}|).
\end{equation}

The negative definiteness of the intersection
matrix of $E$ (see Theorem \ref{critcontr}) shows that there exists $i
\in \{1,...,r\}$ such that $|e_i|\geq 3$. As $p$ is supposed to be
the minimal resolution of $(\mathcal{S},0)$, we have also $e_j \geq 2,
\: \forall \: j \in \{1,...,r\}$. Using equation (\ref{formtrace}), we
deduce then easily by induction on $r$ that $tr \: m \geq 3$. 
\medskip

$\bullet$ {\em Consider now the case $r=1$}. Then, by Lemma
\ref{whynot}, $E$ is a 
rational curve with one 
singular point $P$, where $E$ has normal crossings. Let
$p':(\mathcal{R}', E')\rightarrow (\mathcal{S},0)$ be the resolution
of $(\mathcal{S},0)$ obtained by blowing up $P \in \mathcal{R}$. Then
$E'$ is a normal crossings resolution with smooth components $E_1,
E_2$, where $E_1^2= -1$ and $E_2$ is the strict transform of $E$. As
$(p')^*E= 2E_1 + E_2$ and $((p')^*E)^2=E^2$, we deduce that $E_2^2
=E^2-4 \leq -5$. Now we apply the same argument as in the case $r\geq
2$, but for the resolution $p'$. 

An alternative proof could use Lemma \ref{compeul}.

The fact that each oriented torus fibration with $tr \: m\geq 3$
appears like this is a consequence of the study done in section
\ref{torfibr}. Indeed, there we show how to extract the numbers
$(e_1,...,e_r)$ from the oriented topological type of
$M(\mathcal{S})$. 
\hfill $\Box$
\medskip

By Neumann \cite{N 81}, there exist also abstract boundaries
$M(\mathcal{S})$ which are 
torus fibrations with algebraic monodromy of trace $2$. But in that
case the exceptional divisor of the minimal resolution is an elliptic
curve (then, following Saito \cite{S 74}, one speaks about
\textit{simple elliptic singularities}, which are other particular
cases of minimally elliptic ones).  
\medskip

\subsection{Construction of the canonical graph structure}
\label{construct} $\:$
\medskip

Consider again an arbitrary normal surface singularity
$(\mathcal{S},0)$ and a normal crossings resolution $p$ of it. 

\begin{definition} 
\textit{Suppose  that the set of nodes $\mathcal{N}(p)$ is
  non-empty}. Conceive the graph 
$\Gamma(p)$ as a 1-dimensional CW-complex and take the complement
$\Gamma(p)- \mathcal{N}(p)$. This complement is the disjoint union
of segments, which we call \textbf{chains}. If a chain is open at both
extremities we call it an \textbf{interior chain}. If it is half-open
we call it a \textbf{terminal  chain}. 
\end{definition}

In Figure 20 we represent the chains of Figure 17, with the
hypothesis that $E_4, E_5, E_7 \notin \mathcal{N}(p)$ and $E_6 \in
\mathcal{N}(p)$. That is, we suppose that $E_4, E_5, E_6, E_7$ are
smooth and that $g(E_4)=g(E_5)=g(E_7)=0, \:
g(E_6)\geq 1$. There is only one terminal chain, which contains the 
terminal chain vertex $E_7$.

Denote by $\mathcal{C}(p)$ the set of chains. This set can
be written as a disjoint union
  $$ \mathcal{C}(p)= \mathcal{C}_i(p)\sqcup  \mathcal{C}_t(p)$$
where $\mathcal{C}_i(p)$ denotes the set of interior chains and
$\mathcal{C}_t(p)$ the set of terminal chains. The edges of $\Gamma(p)$
contained in a chain $C\in \mathcal{C}(p)$ correspond to a set of
parallel tori in $M(\mathcal{S})$. Choose one torus $T_C$ among them
and define:
  $$\mathcal{T}'(p):= \bigsqcup _{C \in \mathcal{C}_i(p)}T_C.$$

{\tt    \setlength{\unitlength}{0.92pt}}
\begin{figure}
  \epsfig{file=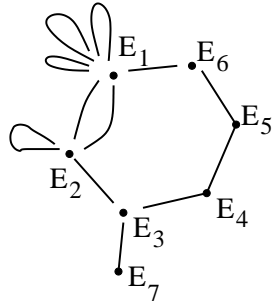, height= 40 mm} 
  \caption{The chains of Figure 17 when $E_6$ is a node}
\end{figure}

By construction, each piece of $M(\mathcal{S})_{\mathcal{T}'(p)}$
contains a unique piece $M(E_i)$ of
$M(\mathcal{S})_{\mathcal{T}(p)}$ such that $E_i$ is a node of
$\Gamma(p)$. If $E_i\in \mathcal{N}(p)$, denote by $M'(E_i)$ the piece
of $M(\mathcal{S})_{\mathcal{T}'(p)}$ which contains $M(E_i)$. One
can extend in a unique way up to isotopy the natural Seifert structure
without exceptional fibres on 
$M(E_i)$ to a Seifert structure on $M'(E_i)$. One obtains like this a
graph structure $(\mathcal{T}'(p), \mathcal{F}'(p))$ on
$M(\mathcal{S})$. 

Till now we have worked with any normal crossings resolution $p$. We
consider now a special one, the minimal normal crossings resolution
$p_{mnc}$.

\begin{proposition} \label{canid}
  Suppose that $(\mathcal{S},0)$ is neither a Hirzebruch-Jung
  singularity, nor a cusp singularity. Then 
  the graph structure $(\mathcal{T}'(p_{mnc}), \mathcal{F}'(p_{mnc}))$ is the
  canonical graph structure on $M(\mathcal{S})$.
\end{proposition}

\textbf{Proof:} If $\mathcal{T}'(p_{mnc})$ \textit{is empty}, as
$(\mathcal{S},0)$ is not a cusp singularity we deduce that
$(\mathcal{T}'(p_{mnc}), \mathcal{F}'(p_{mnc}))$ is a Seifert
structure. By Proposition \ref{nonun}, we see that it is the canonical
graph structure on $M(\mathcal{S})$.

Suppose now that
$\mathcal{T}'(p_{mnc})$ \textit{is non-empty}. 
One has to verify two facts (see Definition \ref{cangraph}):
  
$\bullet$ first, that all the fibrations induced by
$\mathcal{F}'(p_{mnc})$ on the pieces which are thick Klein bottles
have orientable basis;

$\bullet$ second, that by taking the various choices of Seifert
structures on the pieces of $M(\mathcal{S})_{\mathcal{T}'(p_{mnc})}$,
  one does not obtain isotopic fibres coming from different sides on
  one of the tori of $\mathcal{T}'(p_{mnc})$.

The first fact is immediate, as one starts from Seifert structures with
orientable basis on the pieces of
$M(\mathcal{S})_{\mathcal{T}(p_{mnc})}$ before eliminating tori of
$\mathcal{T}(p_{mnc})$ in order to remain with $\mathcal{T}'(p_{mnc})$. 

In what concerns the second fact, the idea is to look at the fibres
corresponding to the chain vertices of any interior chain $C$. The union
of the pieces of $M(\mathcal{S})_{\mathcal{T}(p_{mnc})}$ which are
associated to those vertices is a thick torus $N_R$. Take a fibre in
each piece (remember that they are naturally oriented as boundaries of
holomorphic discs) and look at their images in $L=H_1(N_R, \Z)$. One gets like
this a sequence of vectors $v_1,...,v_s \in L$. Consider also the
images $v_0$ and $v_{s+1}$ of the fibres coming from the nodes 
of $\Gamma(p_{mnc})$ to which $C$ is adjacent, the order of the
indices respecting the order of the vertices along the chain. 

By Lemma \ref{reltor}, $v_{k+1}= \alpha_k v_k -v_{k-1}$ for any
$k \in \{1,...,s\}$, where $\alpha_k$ is the absolute value of the
self-intersection of the component $E_i$ of $\Gamma(p_{mnc})$ which
gave rise to the vector $v_k$. Here plays
the hypothesis that 
$p_{mnc}$ is \textit{minimal}: this implies that $\alpha_k \geq
2$. Then one can conclude by using Proposition \ref{seifbound}. 

The analysis of thick Klein bottles is similar. It is based on the
fact that a thick Klein bottle can appear only from a portion of the
graph $\Gamma(p)$ as in Figure 21, where $E_1,E_2, E_3$ are smooth
rational curves of self-intersections $-2, -2$, respectively $-n$ (see
Neumann \cite[pages 305, 334]{N 81}). The
important point is that $n \geq 2$. Otherwise the complete
sub-graph of $\Gamma(p)$ with vertices $E_1, E_2, E_3$ would have a
non-definite intersection matrix, which contradicts Theorem
\ref{critcontr}.\hfill $\Box$ 

{\tt    \setlength{\unitlength}{0.92pt}}
\begin{figure}
  \epsfig{file=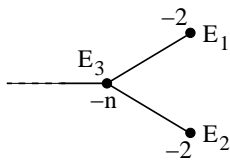, height= 20 mm} 
  \caption{The appearance of a thick Klein bottle}
\end{figure}
\medskip

The plumbing structure $(\mathcal{T}(p_{mnc}), \mathcal{F}(p_{mnc}))$
on $N(\mathcal{S})$ is associated to the resolution $p_{mnc}$ of
$(\mathcal{S},0)$. One can wonder if the canonical graph structure
$(\mathcal{T}'(p_{mnc}),$ $\mathcal{F}'(p_{mnc}))$ is also associated to
some analytic morphism with target $(\mathcal{S},0)$.

This is indeed the case. In order to see it, start from $p_{mnc}$ and
its exceptional divisor $E$. Then contract all the components of $E$
which correspond to chain vertices. One gets like this a normal surface
with only Hirzebruch-Jung singularities. The image of $E$ on it is a
divisor $F$ with again only normal crossings when seen as an abstract
curve. Take then as a representative of $M(\mathcal{S})$ the boundary
of a tubular neighborhood of $F$ in the new surface and split it into
pieces which project into the various components of $F$. The splitting
is done using tori which are associated bijectively to the singular
points of $F$. Namely, in a system of (toric) local coordinates
$(x,y)$ such that $F$ is defined by $xy=0$, one proceeds as for the
definition of the plumbing structure associated to a normal crossings
resolution (see Section \ref{plumcros}). Then this system of tori is
isotopic to $\mathcal{T}'(p_{mnc})$.

\medskip
\section{Invariance of the canonical plumbing  structure on the
  boundary of a normal surface singularity} \label{canplumb}

In this section  we describe how to reconstruct the plumbing
structure \linebreak $(\mathcal{T}(p_{mnc}), \mathcal{F}(p_{mnc}))$ on
$M(\mathcal{S})$ associated to the minimal normal crossings 
resolution of $(\mathcal{S},0)$, only from the abstract oriented manifold
$M(\mathcal{S})$. Namely,  using the classes of
plumbing structures on thick tori defined in section
\ref{HJplumb}, we define a plumbing structure
$\mathcal{P}(M(\mathcal{S}))$ on $M(\mathcal{S})$ and we prove: 

\begin{theorem} \label{canlens}
  1) When considered as an unoriented structure, the plumbing
     structure $\mathcal{P}(M(\mathcal{S}))$ depends up to
     isotopy only on
     the natural orientation of $M(\mathcal{S})$. We call it the
     \textbf{canonical plumbing structure} on 
     $M(\mathcal{S})$.   

  2) The plumbing
  structure $(\mathcal{T}(p_{mnc}), \mathcal{F}(p_{mnc}))$
  associated to the minimal normal crossings resolution of
  $(\mathcal{S},0)$ is isotopic to the canonical
  plumbing structure $\mathcal{P}(M(\mathcal{S}))$.
\end{theorem}

As a corollary we get the theorem of invariance of the plumbing
structure 
$(\mathcal{T}(p_{mnc}), \mathcal{F}(p_{mnc}))$ announced in the
introduction (see Theorem \ref{invplumb}). We also explain how the
orientation reversal on the boundary of a Hirzebruch-Jung or cusp
singularity reflects the duality between supplementary cones explained
in section \ref{comp1} (see Propositions \ref{isolens} and
\ref{dualcusp}). 
\medskip

 In order to prove Theorem \ref{canlens}, we  consider three
 cases, according to the 
 nature of $M(\mathcal{S})$. In the first one it is supposed to be a
 lens space, in the second one a torus fibration with algebraic
 monodromy of trace $\geq 3$ and in the last one
 none of the two (so, by Proposition \ref{HJcusp}, this corresponds to
 the trichotomy: $(\mathcal{S},0)$ is a Hirzebruch-Jung singularity/ a
 cusp singularity/ none of the two). 
 
The idea is to start from some structure on
 $M(\mathcal{S})$ which is well-defined up to isotopy, and to enrich
 it by 
 canonical constructions of Hirzebruch-Jung plumbing structures
 (defined in section \ref{HJplumb}). When
 $M(\mathcal{S})$ is 
 neither a lens space nor a torus fibration with algebraic
 monodromy of trace $\geq 3$, this starting structure
 will be the canonical graph structure (see Definition
 \ref{cangraph}). Otherwise we need some special  theorems of
 structure (Theorems \ref{Bonahon} and \ref{Wald}).

\subsection{The case of lens spaces} \label{caselens}             $\:$
\medskip

Notice that by Proposition \ref{HJcusp} 1), $M(\mathcal{S})$ is a
lens space if and only if 
$(\mathcal{S},0)$ is a Hirzebruch-Jung singularity. 

The following theorem was proved by Bonahon \cite{B 83}:

\begin{theorem} \label{Bonahon}
  Up to isotopy, a lens space contains a unique torus which splits it
  into two solid tori.
\end{theorem}

We say that a torus embedded in a lens space and splitting it
into two solid tori is a \textit{central torus}. By the previous
theorem, a central torus is well-defined up to isotopy. 

Let $M$ be an \textit{oriented} lens space and $T$ a central torus in
$M$. Consider a tubular neighborhood $N(T)$ of $T$ in $M$, whose
boundary components we denote by $T_-$ and $T_+$, ordered in an
arbitrary way. Then $M_{T_- \sqcup T_+}$ has three pieces, one being
sent diffeomorphically by the reconstruction map $r_{M, T_- \sqcup
  T_+}$ on $N(T)$ - by a slight abuse of notations, we keep calling it $N(T)$ -
 and the others, $M_-$ and $M_+$, having boundaries sent by $r_{M, T_- \sqcup
  T_+}$ on $T_-$, respectively $T_+$ (see Figure 33). The manifolds 
$M_-$ and $M_+$ are solid tori, as $T$ was supposed to be a central
torus. Let $\gamma_-$ and $\gamma_+$ be meridians of $M_-$,
respectively $M_+$, oriented compatibly with the orientation of $N(T)$
(see Definition \ref{compat1}). Consider the Hirzebruch-Jung plumbing structure
$\mathcal{P}(N(T), \gamma_-, \gamma_+)$ on $N(T)$, whose tori are
denoted by $T_0=T_-, T_1,..., T_{r+2}=T_+$, as explained in Section
\ref{HJplumb}.

{\tt    \setlength{\unitlength}{0.92pt}}
\begin{figure}
  \epsfig{file=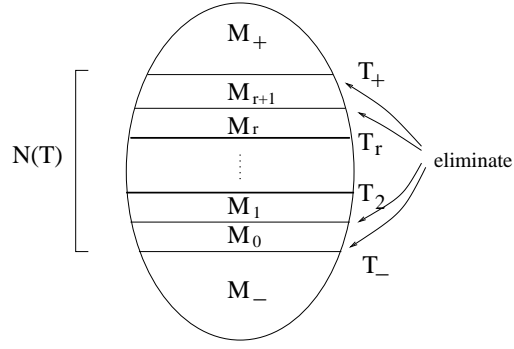, height= 45 mm} 
  \caption{Construction of the canonical plumbing structure on a lens space}
\end{figure}

Denote $\mathcal{T}_M:= T_2 \sqcup \cdots \sqcup T_{r}$. Then
$M_{\mathcal{T}_M}$ contains four pieces less than the manifold $M_{T_- \sqcup
  \mathcal{T}_M \sqcup T_+}$. Denote by $M_-'$ and $M_+'$ the piece
which ``contains" $M_-$, respectively $M_+$. On $M_-'$ we consider the
Seifert structure  which
extends the Seifert structure of $M_1$ and on $M_+'$ the one which
extends the Seifert structure of $M_{r}$. By applying the intersection
theoretical criterion of Proposition \ref{seifbound} b), we see that
those Seifert structures have no exceptional fibres (we used a similar
argument to construct in Section \ref{HJplumb} the Hirzebruch-Jung
plumbing structure on solid tori). 
On the other pieces of
$M_{\mathcal{T}_M}$ we consider the Seifert structure coming from the
plumbing structure $\mathcal{P}(M, \gamma_-, \gamma_+)$. Denote by
$\mathcal{P}(M)$ the plumbing structure constructed like
this on the oriented manifold $M$. 
\medskip

\textbf{Proof of Theorem \ref{canlens}:} 
 
1) This is obvious by construction (we use Theorem \ref{Bonahon}).

2) In the construction of $\mathcal{P}(M(\mathcal{S}))$, one can take
   as central torus $T$ any torus \linebreak $(\Pi\circ \Psi)^{-1}(c)$, with $c
   \in (0,1)$, in the notations of the proof of Proposition
   \ref{HJcusp}, 1). Then one sees that $[\gamma_-]=[f_0]$ and
   $[\gamma_+]= [f_{r+1}]$ in the lattice $L= H_1(M(\mathcal{S})-(\Pi\circ
   \Psi)^{-1}\{0,1\}, \Z)=H_1(T, \Z)$. Using the relations
   (\ref{relens}) and the 
   definition of a Hirzebruch-Jung plumbing structure on a thick torus
   (see section \ref{HJplumb}), we deduce that the images of the
   fibres $f_i$ in $L$ are equal to the images of the fibres of
   $\mathcal{P}(M(\mathcal{S}))$ (see also Proposition
   \ref{converse}). The proposition follows by the fact that on a
   2-torus, any oriented essential curve is well-defined up to isotopy
   by its homology class. 
\hfill $\Box$
\medskip 

Let $\sigma$ be the strictly convex cone of $L_{\R}$ whose edges are
generated by $[\gamma_-]$ and $[\gamma_+]$. If one changes the
ordering of the components of $\partial N(T)$, then one gets the same
cone $\sigma$, and if one changes simultaneously the orientations of
$\gamma_-$ and $\gamma_+$, then one gets the opposite cone. But if one
\textit{changes the orientation of} $M$, then the cone $\sigma$ is
replaced by a supplementary cone.  
So, in view of Section \ref{reldual}, the two cones are in
duality. In this sense, the canonical plumbing structure
$\mathcal{P}(-M(\mathcal{S}))$ is dual to
$\mathcal{P}(M(\mathcal{S}))$. We get:

\begin{proposition} \label{orientchange}
   Let $(\mathcal{S},0)$ be a Hirzebruch-Jung singularity. Then the
   canonical plumbing structures with respect to the two possible 
   orientations of $M(\mathcal{S})$ are dual to each other. More
   precisely, if $(\mathcal{S},0) \simeq (\mathcal{Z}(L,
   \sigma),0)\simeq \mathcal{A}_{p,q}$, then $-M(\mathcal{S})$ is
   orientation preserving diffeomorphic to $M(\check{\mathcal{S}})$,
   where, with the notations of section \ref{geomHJ}, 
 $(\check{\mathcal{S}},0)\simeq (\mathcal{Z}(\check{L},
   \check{\sigma}),0)\simeq \mathcal{A}_{p,p-q}$.
\end{proposition}

Let $\lambda:=\dfrac{p}{q}$ be the type of the cone $(L, \sigma)$ in
the sense of Definition \ref{typecone}, where $0<q<p$ and
$gcd(p,q)=1$. The oriented lens space $M(\mathcal{S})$, where
$(\mathcal{S},0) \simeq (\mathcal{Z}(L,  \sigma),0)\simeq
\mathcal{A}_{p,q}$, is said classically to be \textit{of type}
$L(p,q)$. By Propositions \ref{changeorder} and \ref{dualtype},
combined with Theorem \ref{Bonahon}, we get the following classical fact:

\begin{proposition} \label{isolens}
  1) The lens spaces $L(p,q)$ and $L(p,q')$ are orientation-preserving
     diffeomorphic if and only if $p=p'$ and $q'\in \{ q,
     \overline{q}\}$, where $0<\overline{q}<p, \: q\overline{q}\equiv
     1 (mod \: p)$.

  2) The lens spaces $L(p,q)$ and $L(p,q')$ are orientation-reversing
     diffeomorphic if and only if $p=p'$ and $q'\in \{ p-q, p-
     \overline{q}\}$.
\end{proposition}

\medskip

\subsection{The case of torus fibrations with $tr \: m \geq 3$}
\label{torfibr}$\:$ 
\medskip

Notice that by Proposition \ref{HJcusp} 2), $M(\mathcal{S})$ is a
torus fibration whose algebraic monodromy verifies $tr \: m \geq 3$ if
and only if $(\mathcal{S},0)$ is a cusp singularity. First we 
study with a little more detail torus fibrations. 

Let $M$ be an orientable torus fibration. Take a fibre torus $T$. Then
consider the lattice  
$L=H_1(T, \Z)$ and the algebraic monodromy operator $m \in SL(L)$ (see
Definition \ref{algmon}) associated with one of the two possible
orientations of the base. 

The following theorem is a consequence of Waldhausen
\cite[section 3]{W 68} (see also Hatcher \cite[section 5]{H 00}): 
 
\begin{theorem} \label{Wald}
  Up to isotopy, an orientable torus fibration  $M$ such that $tr \: m \geq 3$ 
  contains a unique torus which splits it  into a thick torus (see
  Definition \ref{split}).
\end{theorem}

We say that a torus embedded in an orientable torus fibration
whose algebraic monodromy $m$ verifies $tr \: m \geq 3$ 
and which splits it into a thick torus is
a \textit{fibre torus}. By the previous theorem, a fibre torus is
well-defined up to isotopy. 

 From now on, we
suppose that indeed  $tr \: m \geq 3$ (see Proposition \ref{HJcusp},
2)).  As $M$ is orientable, $m$ preserves the orientation of $L$,
which shows that $det \:  
m=1$. This implies that  the characteristic polynomial of $m$ is $X^2
- (tr \: m) X 
+1$. We deduce that $m$ has two strictly positive eigenvalues with product
$1$, and so the eigenspaces are two distinct real lines in $L_{\R}$.

But the most important point is that these lines are
\textit{irrational}. Indeed, the eigenvalues are $\nu_{\pm}
:=\frac{1}{2}(tr\: m \pm \sqrt{(tr\: m)^2 -4})$ and $(tr\: m)^2 -4$ is
never a square if $tr\: m \geq 3$. 

Denote by $d_-$ and $d_+$ the eigenspaces corresponding to $\nu_-$,
respectively $\nu_+$. Then $m$ is strictly contracting when restricted
to $d_-$ and strictly expanding when restricted to $d_+$. Choose 
arbitrarily one of the two half-lines in which $0$ divides the line
$d_-$, and call it $l_-$. 

At this point we have not used any
orientation of $M$. \textit{Suppose now 
that $M$ is oriented}. Then the chosen orientation
on the basis of the torus fibration induces an orientation of the
fibre 
torus $T$, by deciding that this orientation, followed by the
transversal orientation which projects on the orientation of the base
induces the ambient orientation on $M$. 

Denote by $l_+$ the half-line
bounded by $0$ 
on $d_+$ into which $l_-$ arrives first when turned in the
\textit{negative} 
direction. Let $\sigma$ be the strictly convex cone bounded by these
two half-lines (see Figure 24).

{\tt    \setlength{\unitlength}{0.92pt}}
\begin{figure}
  \epsfig{file=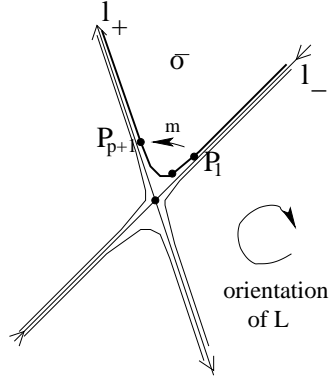, height= 50 mm} 
  \caption{The case of torus bundles with $tr \: m \geq 3$}
\end{figure}

We arrive like this at a pair $(L, \sigma)$ where both edges of
$\sigma$ are irrational. As $m$ preserves $L$ and $\sigma$, \textit{it
preserves also the polygonal line $P(\sigma)$}. 

Let $P_1$ be an
arbitrary integral point of $P(\sigma)$. Consider the sequence $(P_n)_{n \geq
  1}$ of integral points of $P(\sigma)$ read in the positive direction
along $P(\sigma)$, starting from $P_1$. There exists an index $t \geq
1$ such that $P_{t+1} = m(P_1)$. It is the \textit{period} of the
action of $m$ on the linearly ordered set of integral points of
$P(\sigma)$. 

Consider  $t$ parallel tori $T_1,..., T_t$ inside $M$, where
$T_1=T$ and the indices form an increasing function of the orders of
appearance of the tori when one turns in the positive
direction. Denote $\mathcal{T}:= \bigsqcup_{1 \leq k \leq t}T_k$ and
$T_{t+1}:= T_1$. For each $k \in \{1,...,t\}$, denote by $M_k$ the
piece of $M_{\mathcal{T}}$ whose boundary components project by $r_{M,
  \mathcal{T}}$ on $T_k$ and $T_{k+1}$ (see Figure 25). Then look at
the thick torus 
$M_T$. Let $T_-$ be its boundary component through which one ``enters
inside" $M_T$ when one turns in the positive direction, and $T_+$ be
the one by which one ``leaves" $M_T$. Identify then $H_1(M_T, \Z)$ with
$H_1(T_-, \Z)$ through the inclusion $T_- \subset M_{T_1}$ and
$H_1(T_-, \Z)$ with $H_1(T, \Z)=L$ through the reconstruction mapping
$r_{M,T}|_{T_-}: T_- \rightarrow T$.

{\tt    \setlength{\unitlength}{0.92pt}}
\begin{figure}
  \epsfig{file=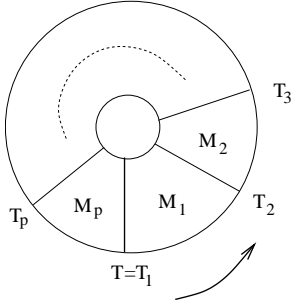, height= 40 mm} 
  \caption{Construction of the canonical plumbing structure on a torus
  fibration with $tr \: m \geq 3$}
\end{figure}

Consider now on each piece $M_k$ an oriented Seifert fibration
$\mathcal{F}_k$ such that the class of a fibre in $L$ (after
projection in $M_T$ and identification of $H_1(M_T, \Z)$ with $L$, as
explained before) is equal to $OP_k$. Denote by $\mathcal{F}$ the
Seifert structure on $M_{\mathcal{T}}$ obtained by taking the union of
  the structures $\mathcal{F}_k$. We get like this a plumbing
  structure on $M$. Denote it by $\mathcal{P}(M)$. 

This plumbing structure does not depend, up to isotopy, on the choice of
the initial 
integral point on $P(\sigma)$. Indeed, by lifting to $M$ a vector
field of the form $\frac{\partial}{\partial \theta}$ on the base of
the torus fibration and by considering its flow, one sees that one
gets isotopic torus fibrations by starting from any integral point of
$P(\sigma)$.

Notice that it does neither depend on the choice of the half-line
$l_-$. An opposite choice would lead to the choice of an opposite
cone, that is to the same unoriented plumbing structure. 
\medskip

\textbf{Proof of Theorem \ref{canlens}:} 

1) This is obvious by construction (we use Theorem \ref{Wald}).

2) In the construction of $\mathcal{P}(M(\mathcal{S}))$, one can take
   as fibre torus $T$ the torus $T_{r,1}$, with the notations of the
   proof of Proposition \ref{HJcusp}, 2). Using the relations
   (\ref{recfibr}) and Proposition \ref{converse}, we get the
   Proposition. 
\hfill $\Box$
\medskip

By Theorem \ref{taut}, cusp singularities are determined up to
analytic isomorphism by the topological type of the oriented manifold
$M(\mathcal{S})$. By Theorem \ref{Wald}, this manifold can be encoded
by a pair $(T,\mu)$, where $T$ is an oriented fibre and $\mu$ is a 
geometric monodromy diffeomorphism of $T$ obtained by turning in the
positive direction 
determined by the chosen orientation of $T$ (recall that this is
precisely the point were we use the given orientation of
$M(\mathcal{S})$). But it is known that $\mu$ can be reconstructed up
to isotopy by its action on $L=H_1(T, \mathbf{Z})$, that is, by the
algebraic monodromy operator $m \in SL(L)$. Moreover, to fix an
orientation of $T$ is the same as to fix an orientation of $L$. As
explained in section \ref{reldual}, such an orientation can be encoded
in a symplectic isomorphism $\omega:L \rightarrow \check{L}$.

Denote by $\mathcal{C}(L, \omega, m)$ the cusp singularity associated
to an oriented lattice $(L, \omega)$ and an algebraic monodromy
operator $m \in SL(L)$ with $tr \: m \geq 3$. If one changes the
orientation of the base of the torus fibration, one gets the triple
$(L, -\omega, m^{-1})$. This shows that:
$$\mathcal{C}(L, \omega, m)\simeq \mathcal{C}(L, -\omega, m^{-1}).$$

When one changes the orientation of $M(\mathcal{S})$, we see that the
  cone $(L, \sigma)$ is replaced by a supplemetary one. In view of
  Section \ref{reldual}, we deduce that the two cones are dual to each
  other. In this sense, 
  we get the following analog of Proposition \ref{orientchange}:

\begin{proposition} \label{dualcusp}
  Let $(\mathcal{S},0)$ be a cusp singularity. Then the
   canonical plumbing structures with respect to the two possible 
   orientations of $M(\mathcal{S})$ are dual to each other. More
   precisely, if $(\mathcal{S},0) \simeq (\mathcal{C}(L, \omega,
   m),0)$, then $-M(\mathcal{S})$ is 
   orientation preserving diffeomorphic to $M(\check{\mathcal{S}})$,
   where $(\check{\mathcal{S}},0)\simeq (\mathcal{C}(L, -\omega,
   m),0)$.
\end{proposition}

\medskip

\subsection{The other singularity boundaries} $\:$
\medskip

As in the two previous cases, we first define the plumbing structure
$\mathcal{P}(M(\mathcal{S}))$. 

Consider the canonical graph structure $(\mathcal{T}_{can},
\mathcal{F}_{can})$ on $M(\mathcal{S})$. We  do our construction
starting from the neighborhoods of the JSJ tori (the elements of
$\mathcal{T}_{can}$) and the exceptional fibres in
$\mathcal{F}_{can}$.

$\bullet$ For each component $T$ of $\mathcal{T}_{can}$, consider a saturated
   tubular neighborhood $N(T)$. We choose them pairwise disjoint. So,
   each manifold $N(T)$ is a thick torus. 
We consider on each one of its
boundary components a fibre of $\mathcal{F}_{can}$. Denote these
   fibres by $\gamma(T), \delta(T)$. We consider on  $N(T)$ the
   restriction  of the orientation of $M(\mathcal{S})$. 
Consider the associated Hirzebruch-Jung plumbing
structure $\mathcal{P}(N(T), \gamma(T), \delta(T))$ (see
Definition \ref{HJPS1}). Replace the Seifert structure on $N(T)$
induced from $\mathcal{F}_{can}$ with this plumbing structure. Then
eliminate the boundary components of $N(T)$ from the tori present in
$M(\mathcal{S})$ (by construction, the fibrations coming from both
   sides agree on them up to isotopy). 

$\bullet$  For each exceptional fibre $F$, 
consider a solid torus $N(F)$, which is a saturated tubular neighborhood
of $F$. Choose those neighborhoods pairwise disjoint.  On the 
boundary of $N(F)$, take a fiber $\gamma(F)$ of
$\mathcal{F}_{can}$. Consider the associated Hirzebruch-Jung plumbing
structure $\mathcal{P}(N(F), \gamma(F))$ (see Definition
\ref{HJPS2}). Replace the Seifert structure on $N(F)$ induced from
$\mathcal{F}_{can}$ with this plumbing structure. Then eliminate the
boundary component of $N(F)$ from the tori present inside
$M(\mathcal{S})$ (by construction, the fibrations coming from both
sides agree on it up to isotopy). Denote by
$\mathcal{P}(M(\mathcal{S}))$ the plumbing structure constructed like
this on $M(\mathcal{S})$. 

\medskip
\textbf{Proof of Theorem \ref{canlens}:}  The proof is very
similar to the ones explained in the two previous cases, but starting this time
from the canonical graph structure on $M(\mathcal{S})$. The main point
is Proposition \ref{canid}. We leave the details to the reader. 
  \hfill $\Box$
\medskip

\subsection{The invariance theorem} $\:$
\medskip

Let $(\mathcal{S},0)$ be a normal surface singularity. In \cite{N 81},
Neumann proved that the weighted dual graph $\Gamma(p_{mnc})$ of the
exceptional divisor of its minimal normal crossings resolution
$p_{mnc}$ is determined by the oriented manifold $M(\mathcal{S})$. But
he says nothing about the action of the group
$\operatorname{Diff}^+(M(\mathcal{S}))$ on $(\mathcal{T}(p_{mnc}),
  \mathcal{F}(p_{mnc}))$. As a corollary of Theorem \ref{canlens} we
  get:

\begin{theorem} \label{invplumb}
  The plumbing structure $(\mathcal{T}(p_{mnc}),
  \mathcal{F}(p_{mnc}))$  is invariant up to isotopy by the group
  $\operatorname{Diff}^+(M(\mathcal{S}))$. 
\end{theorem}

\textbf{Proof: } Suppose first that $M(\mathcal{S})$ is not a lens
space or a torus fibration. As the canonical graph structure on it is
invariant by the group $\operatorname{Diff}^+(M(\mathcal{S}))$ up to
isotopy, we deduce that 
the canonical plumbing structure is also invariant up to isotopy by
this group. This conclusion is also true when $M(\mathcal{S})$ is a
lens space or a torus fibration, as one starts in the construction of
$\mathcal{P}(M(\mathcal{S}))$ from tori which are invariant up to
isotopy. Then we apply Theorem 
\ref{canlens}. 
  \hfill $\Box$

An easy study of the fibres of $\mathcal{F}(p_{mnc})$ in the
neighborhoods of the tori of $\mathcal{T}(p_{mnc})$ which correspond
to self-intersection points of components of $E_{mnc}$ show that the 
analogous statement about the minimal \textit{good} normal crossings
resolution of $\mathcal{S}$ is also true. 
\medskip

We arrived at the conclusion that the affirmation of 
Theorem \ref{invplumb} was true while we were thinking about the
natural contact structure on 
$M(\mathcal{S})$ (see Caubel, N{\'e}methi \& Popescu-Pampu \cite{CNP
  05}). Indeed, in that paper we prove that for normal surface singularities,
the natural contact structure depends only on the topology of
$M(\mathcal{S})$ \textit{up to contactomorphisms}. It was then natural
to look at the subgroup of $\operatorname{Diff}^+(M(\mathcal{S}))$
which leaves it 
invariant up to isotopy. Presently, we do not know how to characterize
it. But we realized that the homotopy type of the underlying
unoriented plane
field was invariant by the full group $\operatorname{Diff}^+(M(\mathcal{S}))$,
provided that Theorem \ref{invplumb} was true (see \cite[section 5]{CNP
  05}).

\par\medskip\centerline{\rule{2cm}{0.2mm}}\medskip
\setcounter{section}{0}
\medskip


{\small
\begin{thebibliography}{99}
%
\bibitem{A 98} Arnold, V.I. \textit{Higher dimensional continued
    fractions.} Regular and chaotic dynamics \textbf{3}, 3 (1998),
    10-17. 

\bibitem{A 62} Artin, M. \textit{Some numerical criteria for
    contractability of curves on algebraic surfaces.} Am. J. 
    of Maths. \textbf{84} (1962), 485-496.

\bibitem{A 66} Artin, M. \textit{On isolated rational singularities of
    surfaces.} Am. J.  
    of Maths. \textbf{88} (1966), 129-136.

\bibitem{BHPV 04} Barth, W.P., Hulek, K., Peters, C.A.M., Van de Ven,
  A. \textit{Compact complex surfaces.} Second enlarged edition,
  Springer, 2004. 

\bibitem{B 83} Bonahon, F. \textit{Diff{\'e}otopies des espaces
    lenticulaires.} Topology \textbf{22} (1983), 305-314. 

\bibitem{B 02} Bonahon, F. \textit{Geometric structures on
    3-manifolds.} In \textit{Handbook of geometric topology.}
    R.J. Daverman, R.B.Sher eds., Elsevier, 2002, 93-164. 

\bibitem{B 91} Brezinski, C. \textit{History of Continued Fractions
    and Pad{\'e} Approximants.} Vol. \textbf{12} of Springer Series in
    Computational Mathematics, Springer, 1991. 

\bibitem{B 00} Brieskorn, E. \textit{Singularities in the work of
    Friedrich Hirzebruch.} Surveys in Differential Geometry
    \textbf{VII} (2000), 17-60. 

\bibitem{BK 86} Brieskorn, E., Kn{\"o}rrer, H. \textit{Plane algebraic
    curves.} Birkh{\"a}user, 1986. 



\bibitem{CA 00} Casas-Alvero, E. \textit{Singularities of plane
    curves.} London Mathematical Society Lecture Note Series
    \textbf{276}, Cambridge U.P., 2000.

\bibitem{CNP 05} Caubel, C., N{\'e}methi, A., Popescu-Pampu,
  P. \textit{Milnor open books and Milnor fillable contact
  3-manifolds.} To appear in Topology.

\bibitem{C 73} Cohn, H. \textit{Support polygons and the resolution of
    modular functional singularities.} Acta Arithmetica \textbf{XXIV} (1973),
    261-278. 

\bibitem{C 00} Cossart, V. \textit{Uniformisation et
    d{\'e}singularisation des surfaces d'apr{\`e}s Zariski.} Progress in
    Math. \textbf{181}, 2000, 239-258. 

\bibitem{DHH 98} Dais, D.I., Haus, U.U., Henk, M. \textit{On crepant
    resolutions of 2-parameter series of Gorenstein cyclic quotient
    singularities.} Results Math. \textbf{33} (1998) no. 3-4,
    208-265. 

\bibitem{D 99} Davenport, H. \textit{The Higher Arithmetic.} Seventh
  edition,  Cambridge Univ. Press, 1999. 

\bibitem{D 83} Durfee, A.H. \textit{Neighborhoods of algebraic sets.}
  Trans. Am. Math.Soc. \textbf{276} (1983), 517-530. 

\bibitem{EW 85} Ebeling, W., Wall, C.T.C. \textit{Kodaira
    singularities and an extension of Arnold's strange duality.}
    Compositio Mathematica \textbf{56} (1985), 3-77.

\bibitem{EN 85} Eisenbud, D., Neumann, W. \textit{Three-Dimensional Link
    Theory and Invariants of Plane Curve Singularities.} Annals of Math.
    Studies \textbf{110}, Princeton Univ.Press, 1985.

\bibitem{EC 18} Enriques, F., Chisini, O. \textit{Lezioni sulla teoria
    geometrica delle equazioni e delle funzioni algebriche.} vol. II,
    N. Zanichelli, Bologna, 1918. 

\bibitem{E 72} Epstein, D.P.A. \textit{Periodic flows on 3-manifolds.}
  Ann. of Math. \textbf{95} (1972), 66-82. 

\bibitem{FM 97} Fomenko, A.T., Matveev, S.V. \textit{Algorithmic and Computer 
 Methods for Three-Manifolds.} Mathematics and its Applications 
 vol. \textbf{425}, 
 Kluwer Acad. Publishers, 1997. Translated from the russian edition,
 Moscow Univ. Press, 1991. 

\bibitem{F 87} Fowler, D.H. \textit{The mathematics of Plato's
    Academy. A new reconstruction.} Clarendon Press, 1987.

\bibitem{F 93} Fulton, W. \textit{Introduction to Toric Varieties.} 
  Princeton Univ. Press, 1993.

\bibitem{GB 96} Garc{\'\i}a Barroso, E. R. \textit{Invariants des
    singularit{\'e}s de courbes planes et courbure des fibres de
    Milnor.} Tesis doctoral, Univ. de La Laguna, Tenerife (Spain),
    1996. 

\bibitem{G 00} Giroux, E. \textit{Structures de contact en dimension
trois et bifurcations des feuilletages de surfaces,} Inv. Math.
\textbf{141} (2000), 615-689.

\bibitem{GT 00} Goldin, R., Teissier, B. \textit{Resolving plane
    branch singularities with one toric morphism.} In
    \textit{Resolution of singularities, a research textbook in
    tribute to Oscar Zariski.} Birkh{\"a}user, Progress in
    Math. \textbf{181}, 2000, 315-340. 

\bibitem{GP 03} Gonz{\'a}lez P{\'e}rez, P.D. \textit{Toric embedded
    resolutions of quasi-ordinary hypersurface singularities.}
    Ann. Inst. Fourier, Grenoble \textbf{53}, 6 (2003), 1819-1881.

\bibitem{GPGS 04} Gonz{\'a}lez P{\'e}rez, P.D., Gonzalez-Sprinberg,
  G. \textit{Analytical invariants of quasi-ordinary hypersurface
  singularities associated to divisorial valuations.} Kodai
  Math. J. \textbf{27}, 2 (2004), 164-173.

\bibitem{GPT 02} Gonz{\'a}lez P{\'e}rez, P.D., Teissier, B. \textit{Toric embedded
    resolution of non-necessarily normal toric varieties.}
    C. R. Acad. Sci. Paris \textbf{334} (2002) 5, 379-382. 

\bibitem{GS 77} Gonzalez-Sprinberg, G. \textit{Eventails en dimension
    2 et transform{\'e}e de Nash.} Preprint, ENS Paris, 1977. 

\bibitem{G 62} Grauert, H. \textit{{\"U}ber Modifikationen und exzeptionnelle
  analytische Mengen.} Math. Ann. \textbf{146} (1962), 331-368. 


\bibitem{HW 88} Hardy, G.H., Wright, E.M. \textit{An introduction to
    the theory of numbers.} V-th edition, Clarendon Press, Oxford,
    1988. 

\bibitem{H 00} Hatcher, A. \textit{Notes on Basic 3-manifold
  Topology.} Available on the web-page  http://math.cornell.edu/~hatcher.

\bibitem{H 53} Hirzebruch, F. \textit{{\"U}ber vierdimensionale 
    Riemannsche Fl{\"a}chen Mehrdeutiger analytischer Funktionen von
    zwei komplexen Ver{\"a}nderlichen.} Math. Ann. \textbf{126}, 1-22
    (1953). 

\bibitem{H 62} Hirzebruch, F. \textit{The Topology of Normal
  Singularities of an Algebraic Surface (d'apr{\`e}s un article de
  D.Mumford).} Sem. Bourbaki 1962/63, Exp. 250; Gesammelte Abhandlungen,
  Band II, Springer-Verlag, 1987, 1-7.

\bibitem{H 73} Hirzebruch, F. \textit{Hilbert modular surfaces.}
  Enseign. Math. \textbf{19} (1973), 183-281. 

\bibitem{HNK 71} Hirzebruch, F., Neumann, W., Koh, S.S.
  \textit{Differentiable Manifolds and Quadratic Forms.} Marcel
  Dekker, Inc. 1971.

\bibitem{J 80} Jaco, W.H. \textit{Lectures on Three-Manifold Topology.}
  Regional Conference Series in Mathematics no.\textbf{43}, AMS,
  1980. 

\bibitem{JS 79} Jaco, W.H., Shalen, P.B. \textit{Seifert Fibered
  Spaces in Three-manifolds.} Memoirs of the AMS vol.\textbf{21}, 220,
  1979.



\bibitem{J 79} Johannson, K. \textit{Homotopy Equivalences of
    3-Manifolds with Boundaries.} LNM \textbf{761}, Springer-Verlag, 1979.


\bibitem{J 08} Jung, H.W.E. \textit{Darstellung der Funktionen eines 
algebraischen K{\"o}rpers zweier unabh{\"a}ngigen Ver{\"a}nderlichen x,y in
  der Umgebung einer Stelle x=a, y=b.} J. Reine
  Angew. Math. \textbf{133} (1908), 289-314. 

\bibitem{J 85} Jurkiewicz, J. \textit{Torus embeddings, polyhedra,
    k$^*$-actions and homology.} Dissertationes Mathematicae (Rozprawy
    Matematyczne) \textbf{236} (1985).

\bibitem{KKMD 73} Kempf, G., Knudsen, F., Mumford, D., St. Donat,
  B. \textit{Toroidal embeddings.} Springer LNM \textbf{339}, 1973. 


\bibitem{K 96} Klein, F. \textit{{\"U}ber eine geometrische Auffassung
    der gew{\"o}hlischen Kettenbuchentwicklung.}
    Nachr. Ges. Wiss. G{\"o}ttingen. Math.-Phys. Kl. \textbf{3} (1895),
    357-359. French translation: \textit{Sur une repr{\'e}sentation
    g{\'e}om{\'e}trique du d{\'e}veloppement en fraction continue
    ordinaire.} Nouvelles Annales de Math{\'e}matiques (3) \textbf{15}
    (1896), 327-331.

\bibitem{K 04} Klein, F. \textit{Elementary mathematics from an
    advanced standpoint. Arithmetic, algebra, analysis.}
    Dover, 2004. Translation of the third german edition of 1925. 

\bibitem{L 71} Laufer, H.B. \textit{Normal two-dimensional
    singularities.} Princeton Univ. Press, 1971. 

\bibitem{L 73} Laufer, H.B. \textit{Taut two-dimensional
    singularities.} Math. Ann. \textbf{205}, 131-164 (1973). 

\bibitem{L 77} Laufer, H.B. \textit{On minimally elliptic
    singularities.} Am. J. of Math. \textbf{99} (6), 1257-1295 (1977).

\bibitem{L 00} L{\^e}, D.T. \textit{Geometry of complex surface
    singularities.} In \textit{Singularities-Sapporo 1998.} Advanced
    Studies in Pure Maths. \textbf{29} (2000), 163-180.

\bibitem{LMW 89} L{\^e}, D.T., Michel, F., Weber, C. \textit{Sur le
    comportement des polaires 
                 associ{\'e}es aux germes de courbes planes.} Compositio 
                 Mathematica \textbf{72}, (1989), 87-113.

\bibitem{L 75} Lipman, J. \textit{Introduction to resolution of
    singularities.} Proceedings of Symposia in Pure Mathematics
    \textbf{29} (1975), 187- 230. 

\bibitem{MW 85} Michel, F., Weber, C. \textit{Topologie des germes de courbes
    planes {\`a} plusieurs branches.} Preprint Univ. Gen{\`e}ve,
    1985.

\bibitem{M 00} Moussafir, J.-O. \textit{Voiles et poly{\`e}dres de
    Klein. G{\'e}om{\'e}trie, algorithmes et statistiques.} Th{\`e}se,
    Univ. Paris IX-Dauphine, 2000. Available on the page:
    http://www.ceremade.dauphine.fr/ $\tilde{}$ msfr/ 

\bibitem{M 61} Mumford, D. \textit{The Topology of Normal
    Singularities of an Algebraic Surface and a Criterion for
    Simplicity.} Publ. Math. IHES \textbf{9} (1961), 229-246.

\bibitem{N 04} N{\'e}methi, A. \textit{Invariants of normal surface
    singularities.} Contemporary Mathematics \textbf{354}, 161-208
    AMS (2004). 

\bibitem{N 81} Neumann, W. \textit{A calculus for plumbing applied to
    the topology of complex surface singularities and degenerating
    complex curves.} Transactions of the AMS \textbf{268}, 2 (1981),
    299-344. 

\bibitem{NR 78} Neumann, W.D., Raymond, F. \textit{Seifert Manifolds,
  Plumbing, $\mu$-Invariant and Orientation Reversing Maps.} Algebraic
  and Geometric Topology Proceedings, Santa Barbara (1977), LNM
  \textbf{664}, Springer-Verlag, 1978, 162-195.

\bibitem{NS 97} Neumann, W.D., Swarup, G.A. \textit{Canonical
    Decompositions of 3-Manifolds.} Geometry and Topology,
    \textbf{1} (1997), 21-40.

\bibitem{O 88} Oda, T. \textit{Convex bodies and algebraic geometry, 
  an introduction to the theory of toric varieties.} Ergebnisse der
  Math. \textbf{15}, Springer-Verlag, 1988.

\bibitem{O 72} Orlik, P. \textit{Seifert manifolds.} Lecture Notes in
  Mathematics, Vol. \textbf{291}. Springer-Verlag, Berlin-New York, 1972.  


\bibitem{PP 01} Popescu-Pampu, P. \textit{Arbres de contact des 
  singularit{\'e}s quasi-ordinaires et graphes d'adjacence pour les 
  3-vari{\'e}t{\'e}s r{\'e}elles.} Th{\`e}se, Univ. Paris 7, 2001. Available on the
  page 
http://tel.ccsd.cnrs.fr/documents/archives0/00/00/28/00/index\_fr.html 

\bibitem{PP 02} Popescu-Pampu, P. \textit{On a canonical placement of
    knots in irreducible 3-manifolds.} C.R. Acad. Sci. Paris, S{\'e}rie
    I \textbf{334} (2002), 677-682. 

\bibitem{PP 05} Popescu-Pampu, P. \textit{On higher dimensional
    Hirzebruch-Jung singularities.} Rev. Mat. Complut. \textbf{18}
    (2005), no.1, 209-232. 

\bibitem{R 62} Randow, R. von \textit{Zur Topologie von
  dreidimensionalen Baummannigfaltigkeiten.} Bonner Mathematischer Schriften
  \textbf{14}, 1962.

\bibitem{R 74} Riemenschneider, O. \textit{Deformationen von
    Quotientensingularit{\"a}ten (nach Zyklischen Gruppen).}
    Math. Ann. \textbf{209}, 211-248 (1974). 

\bibitem{S 74} Saito, K. \textit{Einfach-elliptische Singularit{\"a}ten.}
  Invent. Math. \textbf{23} (1974), 289-325.

\bibitem{S 75} Scharf, A. \textit{Zur Faserung von
    Graphenmannigfaltigkeiten.} Math. Ann. \textbf{215}, 35-45
    (1975). 

\bibitem{S 83} Scott, P. \textit{The geometries of 3-manifolds.}
  Bull. London Math. Soc. \textbf{15} (1983), 401-487. 

\bibitem{S 33} Seifert, H. \textit{Topologie dreidimensionaler
    gefaserter R{\"a}ume.} Acta Math. \textbf{60} (1933), 147-238. English
    translation in \cite{ST 80}.

\bibitem{ST 80} Seifert, H., Threlfall, W. \textit{A Textbook of
    Topology.} Academic Press, 1980. English translation of \textit{Lehrbuch
    der Topologie.}  Teubner, 1934.

\bibitem{S 92} Shallit, J. \textit{Real numbers with bounded partial
    quotients.} L'Enseignement Math{\'e}matique \textbf{38} (1992),
    151-187. 


\bibitem{S 76} Smith, H.J.S. \textit{Note on continued fractions.}
  Messenger of Mathematics, Ser. II, vol. vi (1876), 1-14. Reedited in
  Collected Papers II, 135-147.  

\bibitem{S 90} Spivakovsky, M. \textit{Valuations in function fields
    of surfaces.} Am. J. of Maths. \textbf{112} (1990),
    107-156. 

\bibitem{S 78} Stark, H.M. \textit{An introduction to number theory.}
  The MIT Press, 1978.

\bibitem{S 96} Sturmfels, B. \textit{Gr{\"o}bner bases and convex
    polytopes.} Univ. Lecture Series \textbf{8},  AMS, 1996.

\bibitem{T 04} Teissier, B. \textit{Monomial ideals, binomial ideals,
    polynomial ideals.} In \textit{Trends in commutative algebra.}
    L.L. Avramov, M.Green, C. Huneke, K.E. Smith, B. Sturmfels eds.,
    Cambridge U.P., 2004, 211-246.  

\bibitem{T 59} Thom, R. \textit{Les structures diff{\'e}rentiables des
    boules et des sph{\`e}res.} Centre Belge Rech. Math., Colloque
    G{\'e}om. Diff{\'e}r. Globale, Bruxelles, du 
    19 au 22 D{\'e}c. 1958, 27-35 (1959). 

\bibitem{T 68} Tjurina, G.N. \textit{Absolute isolatedness of rational
    singularities and triple rational points.} Functional Analysis and
    its Applications \textbf{2}, 4 (1968), 324-332. 

\bibitem{V 44} Du Val, P. \textit{On absolute and non-absolute
    singularities of algebraic surfaces.} Revue de la Facult{\'e} des
    Sciences de l'Univ. d'Istanbul (A) \textbf{91} (1944), 159-215.

\bibitem{V 89} Vogt, E. \textit{A Foliation of $\mathbf{R}^3$ and
  other Punctured 3-Manifolds by Circles.} Publ. Math. IHES
  \textbf{69} (1989), 215-232.

\bibitem{W 67} Waldhausen, F. \textit{Eine Klasse von 3-dimensionalen 
               Mannigfaltigkeiten.} Inv. Math. \textbf{3} (1967),
               308-333 and     Inv. Math. \textbf{4} (1967), 87-117.

\bibitem{W 68} Waldhausen, F. \textit{On Irreducible 3-Manifolds which
  are Sufficiently Large}, Annals of Maths. \textbf{87} (1968), 56-88.


\bibitem{W 00} Wall, C.T.C. \textit{Quadratic forms and normal surface
    singularities.} In \textit{Quadratic forms and their
    applications. (Dublin 1999)}, Contemp. Math. \textbf{272},
    293-311, AMS (2000).

\bibitem{W 04} Wall, C.T.C. \textit{Singular points of plane curves.}
  Cambridge UP, 2004. 

\bibitem{W 39} Whitehead, J.H.C. \textit{Simplicial spaces, nuclei,
    and $m$-groups.} Proc. London Math. Soc. \textbf{45} (1939),
    243-327. 

\bibitem{Z 39} Zariski, O. \textit{The reduction of the singularities
    of an algebraic surface.} Ann. of Math. \textbf{40} (1939), 639-689.

\bibitem{Z 99} Zisman, M. \textit{Fibre bundles, fibre maps.} In
  \textit{History of Topology.} I.M. James ed., Elsevier, 1999,
  605-629. 

\end{thebibliography} }
\end{document}